\documentclass[12pt]{amsart}
\usepackage{mathtools, amsmath, xcolor, marginnote,graphics,latexsym}
\usepackage{natbib}
\usepackage{amssymb, esint, upgreek}
\usepackage{amsthm}
\usepackage{hyperref}
\usepackage{cleveref}
\usepackage{float}
\usepackage{bm, enumerate}
\graphicspath{{./images/}}
\newtheorem{definition}{\bf Definition}[section]
\newtheorem{proposition}{\bf Proposition}[section]
\newtheorem{theorem}{\bf Theorem}[section]
\newtheorem{lemma}{\bf Lemma}[section]
\newtheorem{corollary}{\bf Corollary}[section]

\newtheorem{result}{\bf Result}[section]
\numberwithin{equation}{section}
\begin{document}
\title[Higher Integrability of Differential Form]{Higher Integrability of Differential Form Globally Solving Certain Minimization Problem}
\author[Adhikari]{Sayan Adhikari}
    \address{Indian Institute of Science, Department of Mathematics, Bengaluru, 560012, India}
    \email{sayana@iisc.ac.in, sayanadhikari7998@gmail.com}
	\date{\today}
\begin{abstract}
Let $0\leq k\leq n-1,1<p<\infty$ and $u\in W^{1,p}\left(B;\Lambda^k\right)$, where $\Lambda^k$ is the vector space of all $k$-forms in $\mathbb{R}^n$ and $B$ is an open ball in $\mathbb{R}^n$. In this article, with one global higher integrability assumption on $du$, we study the global higher integrability of $dv$ where $v$ is the unique minimizer of the following functional:
\begin{align*}
\overline{w}\mapsto \int_{B} \left \lvert d\overline{w} \right \rvert^p \text{ where } \overline{w} \in u+W^{1,p}_{d^*,T}\left(B;\Lambda^{k}\right).
\end{align*}
As a machinery, with a similar higher integrability assumption on the $k+1$-differential form $\widetilde{F}$, we prove one local higher integrability of $d\widetilde{W}$ where $\widetilde{W}\in W^{1,p}\left(B_r\left(0\right);\Lambda^{k}\right)$ weakly solves:
\begin{align*}
d^*\left(\left \langle G \left(d\widetilde{W}+\widetilde{F}\right), d\widetilde{W}+\widetilde{F}\right  \rangle^{\frac{p-2}{2}}G \left(d \widetilde{W}+\widetilde{F}\right)\right)=0 \text{ in } B_r\left(0\right),
\end{align*}
for $G:B_r\left(0\right)\rightarrow GL\left({n\choose k+1};\mathbb{R}\right)$ symmetric, bounded, and uniformly elliptic  measurable map, which is sufficiently close to the identity map $I:\Lambda^{k+1}\rightarrow \Lambda^{k+1}$.
For each of the above mentioned cases, we prove one strong $\left(q,q\right)$ type estimate. As the primary tools, we use the Hardy-Littlewood maximal inequality, Fefferman-Stein inequality, and one classical regularity result of K. Ulhenbeck. 
\end{abstract}
\maketitle
\section{Introduction}
Fix an integer $n\geq 2$, a real number $p$ in $\left(1,\infty\right)$ and 
\begin{align}\label{qtilde}
\widetilde{q}\in
\begin{cases} \left[p,\infty\right)  &\text{ for } n=2,\\ 
\left[p,\frac{np}{n-2}\right] &\text{ for } n\geq 3.
\end{cases}
\end{align}
 For $k\in \left\{0,1,\cdots,n-1\right\}$, we assume $G:B_{r}\left(0\right)\rightarrow GL\left({n\choose k+1},\mathbb{R}\right)$ is a measurable, symmetric, bounded function with $\lambda_G,\Lambda_G>0$ such that
\begin{align}\label{ellipticity}
\lambda_G\left \lvert \xi\right \rvert^2\leq \left \langle G\left(x\right) \xi,\xi \right \rangle\leq \Lambda_G\left \lvert \xi\right \rvert^2 \text{ for almost all $x\in B_r\left(0\right)$}.    
\end{align}
For $j\in \left\{0,1,\cdots,n-1\right\}$, let $\Lambda^{j}$ be the space of all $j$-forms on $\mathbb{R}^n$. Thus, for $x\in B_r\left(0\right)$ we understand $G\left(x\right)$ as an element of $\mathcal{L}\left(\Lambda^{k+1};\Lambda^{k+1}\right)$ which is the space of all linear map from the $n\choose k+1$-dimensional space to itself.\\ In the first part of this article we study the weak solution $ \widetilde{W}\in W^{1,p}\left(B_r\left(0\right);\Lambda^k\right)$ of the following:
    \begin{align}\label{First mention of G eqn}
        d^*\left(\left \langle G \left(d\widetilde{W}+\widetilde{F}\right), d\widetilde{W}+\widetilde{F}\right  \rangle^{\frac{p-2}{2}}G \left(d \widetilde{W}+\widetilde{F}\right)\right)=0 \text{ in } B_r\left(0\right)
    \end{align}
  under the assumption $\widetilde{F}\in L^{\widetilde{q}}\left(B_r\left(0\right);\Lambda^{k+1}\right)$.
In this article we use Hardy-Littlewood Maximal Inequality and Fefferman-Stein Inequality; we briefly recall the results:\\
     \begin{result}\textbf{Hardy-Littlewood Maximal Inequality:}  Let, $f$ be locally integrable in $\mathbb{R}^n$, then the maximal function of $f$ say $Mf$ is defined by
        \begin{align*}
            Mf\left(x\right):=\sup_{r>0}\frac{1}{\left \lvert B_r\left(x\right)\right \rvert}\int_{B_r\left(x\right)}\left \lvert f
            \right \rvert.
        \end{align*}
        If $1<\gamma\leq \infty$ and $f\in L^\gamma\left(\mathbb{R}^n\right)$, then there exists $A_0\left(n,\gamma\right)>0$  such that
    \begin{align}\label{Mf <f}
        \int_{\mathbb{R}^n}\left\lvert Mf\right\rvert^\gamma\leq A_0\left(n,\gamma\right)\int_{\mathbb{R}^n}\left\lvert f\right\rvert^\gamma \text{ for all $f\in L^\gamma\left(\mathbb{R}^n\right)$ }.
    \end{align}
    (see \cite[Theorem $1 \left(c\right)$]{8bbc344b-39d3-31c4-89fe-4b604b1a5849}).
    \end{result}
    \begin{result}
    \textbf{Fefferman-Stein Inequality:}
        Suppose that, $Q_0$  is a cube in $\mathbb{R}^n$ with sides parallel to the co-ordinate axes and $g\in  L^1\left(Q_0\right)$, then we define the (centered) sharp maximal function of $g$, say $g^{\#}$ by
        \begin{align*}
            g^{\#}\left(x\right):=\sup_{Q\left(x\right)\subseteq Q_0} \frac{1}{\left \lvert Q\left(x\right)\right \rvert}\int_{Q\left(x\right)} \left \lvert g-\left(g\right)_{Q\left(x\right)}\right \rvert.
        \end{align*}
        where $\left(g\right)_{Q\left(x\right)}:=\frac{1}{\left \lvert Q\left(x\right)\right\rvert}\int_{Q\left(x\right)}g$, and the supremum is taken over all cubes $Q\left(x\right)\left(\subseteq Q_0\right)$ centered at $x$.
        If $1<{\gamma}<\infty$ and $g^{\#}\in L^\gamma\left(Q_0\right)$, then $g\in L^{\gamma}\left(Q_0\right)$ and
        \begin{align}\label{Fefferman Stein Cube}
            \int_{Q_0}\left \lvert g \right \rvert^{\gamma}\leq A\left(n,{\gamma}\right)\left[\int_{Q_0}\left \lvert g^{\#} \right  \rvert^{\gamma}+\left \lvert Q_0\right \rvert^{1-\gamma} \left(\int_{Q_0}\left \lvert g \right \rvert\right)^{\gamma}\right],
        \end{align}
(see \cite[Theorem $6.30$]{giaquinta2013introduction}; 
the result stated in there, is for uncentered sharp maximal function but since centered and uncentered sharp maximal functions have similar growth
\eqref{Fefferman Stein Cube} follows.)
\end{result}
Before we state the classical regularity result of K. Ulhenbeck we recall some standard auxiliary functions and related inequalities, which will be used throughout the discussion:\\
For all $z \in \Lambda^{k+1}$ and $1<\gamma<\infty$ we define    
\begin{align}\label{Definition of Vgamma}
V_\gamma\left(z\right):=\left \lvert z\right \rvert^{\frac{\gamma-2}{2}}z.
\end{align} 
Notice that $\left \lvert V_\gamma\left(z\right)\right \rvert^2=\left \lvert z\right \rvert^\gamma$ for $\gamma \in \Lambda^{k+1}$. For $\xi, \zeta \in \Lambda^{k+1}$ we have the following inequalities:
\begin{align}\label{V sim}
    \left(\left \lvert \xi\right \rvert+\left \lvert \zeta\right \rvert\right)^{\frac{\gamma-2}{2}}\left \lvert \xi-\zeta\right\rvert \sim \left \lvert V_\gamma\left(\xi\right)-V_\gamma\left(\zeta\right)\right \rvert,
\end{align}
(see \cite[Lemma $2.1$]{MR1179331} with $\delta=\frac{\gamma-2}{2}, \mu=0$) and
\begin{align}\label{V upsim}
    \left\langle \left \lvert \xi\right \rvert^{\gamma-2}\xi-\left \lvert \zeta\right\rvert^{\gamma-2}\zeta,\xi-\zeta\right \rangle\geq c_0\left(p\right)\left \lvert V_\gamma\left(\xi\right)-V_\gamma\left(\zeta\right)\right \rvert^2,
          \end{align}
          (see $\left(2.5\right)$ from \cite{MR3294408}).\\
         In this article we are mainly concerned with the higher integrability of $d\widetilde{W}$ and as one of the primary tools we mention the classical estimates of K.Ulhenbeck.\\
         \begin{result}
         
\textbf{K. Ulhenbeck's Theorem:} Let, $\Omega$ be an open subset of $\mathbb{R}^n, 1<p<\infty, 0\leq k\leq n-1$ and let $u \in W^{d,1}_{\operatorname{loc}}\left(\Omega; \Lambda^k\right)$ be a local  minimizer of the functional
        $\overline{w}\mapsto \int_{\Omega}\left \lvert d\overline{w}\right \rvert^p$, then $du$ is locally H\"{o}lder continuous. Moreover,
        \begin{enumerate} 
            \item There exists $\bar{c}$ depending on $n,k,p$ such that for all ball $B_\rho \subset \Omega$ with radius $\rho$ we have the following: 
            \begin{align}\label{sup bound}
                \sup_{B_{\frac{\rho}{2}}}\left\lvert du\right \rvert^p \leq \frac{\bar{c}}{\left \lvert B_\rho\right \rvert} \int_{B_\rho}\left \lvert du\right \rvert^p.
            \end{align}
            In the above estimate $B_{\frac{\rho}{2}}$ and $B_\rho$ have the same center.
            \item There are constants $c_U\equiv c_U\left(n,k,p\right)>0 $ and $\alpha_1\equiv \alpha_1\left(n,k,p\right)\in \left(0,1\right)$ such that for all balls $B_\rho \subset\subset B_R \left(\subset \Omega\right)$ with the same center we have the following decay estimate:
            \begin{align}\label{decay estimate}
                \frac{1}{\left \lvert B_\rho\right \rvert}\int_{B_\rho}& \left \lvert V_p\left(du\right)-\left(V_p\left(du\right)\right)_{B_\rho} \right \rvert^2 \leq c_U\left(\frac{\rho}{R}\right)^{2\alpha_1}\frac{1}{\left \lvert B_R\right \rvert}\int_{B_R} \left \lvert V_p\left(du\right)-\left(V_p\left(du\right)\right)_{B_R} \right \rvert^2,
            \end{align}
        \end{enumerate}
         (see \cite[Theorem~$4.1$]{MR1179331}, \cite[Theorem $8$]{MR3989958} and K. Ulhenbeck's \cite{MR474389}).\\
        \end{result}
          We define
          \begin{align}\label{c*0}
            c^*_0\left(n\right):=\frac{\left \lvert B\right \rvert}{\left \lvert Q\right \rvert},  
          \end{align}
          where $B$ is any ball in $\mathbb{R}^n$ and $Q$ is any (maximal) cube inscribed inside $B$.\\
          
          In \cite[Theorem~$3$]{bwmeta1.element.bwnjournal-article-smv75i1p22bwm}, T. Iwaniec proved the local higher integrability of $\nabla u$; where $p\geq 2$ and $u$ is a weak solution of\begin{align*}
              \operatorname{div}\left(\left \langle G\left(x\right)\nabla u,\nabla u\right \rangle^{\frac{p-2}{2}}G\left(x\right)\nabla u\right)=0.
          \end{align*} 
          He showed that, if $p\leq q<\infty$ and $K_G-1$ is small enough  depending on $q$
          then $\nabla u \in L^q_{\operatorname{loc}}\left(\Omega\right)$, where the domain of $G$ is $\Omega$ and $K_G$ is defined by
          \begin{align*}K_G:=\left(1+ \left \lVert G-I\right \rVert_{L^\infty\left(\Omega\right)}\right)^{\frac{p}{2}}.
          \end{align*} 
  In our article to denote a similar smallness of $\left \lVert G-I\right \rVert_{B_r\left(0\right)}$, we define the $*$-condition on $G
  $ by the following three conditions: \\
     We fix $\delta_*>0$ small enough so that
    \begin{align}\label{delta}
    \delta_*=\min\left \{\frac{1}{16}, \left[1000 \times 2^{\tilde{q}}A_0\left(n,\frac{\widetilde{q}}{p}\right)A\left(n,\frac{2\widetilde{q}}{p}\right)\left[8c_0^*\left(n\right)c_U\right]^{\frac{\tilde{q}}{p}} \right]^{-\frac{p}{2\alpha_1 \widetilde{q}}}\right\}.
    \end{align} 
    Depending on $\delta_*>0$,  $\left \lVert G-I\right \rVert_{L^\infty\left(B_r\left(0\right)\right)}$ is small enough such that
    \begin{align}\label{r}
  \left[\frac{2^{p+4}}{c_0\left(p\right)} c_0^*\left(n\right) \delta_*^{-n} K_p\left \lVert G-I\right \rVert_{L^\infty\left(B_r\left(0\right)\right)}\right]^\frac{\tilde{q}}{p} A_0\left(n,\frac{\widetilde{q}}{p}\right)A\left(n,\frac{2\widetilde{q}}{p}\right) \leq  \frac{10^{-3}}{2^{\tilde{q}}},
  \end{align}
  and 
  \begin{align}\label{1<p<2 and p<2lambda first mention}
        &\frac{1}{2}\leq \lambda_G\leq \Lambda_G\leq \frac{3}{2}\text{ and }\left\lvert 1-\Lambda_G^{\frac{p-2}{2}}\lambda_G\right \rvert,\left \lvert 1-\lambda_G^{\frac{p-4}{2}}\Lambda_G^2\right \rvert<\frac{p-1}{6-2p}.
        \end{align}
    Here, $\widetilde{q}, A_0,A, c_0\left(p\right)$ and $c_0^*\left(n\right) $ are the constants mentioned in \eqref{qtilde}\eqref{Mf <f},\eqref{Fefferman Stein Cube},\eqref{V upsim} and \eqref{c*0} respectively. The constants $c_U,\alpha_1$ are from \eqref{decay estimate} of Ulhenbeck's theorem.\\
    
    In our first Theorem we prove the local higher integrability of $d\widetilde{W}$ and one $L^q-L^q$ estimate:
            
            \begin{theorem}\label{G Theorem}
                Let $1<p<\infty, 0\leq k\leq n-1$ and let $G:B_{r}\left(0\right)\rightarrow GL\left({n\choose k+1},\mathbb{R}\right)$ be measurable, symmetric, bounded map that satisfies \eqref{ellipticity} and $*$-condition (see \eqref{delta}, \eqref{r} and \eqref{1<p<2 and p<2lambda first mention}). If $\widetilde{F}\in L^{\widetilde{q}}\left(B_r\left(0\right);\Lambda^{k+1}\right)$ where $\widetilde{q}$ is any real number satisfying \eqref{qtilde} and $\widetilde{W}\in W^{1,p}\left(B_r\left(0\right);\Lambda^k\right)$  is a weak solution of the following: 
\begin{align}\label{G weak formulation in theorem}
        d^*\left(\left \langle G \left(d\widetilde{W}+\widetilde{F}\right), d\widetilde{W}+\widetilde{F}\right  \rangle^{\frac{p-2}{2}}G \left(d \widetilde{W}+\widetilde{F}\right)\right)=0 \text{ in } B_r\left(0\right);
    \end{align}
    then there exists a cube $Q_0\left(\subset\subset B_r\left(0\right)\right)$ with sides parallel to the coordinate axes such that
    \begin{align*}
        \int_{Q_0}\left \lvert d\widetilde{W}+\widetilde{F} \right \rvert^{\widetilde{q}}  &\leq \begin{aligned}[t]
        &c\left(n,\tilde{q},p,k,r\right) \left[\left[\int_{Q_0}\left \lvert d\widetilde{W}+\widetilde{F} \right \rvert^p\right]^{\frac{\tilde{q}}{p}}
        +\int_{Q_0}\left \lvert \widetilde{F}\right\rvert^{\tilde{q}}\right].
        \end{aligned}
    \end{align*}
    Here $\left \lvert Q_0\right \rvert$ can be chosen depending on $n,p,\widetilde{q},k$ and $r$ only.
    \end{theorem}
    In the second part of this paper we deal with the following minimization problem:
    \begin{align}\label{Lq Lq Minimization Problem}
            \overline{w}\mapsto \int_{B} \left \lvert d\overline{w} \right \rvert^p \text{ where } \overline{w} \in u+W^{1,p}_{d^*,T}\left(B;\Lambda^{k}\right),
        \end{align}
        where $B$ is an open ball in $\mathbb{R}^n$. If $v\in u+W^{1,p}\left(B;\Lambda^k\right)$ is the unique minimizer of the functional \eqref{Lq Lq Minimization Problem} then by definition we get
        \begin{align}\label{vp<up}
            \int_B \left \lvert dv\right \rvert^p\leq \int_B \left \lvert du\right \rvert^p.
        \end{align}
        Using K. Ulhenbeck's sup bound estimate in \eqref{sup bound} we have the local higher integrability of $dv$. 
    In the last part of our discussion, as a corollary of Theorem \ref{G Theorem} we prove a global $L^q-L^q$ type estimate corresponding to 
    $dv$ which is the main result of this article:
        \begin{theorem}\label{Goal Theorem}
            Let $B$ be an open ball in $\mathbb{R}^n, n\geq 2, 1<p<\infty,0\leq k\leq n-1$ and $ u\in W^{1,p}\left(B;\Lambda^k\right)$ with $du\in L^{\widetilde{q}}\left(B;\Lambda^{k+1}\right)$ where $\widetilde{q}$ satisfies \eqref{qtilde}. If $v\in u+W^{1,p}_{d^*,T}\left(B;\Lambda^k\right)$ is the unique minimizer of the functional mentioned in \eqref{Lq Lq Minimization Problem}, then
             \begin{align}\label{goal}
            \int_B\left \lvert dv\right \rvert^{\widetilde{q}}\leq c\left(n,p,\widetilde{q},k\right)\int_{B} \left \lvert du\right \rvert^{\widetilde{q}}.
        \end{align} 
        \end{theorem}
\smallskip

\section{Notations and Preliminaries}
\subsection{Spaces of Differential Forms}
We start this section with a brief introduction to the differential form in $\mathbb{R}^n$ and mention some crucial results that will be used throughout this paper; (see \cite[Chapter $3$]{csato2011pullback}).\\
Let $n\geq 2, 0\leq k\leq n$ be integers, and let  $\Lambda^{k}$ be the space of all $k$-forms on $\mathbb{R}^n$. For a measurable subset $\Omega$ in $\mathbb{R}^n$, we define a differential $k$-form $f :\Omega \rightarrow \Lambda^k$ by
\begin{align*}
f:=\sum_{1\leq i_1<\cdots<i_k\leq n}f_{i_1\cdots i_k}dx^{i_1}\wedge \cdots \wedge dx^{i_k}= \sum_{\substack{I \in \mathcal{J}^k}}f_Idx^I,
\end{align*}
where $\mathcal{J}^k:=\left\{\left(i_1,i_2,\cdots, i_k\right):1\leq i_1<i_2<\cdots<i_k\leq n\right\}$ and $f_I=f_{i_1\cdots i_k}: \Omega \rightarrow \mathbb{R}$ is measurable for each $I\in \mathcal{J}^k$. For $f=\sum_{\substack{I \in \mathcal{J}^k}}f_Idx^I,g=\sum_{\substack{I \in \mathcal{J}^k}}g_Idx^I\in \Lambda^k$ we define
\begin{align}\label{norm of differential form}
    \left \langle f,g\right \rangle:=f=\sum_{\substack{I \in \mathcal{J}^k}}f_Ig_I \text{ and } \left \lvert f\right \rvert:=\left \langle f,f\right \rangle^{\frac{1}{2}}.
\end{align}
We say $f \in L^p\left(\Omega;\Lambda^k\right)$ if each $f_I \in L^p\left(\Omega\right)$. Moreover, if  $\Omega \left(\subseteq \mathbb{R}^n\right)$ is open, we define $W^{m,p}\left(\Omega,\Lambda^k\right)$, $C^{r,\alpha}\left(\bar{\Omega};\Lambda^k\right)$, $C^r\left(\Omega; \Lambda^k\right)$ and other well known spaces in a similar way. The corresponding norms are defined component-wise for instance if $f\in W^{m,p}\left(\Omega;\Lambda^k\right)$ then we define
$\left \lVert f\right \rVert_{W^{m,p}\left(\Omega;\Lambda^k\right)}:=\sum_{I\in \mathcal{J}^k}\left \lVert f_I\right\rVert_{W^{m,p}\left(\Omega\right)}$.\\
\begin{definition}\label{Pullback definition}
Suppose, $\Omega,\Omega_0$ are open subsets of $\mathbb{R}^n$ and $\Phi$ is a diffeomorphism from $\Omega$ onto $\Omega_0$. If $\mathcal{G}=\sum_{1\leq i_1<\cdots<i_k\leq n} \mathcal{G}_{i_1\cdots i_k}dx^{i_1}\wedge \cdots \wedge dx^{i_k}$ is a $k$-differential form on $\Omega_0$.
 Then we define
 \begin{enumerate}
 \item $\Phi^* \left(\mathcal{G}\right):\Omega \rightarrow \Lambda^k$ the pullback of $\mathcal{G}$ by $\Phi$ as follows
 \begin{align*}
     \Phi^*\left(\mathcal{G}\right):=\sum_{1\leq i_1<\cdots<i_k\leq n} \left(\mathcal{G} _{i_1\cdots i_k}\circ \Phi \right)d\Phi^{i_1}\wedge \cdots \wedge d\Phi^{i_k}
 \end{align*}
 where
        $d\phi^s=\sum_{i=1}^n \frac{\partial \phi^s}{\partial x_j}dx^j$.
 \item $\mathcal{G}\circ \Phi: \Omega \rightarrow \Lambda^k$, another $k$-form on $\Omega$ by
 \begin{align*}
 \mathcal{G}\circ \Phi:=\sum_{1\leq i_1<\cdots<i_k\leq n} \left(\mathcal{G} _{i_1\cdots i_k}\circ \Phi \right) dx^{i_1}\wedge \cdots \wedge dx^{i_k}.
 \end{align*}
 \end{enumerate}
        \end{definition}
We briefly mention the properties of pullback which are necessary for our article:
\begin{align}
&\text{$\Phi^*$ and $d$ commute with each other} \\
& \text{For two diffeomorphisms $\Phi:\Omega\rightarrow \Omega_0$ and $\Psi:\Omega_0\rightarrow \Omega$ we have $\left(\Phi\circ \Psi\right)^*=\Psi^* \circ \Phi^*$}
\end{align}
\begin{definition}
Let, $0\leq k\leq n-1$ and $\Omega$ be an open bounded subset of $\mathbb{R}^n$. Let $f$ be in $L^1_{\operatorname{loc}}\left(\Omega; \Lambda^k\right)$, and there exists $v \in L^1_{\operatorname{loc}}\left(\Omega; \Lambda^{k+1}\right)$ such that 
\begin{align*}
    \int_{\Omega}\left\langle f, d^*\phi\right \rangle =-\int_{\Omega}\left \langle v, \phi\right\rangle \text{ for all $\phi \in C^\infty_c\left(\Omega, \Lambda^{k+1}\right)$.}
\end{align*}
Then we say $v$ is the weak exterior derivative of $f$. On the other hand, for $g \in L^1_{\operatorname{loc}}\left(\Omega; \Lambda^{k+1}\right)$ we say the weak interior derivative of $g$ exists if there exists $v \in L^1_{\operatorname{loc}}\left(\Omega; \Lambda^{k}\right)$ such that 
\begin{align*}
    \int_{\Omega}\left \langle g, d\phi\right\rangle =-\int_{\Omega}\left \langle v, \phi\right \rangle \text{ for all $\phi \in C^\infty_c\left(\Omega; \Lambda^{k}\right)$.}
\end{align*}
and we define $v$ as the weak interior derivative of $g$.
\end{definition}
Let $\Omega \left(\subset \mathbb{R}^n\right)$ be an open bounded set with smooth boundary. If $\nu\left(x\right)$ is the exterior unit normal vector at $x \in \partial \Omega$; then we identify $\nu$ with the $1$-form $\sum_{i=1}^n \nu_i dx^i$
(see \cite[Section ~ $3.2$]{csato2011pullback}). For such domain we mention some standard spaces involving differential forms (see \cite[Page~$5$]{MR3720379}):
        \begin{align*}&W^{1,p}_T\left(\Omega; \Lambda^k\right) := \left\{w \in W^{1,p}\left(\Omega; \Lambda^k\right): \nu \wedge w=0 \text{ on } \partial \Omega\right\},\\
		&W^{1,p}_{d^*, T}\left(\Omega; \Lambda^k\right):=\left \{w \in W^{1,p}_T\left(\Omega; \Lambda^k\right): d^*w=0 \text{ in } \Omega\right\},\\
		&\mathcal{H}_T\left(\Omega;\Lambda^k\right):=\left\{w \in W^{1,2}_T\left(\Omega; \Lambda^k\right):d^*w=0 \text{ and }dw=0 \text{ in } \Omega\right\},\\
        &W^{d,p}\left(\Omega; \Lambda^k\right) := \left \{w \in L^{p}\left(\Omega; \Lambda^k\right): dw \in L^{p}\left(\Omega; \Lambda^{k+1}\right)\right \},\\
        &W^{1,p}_N\left(\Omega; \Lambda^k\right) := \left\{w \in W^{1,p}\left(\Omega; \Lambda^k\right): \nu \lrcorner w=0 \text{ on } \partial \Omega\right\},\\
		&W^{1,p}_{d^*, N}\left(\Omega; \Lambda^k\right):=\left \{w \in W^{1,p}_N\left(\Omega; \Lambda^k\right): d^*w=0 \text{ in } \Omega\right\},\\
		&\mathcal{H}_N\left(\Omega;\Lambda^k\right):=\left\{w \in W^{1,2}_N\left(\Omega; \Lambda^k\right):d^*w=0 \text{ and }dw=0 \text{ in } \Omega\right\}.
        \end{align*}
\subsection{Inequalities and Standard Lemmas Involving Differential Forms}
In this subsection we recall some standard inequalities and lemmas involving differential forms.\\

For a differential form $\mathcal{G}$ on open bounded $\Omega \left(\subset \mathbb{R}^n\right)$, $K\subset \subset \Omega$ and $h\in {\mathbb{R}^n}$ with $0<\left \lvert h\right \rvert \leq d\left(K, \partial \Omega\right)$ we define
\begin{align*}
    \uptau_h \mathcal{G}\left(x\right):=\mathcal{G}\left(x+h\right)-\mathcal{G}\left(x\right) \text{ and } D_h\mathcal{G}\left(x\right)=\frac{\uptau_h \mathcal{G}\left(x\right)}{\left \lvert h\right \rvert} \text{ for $x$ in }  K.
\end{align*}
We also define
\begin{align*}
    \left(\mathcal{G}\right)_\Omega:=\fint_{\Omega} \mathcal{G}=\frac{1}{\left \lvert \Omega\right \rvert}\int_\Omega \mathcal{G}.
\end{align*}
\begin{lemma}\label{differencelemma1}
Suppose, $\gamma \geq 1, B_{R_0}\subset \subset B_R\subset \mathbb{R}^n, h\in \mathbb{R}^n$ with $\left \lvert h \right \rvert< \frac{R-R_0}{2}$ and $\mathcal{G} \in W^{1,\gamma}\left(B_R\right)$; then
\begin{align}\label{differencequotient}
	\int_{B_{R_0}}\left \lvert \uptau_h\mathcal{G} \right \rvert^{\gamma}\leq \left \lvert h \right \rvert^{\gamma}\int_{B_R}\left \lvert\nabla \mathcal{G}\right \rvert^{\gamma}.
\end{align}
(See \cite[Subsection $5.8.2$]{evans2010partial}).
    \end{lemma}
    The following lemma is an consequence of integration by parts formula on smooth functions and density arguments. 
    \begin{lemma}\label{Lemma <G,H>}
        Suppose, $1<\gamma<\infty,0<R_0<R<\infty, h\in \mathbb{R}^n$ with $\left \lvert h\right \rvert\leq \frac{R-R_0}{10^3}, \mathcal{G}\in L^{\gamma}\left(B_R; \Lambda^{k+1}\right)$ and $H\in W^{1,\frac{\gamma}{\gamma-1}}_0\left(B_r; \Lambda^{k+1}\right)$ then
        \begin{align}\label{<tau h G, H>}
            \int_{B_R}\left \langle \uptau_h \mathcal{G},H\right\rangle =-\sum_{j=1}^n\int_{B_R}h_j\left[\int_0^1\left \langle \mathcal{G}\left(x+th\right),\frac{\partial}{\partial x_j}H\left(x\right) \right \rangle dt\right] dx.
        \end{align}
    \end{lemma}
    \begin{theorem}\label{nabla d d*}
    Suppose, $0\leq k \leq n-1, 1<\gamma< \infty, R>0$ then there exists $ c=c(\gamma,R)>0$ such that, 
		\begin{align*}
		\left\lVert  u \right\rVert_{ W^{1,\gamma}\left( B_{R};\Lambda^k\right)}
            &\leq c\left[ \left\lVert du\right\rVert_{L^{\gamma} \left( B_{R} ; \Lambda^{k+1}\right)}+\left\lVert d^{\ast}u\right\rVert_{L^{\gamma} \left( B_{R} ; \Lambda^{k-1} \right) }\right].
		\end{align*}
		for all $u \in W^{1,\gamma}_T \left(B_R;\Lambda^k\right)\cup W^{1,\gamma}_N \left(B_R;\Lambda^k\right)$.\\
        (See \cite[Corollary $2.41$]{sil2016calculus}.)
    \end{theorem}
    Notice that Theorem \ref{nabla d d*} is the classical Poincar\'{e} inequality for $k=0$. Moreover it implies the following:
     \begin{corollary}\label{equivnorm}
		For $0\leq k \leq n-1, 1<\gamma< \infty, R>0, w \mapsto \left \lVert dw \right \rVert_{L^\gamma\left(B_R\right)}$ is equivalent (to the usual) norm on $W^{1,\gamma}_{d^*,T}\left(B_R;\Lambda^k\right)$ and $W^{1,\gamma}_{d^*,N}\left(B_R;\Lambda^k\right)$.   
	\end{corollary}

We know that for any ball $B_R$ in $\mathbb{R}^n$ is a contractible, bounded and smooth domain; therefore 
        $\mathcal{H}_T\left(B_R; \Lambda^k\right)=\left\{0\right\} \text{ for $0\leq k\leq n$}
        \text{ and }
        \mathcal{H}_N\left(B_R; \Lambda^k\right)= \left\{0\right\} \text{ for $1\leq k\leq n$}$.
    As a consequence we have the following Poincar\'{e} type lemmas.
    \begin{lemma}\label{poincare}
        Let $0\leq k \leq n$ be an integer and $\gamma$ be a real number in $\left(1,\infty\right)$. Let $f \in L^\gamma\left(B_R; \Lambda^{k+1}\right)$ with $df=0$ in $B_{R}$ and $\nu \wedge f=0$ on $\partial B_R$ (where $\nu\left(x\right)$ is the exterior unit normal vector at $x\in B_R$);
        then there exists $\theta \in W^{1,\gamma}_0\left(B_R;\Lambda^k\right)$
        \begin{align}
            \left\lbrace
            \begin{aligned}d\theta=f &\text{ in } B_R,\\
            \theta=0 &\text{ on } \partial B_R.
            \end{aligned}
            \right.
        \end{align}
        (See \cite[Theorem $2.47$]{sil2016calculus})
    \end{lemma}
     
\begin{lemma}\label{1-div-curl}
Let $1\leq k\leq n-1$ be an integer and $1<\gamma< \infty$. If $f \in L^\gamma\left(B_R;\Lambda^{k+1}\right)$ with $df=0$ in $B_R$ and $\nu \wedge f=0$ on $\partial B_R$ (where $\nu\left(x\right)$ is the exterior unit normal vector at $x\in B_R$);
then there exists $\theta \in W^{1,\gamma}_{d^*,T}\left(B_R; \Lambda^k\right)$ such that
\begin{align}
\left \lbrace
\begin{aligned}
             d\theta=f &\text{ in } B_R,\\
             d^*\theta=0 & \text{ in } B_R,\\
             \nu \wedge \theta=0 &\text{ on } \partial B_R.
            \end{aligned}
            \right.
		\end{align}
(See \cite[Theorem $14$]{MR3720379}).
\end{lemma}
\begin{lemma}\label{1-div-curl-2}
Let $1\leq k\leq n-1$ be an integer and $1<\gamma< \infty$. Let $f \in L^\gamma\left(B_R;\Lambda^{k+1}\right)$ with $df=0$ in $B_R$;
then there exists $\theta \in W^{1,\gamma}_{d^*,N}\left(B_R; \Lambda^k\right)$ such that
\begin{align}
\left \lbrace
\begin{aligned}
	d\theta=f &\text{ in } B_R,\\
    d^*\theta=0 & \text{ in } B_R,\\
    \nu \lrcorner \theta=0 &\text{ on } \partial B_R,
\end{aligned}
\right.
\end{align}
where for $x\in \partial B_R$, $\nu\left(x\right)$ is the exterior unit normal at $x$.\\
(See \cite[Theorem $14$]{MR3720379}).
\end{lemma}
\subsection{Inequalities Related to the Auxiliary Functions}
Now we collect the inequalities involving the auxiliary functions, which will be frequently used throughout this article. 
\begin{lemma}\label{Acerbi Fusco Lemma} For $-1<\gamma<0$ we have
\begin{align*}
    \int_{0}^1\left \lvert \zeta+t\left(\xi-\zeta\right)\right \rvert^{\gamma }dt\sim \left(\left \lvert \zeta\right \rvert+\left \lvert \xi\right \rvert\right)^{\gamma } \text{ for } \zeta, \xi \in \mathbb{R}^N.
\end{align*}
(See \cite[Lemma $2.1$]{MR997847}).
\end{lemma}
\begin{lemma}\label{G near to I}
Let us assume that, $1<p<\infty$ and $G:B_r\left(0\right)\rightarrow GL\left(N,\mathbb{R}\right)$ is measurable, symmetric, bounded satisfying \eqref{ellipticity}
and \eqref{1<p<2 and p<2lambda first mention}. Then for almost every $x\in \Omega$ and $\zeta, \xi \in \mathbb{R}^N$ we have
\begin{align}\label{Kp}
    &\left \lvert \left \langle G \xi, \xi\right \rangle^{\frac{p-2}{2}}G\xi-\left \lvert \xi\right\rvert^{p-2}\xi\right \rvert\leq K_p\left \lVert G-I\right \rVert_{L^\infty\left(B_r\left(0\right)\right)} \left \lvert \xi\right \rvert^{p-1}, 
\end{align}
\begin{align}\label{c1p}
    \left \langle  \left \langle G \xi, \xi\right \rangle^{\frac{p-2}{2}}G\xi-\left \langle G \zeta, \zeta\right \rangle^{\frac{p-2}{2}}G\zeta, \xi-\zeta \right \rangle \geq  c_1\left(p\right)  \left \lvert V_p\left(\xi\right)-V_p\left(\zeta\right)\right \rvert^2,
\end{align}
\begin{align}\label{G V upper}
    \left \lvert \left \langle G \xi, \xi\right \rangle^{\frac{p-2}{2}}G\xi-\left \langle G \zeta, \zeta\right \rangle^{\frac{p-2}{2}}G\zeta\right \rvert&\leq c_2\left(p\right)\left(\left \lvert \zeta\right \rvert+\left \lvert \xi\right \rvert\right)^{p-2}\left \lvert \xi-\zeta\right \rvert \sim \left \lvert V_p\left(\xi\right)-V_p\left(\zeta\right)\right \rvert^2.
\end{align}
\end{lemma}
\begin{proof}
We start with the proof of \eqref{Kp}.
\begin{align*}
    \left \lvert \left \langle G \xi, \xi\right \rangle^{\frac{p-2}{2}}G\xi-\left \lvert \xi\right\rvert^{p-2}\xi\right \rvert
    \leq \int_0^1\left \lvert  \frac{d}{dt}\left[\left \langle \left(I+t\left(G-I\right)\right)\xi,\xi \right \rangle^{\frac{p-2}{2}}\left(I+t\left(G-I\right)\xi\right)\right]\right \rvert dt.
\end{align*}
Now
\begin{align*}
    &\left \lvert \frac{d}{dt}\left[\left \langle \left(I+t\left(G-I\right)\right)\xi,\xi \right \rangle^{\frac{p-2}{2}}\left(I+t\left(G-I\right)\xi\right)\right]\right \rvert\\
    &\qquad \qquad\leq  \begin{aligned}[t]\left \lvert\frac{p-2}{2}\right\rvert &\left.\left \langle \left(I+t\left(G-I\right)\right)\xi,\xi \right \rangle^{\frac{p-4}{2}} \left\lvert G-I\right\rvert  \left\lvert I+t\left(G-I\right)\right\rvert \left \lvert \xi\right \rvert^3\right. \\
    &\qquad\qquad+\left \langle \left(I+t\left(G-I\right)\right)\xi,\xi \right \rangle^{\frac{p-2}{2}}\left\lvert G-I\right\rvert\left \lvert \xi\right \rvert,
    \end{aligned}
\end{align*}
and $\left[\left(1-t\right)+t\lambda_G\right] \left \lvert \xi \right \rvert^2\leq \left \langle \left(I+t\left(G-I\right)\right)\xi,\xi \right \rangle \leq \left[\left(1-t\right)+t\Lambda_G\right] \left \lvert \xi \right \rvert^2$. Thus for $\frac{1}{2}\leq \lambda_G\leq \Lambda_G\leq \frac{3}{2}$ we get
    \begin{align*}
        \left[\left(1-t\right)+\frac{t}{2}\right] \left \lvert \xi \right \rvert^2=\left \langle \left(\left(1-t\right)I+tG\right)\xi,\xi \right \rangle\leq \left[\left(1-t\right)+\frac{3t}{2}\right] \left \lvert \xi \right \rvert^2.
    \end{align*}
    In that case for $t\in \left[0,1\right]$, we have
    $\frac{1}{2}\left \lvert \xi \right \rvert^2\leq \left \langle \left(\left(1-t\right)I+tG\right)\xi,\xi \right \rangle\leq \frac{3}{2}\left \lvert \xi \right \rvert^2$.
     As a consequence \eqref{Kp} follows.\\
        Next we  prove \eqref{c1p}.\\
        For $p\in \left[2,\infty\right)$ T. Iwaniec proved in $\left(1.9\right)$ of \cite{bwmeta1.element.bwnjournal-article-smv75i1p22bwm} proved that (without the assumptions in \eqref{1<p<2 and p<2lambda first mention}):
        \begin{align*}
    \left \langle  \left \langle G \xi, \xi\right \rangle^{\frac{p-2}{2}}G\xi-\left \langle G \zeta, \zeta\right \rangle^{\frac{p-2}{2}}G\zeta, \xi-\zeta \right \rangle \geq  c\left(p\right)  \left(\left \lvert \xi\right \rvert+\left \lvert \zeta\right \rvert\right)^{{p-2}}\left \lvert \xi-\zeta\right\rvert^2,
\end{align*} and using \eqref{V sim} with it we get \eqref{c1p} for $2\leq p<\infty$.
Here we present a short proof for the case $1<p<2$ assuming \eqref{1<p<2 and p<2lambda first mention}:
    \begin{align*}
        \frac{d}{dt}&\left[\left \langle  \left \langle G \left(\xi+t\left(\zeta-\xi\right)\right), \xi+t\left(\zeta-\xi\right)\right \rangle^{\frac{p-2}{2}}G\left(\xi+t\left(\zeta-\xi\right)\right),\zeta-\xi\right  \rangle \right] \\
        &\begin{aligned}[t]=&  \left \langle G \left(\xi+t\left(\zeta-\xi\right)\right), \xi+t\left(\zeta-\xi\right)\right \rangle^{\frac{p-2}{2}}\left \langle G\left(\zeta-\xi\right),\zeta-\xi \right  \rangle\\
        &+\left(p-2\right)  \left \langle G \left(\xi+t\left(\zeta-\xi\right)\right), \xi+t\left(\zeta-\xi\right)\right \rangle^{\frac{p-4}{2}}\left \langle G\left(\xi+t\left(\zeta-\xi\right)\right),\zeta-\xi \right  \rangle^2
        \end{aligned}\\
        &\begin{aligned}[t]\geq & \Lambda_G^{\frac{p-2}{2}} \lambda_G\left \lvert \xi+t\left(\zeta-\xi\right)\right \rvert^{p-2}\left\lvert \zeta-\xi \right  \rvert^2\\
        &+\left(p-2\right)\lambda_G^{\frac{p-4}{2}} \Lambda_G^2  \left \lvert \xi+t\left(\zeta-\xi\right)\right \rvert ^{p-2}\left \lvert \zeta-\xi \right  \rvert^2 \text{[Since $p-4<p-2<0$]}.
        \end{aligned}
    \end{align*}
    Define, $\epsilon_p:=\frac{p-1}{6-2p}$; in that case by \eqref{1<p<2 and p<2lambda first mention} we have
    \begin{align*}
        \frac{d}{dt}&\left[\left \langle  \left \langle G \xi+t\left(\zeta-\xi\right), \xi+t\left(\zeta-\xi\right)\right \rangle^{\frac{p-2}{2}}G\left(\xi+t\left(\zeta-\xi\right)\right),\zeta-\xi\right  \rangle \right] \\
        &\begin{aligned}[t]\geq  \left(1-\epsilon_p\right)&\left \lvert \xi+t\left(\zeta-\xi\right)\right \rvert^{p-2}\left\lvert \zeta-\xi \right  \rvert^2\\
        &+\left(p-2\right)\left(1+\epsilon\right)\left \lvert \xi+t\left(\zeta-\xi\right)\right \rvert ^{p-2}\left \lvert \zeta-\xi \right  \rvert^2
       \end{aligned}\\
        &\begin{aligned}[t]\geq  \left[p-1-\epsilon_p\left(3-p\right)\right]\left \lvert \xi+t\left(\zeta-\xi\right)\right \rvert ^{p-2}\left \lvert \zeta-\xi \right  \rvert^2 \qquad \left[\text{Since, $p-2<0$}\right]
       \end{aligned}\\
    \end{align*}
        which gives 
        \begin{align*}
            \frac{d}{dt}&\left[\left \langle  \left \langle G \xi+t\left(\zeta-\xi\right), \xi+t\left(\zeta-\xi\right)\right \rangle^{\frac{p-2}{2}}G\left(\xi+t\left(\zeta-\xi\right)\right),\zeta-\xi\right  \rangle \right] \\
        &\begin{aligned}[t]\geq  \frac{p-1}{2}\left \lvert \xi+t\left(\zeta-\xi\right)\right \rvert ^{p-2}\left \lvert \zeta-\xi \right  \rvert^2.
       \end{aligned}\\
        \end{align*}
By Lemma \ref{Acerbi Fusco Lemma} with $\gamma=p-2$ and \eqref{V sim} we complete the proof of \eqref{c1p}.\\
Finally, notice that
\begin{align*}
&\left \lvert \left \langle G \xi, \xi\right \rangle^{\frac{p-2}{2}}G\xi-\left \langle G \zeta, \zeta\right \rangle^{\frac{p-2}{2}}G\zeta\right \rvert\\
&= \int_0^1 \left \lvert \frac{d}{dt}\left[\left \langle G \left(\zeta +t\left(\xi-\zeta\right)\right), \left(\zeta +t\left(\xi-\zeta\right)\right)\right \rangle^{\frac{p-2}{2}}G\left(\zeta +t\left(\xi-\zeta\right)\right)\right]\right \rvert dt.
\end{align*}
Straight forward calculations with \eqref{ellipticity} gives the follwoing:
        \begin{align*}
            &\left \lvert \frac{d}{dt}\left[\left \langle  \left(\zeta +t\left(\xi-\zeta\right)\right), \left(\zeta +t\left(\xi-\zeta\right)\right)\right \rangle^{\frac{p-2}{2}}G\left(\zeta +t\left(\xi-\zeta\right)\right)\right]\right\rvert\\
            &\qquad \leq \begin{cases}\left(p-1\right)\Lambda_G^{\frac{p}{2}}\left(\left \lvert \zeta\right \rvert+\left \lvert \xi\right \rvert\right)^{p-2}\left \lvert \xi-\zeta\right \rvert &\text{ if } p\geq 4,\\
            \left[\Lambda_G^{\frac{p}{2}}+\left(p-2\right)\Lambda_G^{\frac{p-4}{2}}\lambda_G^2\right]\left(\left \lvert \zeta\right \rvert+\left \lvert \xi\right \rvert\right)^{p-2}\left \lvert \xi-\zeta\right \rvert &\text{ if } 2\leq p< 4,\\
            \left[\Lambda_G\lambda_G^{\frac{p-2}{2}}+\left(2-p\right)\Lambda_G^{\frac{p-4}{2}}\lambda_G^2\right]\left \lvert \zeta +t\left(\xi-\zeta\right)\right \rvert^{p-2}\left \lvert \xi-\zeta\right \rvert &\text{ if } 1< p< 2.
            \end{cases}
        \end{align*}
        Combining all three cases mentioned above and using Lemma \ref{Acerbi Fusco Lemma} we derive \eqref{G V upper}.
        \end{proof}
        We end this subsection mentioning one standard result from calculus of variation (see~\cite[Theorem~$3.20$]{dacorogna2007direct} with $a=b= 0, q=\gamma \in \left(1,\infty\right) $ and $f=\mathcal{H}$):
        \begin{lemma}\label{Convex Functional Lemma from Daco's Book}
		Let $\Omega$ be an open subset of $\mathbb{R}^n$, $1<\gamma<\infty,M\in \mathbb{N}$. If $\mathcal{H}:\Omega\times \mathbb{R}^N \rightarrow \left[0, \infty\right]$ is convex with respect to the second variable and $\zeta_m \rightharpoonup \zeta$ in $L^\gamma\left(\Omega; \mathbb{R}^M\right)$, then 
		\begin{align}
			\int_\Omega \mathcal{H}\left(x,\zeta\left(x\right)\right)dx \leq \liminf_{m\rightarrow \infty} \int_\Omega \mathcal{H}\left(x,\zeta_m\left(x\right)\right)dx.
		\end{align}    
	\end{lemma}
\section{Proof of Theorem \ref{G Theorem}}
We divide the proof of Theorem \ref{G Theorem} into three steps. First with some higher regularity condition on $\widetilde{F}$, we prove the higher integrability of $d\widetilde{W}$. Then with the help of Hardy-Littlewood maximal inequality, Fefferman-Stein inequality and the classical regularity estimate of Ulhenbeck we prove strong $\left(q,q\right)$ type estimate locally.
Finally with density argument we prove Theorem \ref{G Theorem}.
Before we start our proof we mention few remarks about the weak formulation of \eqref{G weak formulation in theorem}. 
We say $\widetilde{W} \in W^{1,p}\left(B_r\left(0\right);\Lambda^k\right)$ weakly solves \eqref{G weak formulation in theorem} if for all $\chi\in C_c^\infty\left(B_{r}\left(0\right);\Lambda^k\right)$ we have the following
    \begin{align}\label{final weak formulation Cc}
        \int_{B_r\left(0\right)}\left \langle \left \langle G \left(d\widetilde{W}+\widetilde{F}\right), d\widetilde{W}+\widetilde{F}\right  \rangle^{\frac{p-2}{2}}G \left(d \widetilde{W}+\widetilde{F}\right),d\chi\right \rangle=0
    \end{align}
    Using the fact that $C_c^\infty\left(B_r\left(0\right);\Lambda^k\right)$ is a dense subspace of $W^{1,p}_0\left(B_r\left(0\right);\Lambda^k\right)$ we derive another equivalent statement:\\
    For all $\chi\in W^{1,p}_0\left(B_{r}\left(0\right)\right)$ we have the following
    \begin{align}\label{final weak formulation}
        \int_{B_r\left(0\right)}\left \langle \left \langle G \left(d\widetilde{W}+\widetilde{F}\right), d\widetilde{W}+\widetilde{F}\right  \rangle^{\frac{p-2}{2}}G \left(d \widetilde{W}+\widetilde{F}\right),d\chi\right \rangle=0.
    \end{align}
    Finally using Lemma \ref{poincare} and Lemma \ref{1-div-curl} we derive the following equivalent weak formulation of \eqref{G weak formulation in theorem}:\\
    For all $\chi\in W^{1,p}_{d^*,T}\left(B_{r}\left(0\right)\right)$ we have the following
    \begin{align}\label{d* T weal formulation}
        \int_{B_r\left(0\right)}\left \langle \left \langle G \left(d\widetilde{W}+\widetilde{F}\right), d\widetilde{W}+\widetilde{F}\right  \rangle^{\frac{p-2}{2}}G \left(d \widetilde{W}+\widetilde{F}\right),d\chi\right \rangle=0.
    \end{align}
    We divide the proof of Theorem \ref{G Theorem} into following propositions:
         \begin{proposition}\label{Proposition with F in W2,p}
         Let $0\leq k\leq n-1, 1<p<\infty$ and let $G:B_{r}\left(0\right)\rightarrow GL\left({n\choose k},\mathbb{R}\right)$ be a measurable, symmetric, bounded function satisfying \eqref{ellipticity} and \eqref{1<p<2 and p<2lambda first mention}. We assume that $\widetilde{F}$ is in $W^{2,\widetilde{q}}\left(B_r\left(0\right);\Lambda^{k+1}\right)$ with $\widetilde{q}$ satisfying \eqref{qtilde}. If $\widetilde{W}\in W^{1,p}\left(B_r\left(0\right);\Lambda^k\right)$ weakly solves
    \begin{align}\label{el fo W}
        d^*\left(\left \langle G \left(d\widetilde{W}+\widetilde{F}\right), d\widetilde{W}+\widetilde{F}\right  \rangle^{\frac{p-2}{2}}G \left(d \widetilde{W}+\widetilde{F}\right)\right)=0 \text{ in } B_r\left(0\right).
    \end{align}
    Then $d\widetilde{W} \in L^{\widetilde{q}}_{\operatorname{loc}}\left(B_r\left(0\right)\right)$.
    \end{proposition}
    \begin{proposition}\label{Application of HL FS KU Proposition}
    Suppose, $0\leq k\leq n-1, 1<p<\widetilde{q}<\infty$ with $\widetilde{q}$ satisfying \eqref{qtilde}, $\widetilde{F}\in L^{\widetilde{q}}\left(B_r\left(0\right);\Lambda^{k+1}\right)$ and $G:B_{r}\left(0\right)\rightarrow GL\left({n\choose k},\mathbb{R}\right)$ is a measurable, symmetric, bounded function satisfying \eqref{ellipticity} and $*$-condition (see \eqref{delta}, \eqref{r}, \eqref{1<p<2 and p<2lambda first mention}). If $\widetilde{W}\in W^{1,p}\left(B_r\left(0\right);\Lambda^k\right)$ with $d\widetilde{W} \in L^{\widetilde{q}}_{\operatorname{loc}}\left(B_r\left(0\right);\Lambda^{k+1}\right)$ weakly solves
    \begin{align}\label{2nd Proposition eqn. mention}
        d^*\left(\left \langle G \left(d\widetilde{W}+\widetilde{F}\right), d\widetilde{W}+\widetilde{F}\right  \rangle^{\frac{p-2}{2}}G \left(d \widetilde{W}+\widetilde{F}\right)\right)=0 \text{ in } B_r\left(0\right).
    \end{align}
Then there exists a cube $Q_0\subset\subset B_r\left(0\right)$ with sides parallel to the coordinate axes such that
    \begin{align*}
        \int_{Q_0}\left \lvert d\widetilde{W}+\widetilde{F} \right \rvert^{\widetilde{q}}  &\leq \begin{aligned}[t]
        &c\left(n,\tilde{q},p,\alpha_1,c_U, r\right) \left[\left[\int_{Q_0}\left \lvert d\widetilde{W}+\widetilde{F} \right \rvert^p\right]^{\frac{\tilde{q}}{p}}
        +\int_{Q_0}\left \lvert \widetilde{F}\right\rvert^{\tilde{q}}\right].
        \end{aligned}
    \end{align*}
    Here $\left \lvert Q_0\right \rvert$ can be chosen depending on $n,p,\widetilde{q},k$ and $r$ only.
\end{proposition}
We start with the proof of Proposition \ref{Proposition with F in W2,p}:
    \begin{proof}
    Choose non-negative $\eta \in C^\infty_c\left(\mathbb{R}^n\right)$ with
    $\operatorname{supp}\left(\eta\right)\subseteq B_{\frac{r}{2}}\left(0\right)$ and\\ $\eta \equiv 1$ in $B_{\frac{r}{4}}\left(0\right)$. In distributional  sense we have $d\left(d\widetilde{W}\right)=0$ in $B_r\left(0\right)$; so using Lemma \ref{1-div-curl-2} we have $\widetilde{\widetilde{W}}\in W^{1,p}_{d^*,N}\left(B_r\left(0\right)\right)$ such that the $d\widetilde{\widetilde{W}}=d\widetilde{W} \text{ in } B_r\left(0\right)$.
    Taking  $\chi=\uptau_{-h}\left(\eta^2\uptau_h \widetilde{\widetilde{W}}\right)$ as test function in \eqref{final weak formulation} for $\left \lvert h\right \rvert\leq \frac{r}{1000}$, we arrive at the following: 
        \begin{align*}
            -2&\int_{B_r\left(0\right)}\left \langle \uptau_h\left[\left \langle G \left(d\widetilde{W}+\widetilde{F}\right), d\widetilde{W}+\widetilde{F}\right  \rangle^{\frac{p-2}{2}}G \left(d \widetilde{W}+\widetilde{F}\right)\right],d\eta \wedge \uptau_h \widetilde{\widetilde{W}}\right \rangle\eta  \\
            &=\int_{B_r\left(0\right)}\left \langle \uptau_h\left[\left \langle G \left(d\widetilde{W}+\widetilde{F}\right), d\widetilde{W}+\widetilde{F}\right  \rangle^{\frac{p-2}{2}}G \left(d \widetilde{W}+\widetilde{F}\right)\right],\uptau_h d \widetilde{\widetilde{W}}\right \rangle \eta^2\\
            &=\begin{aligned}[t]\int_{B_r\left(0\right)}&\left \langle \uptau_h\left[\left \langle G \left(d\widetilde{W}+\widetilde{F}\right), d\widetilde{W}+\widetilde{F}\right  \rangle^{\frac{p-2}{2}}G \left(d \widetilde{W}+\widetilde{F}\right)\right],\uptau_h d\widetilde{{{W}}}\right \rangle \eta^2.\\
            & \qquad\qquad\qquad\qquad\qquad\qquad\qquad \qquad \text{ [Since, we have $d\widetilde{\widetilde{W}}=d\widetilde{W}$ in $B_r\left(0\right)$.]}
            \end{aligned}
        \end{align*}
        On the other hand
        \begin{align*}
        &\begin{aligned}[t]\int_{B_r\left(0\right)} \left \langle \uptau_h\left[\left \langle G\left(d\widetilde{W}+\widetilde{F}\right), d\widetilde{W}+\widetilde{F}\right \rangle^{\frac{p-2}{2}}G \left(d\widetilde{W}+\widetilde{F}\right)\right],\uptau_hd\widetilde{W}\right \rangle\eta^2
        \end{aligned}\\
        &=\begin{aligned}[t] &
        \begin{aligned}[t]\int_{B_r\left(0\right)} \left \langle \uptau_h\left[\left \langle G\left(d\widetilde{W}+\widetilde{F}\right), d\widetilde{W}+\widetilde{F}\right \rangle^{\frac{p-2}{2}}G \left(d\widetilde{W}+\widetilde{F}\right)\right],\uptau_h\left(d\widetilde{W}+\widetilde{F}\right)\right \rangle\eta^2
        \end{aligned}\\
        &-\int_{B_r\left(0\right)} \left \langle \uptau_h\left[\left \langle G\left(d\widetilde{W}+\widetilde{F}\right), d\widetilde{W}+\widetilde{F}\right \rangle^{\frac{p-2}{2}}G \left(d\widetilde{W}+\widetilde{F}\right)\right],\uptau_h\widetilde{F}\right \rangle\eta^2.
        \end{aligned}
    \end{align*}
    Define,
    \begin{align*}
        &A_1:=\int_{B_r\left(0\right)} \left \langle \uptau_h\left[\left \langle G\left(d\widetilde{W}+\widetilde{F}\right), d\widetilde{W}+\widetilde{F}\right \rangle^{\frac{p-2}{2}}G \left(d\widetilde{W}+\widetilde{F}\right)\right],\uptau_hd\widetilde{W}\right \rangle\eta^2,\\
        &A_2:=\int_{B_r\left(0\right)} \left \langle \uptau_h\left[\left \langle G\left(d\widetilde{W}+\widetilde{F}\right), d\widetilde{W}+\widetilde{F}\right \rangle^{\frac{p-2}{2}}G \left(d\widetilde{W}+\widetilde{F}\right)\right],\uptau_h\left(d\widetilde{W}+\widetilde{F}\right)\right \rangle\eta^2,\\
        &A_3:=\int_{B_r\left(0\right)} \left \langle \uptau_h\left[\left \langle G\left(d\widetilde{W}+\widetilde{F}\right), d\widetilde{W}+\widetilde{F}\right \rangle^{\frac{p-2}{2}}G \left(d\widetilde{W}+\widetilde{F}\right)\right],\uptau_h\widetilde{F}\right \rangle\eta^2.
    \end{align*}
i.e. $A_1=A_2-A_3$ which implies $A_1\leq \left \lvert A_2\right \rvert+\left \lvert A_3\right \rvert$.\\
By our assumption \eqref{ellipticity},\eqref{1<p<2 and p<2lambda first mention} are true, hence by Lemma \ref{G near to I} for $x\in B_r\left(0\right)$ the inequality \eqref{c1p} holds and  we have the following estimate:
\begin{align}\label{A_1}
 c_1\left(p\right) \int_{B_r\left(0\right)}\left \lvert \uptau_h V_p\left(d\widetilde{W}+\widetilde{F}\right)\right \rvert^2 \eta^2 \leq A_1.
\end{align}
We split the proof of the Proposition \ref{Proposition with F in W2,p} in two cases $2\leq p<\infty$ and $1<p<2$.\\ 
Let us first assume that $2\leq p<\infty$. 
\begin{align*}
A_2&\begin{aligned}[t]\leq \frac{c}{r}\int_{B_r\left(0\right)} \left \lvert  \uptau_h\left[\left \langle G\left(d\widetilde{W}+\widetilde{F}\right), d\widetilde{W}+\widetilde{F}\right \rangle^{\frac{p-2}{2}}G \left(d\widetilde{W}+\widetilde{F}\right)\right]\right \rvert \left \lvert \uptau_h \widetilde{\widetilde{W}}\right \rvert \eta \left[\text{Since, $\left \lvert \nabla \eta\right \rvert \lesssim \frac{1}{r}$}\right]
\end{aligned}\\
&\begin{aligned}[t]\leq \frac{c}{r}&\int_{B_r\left(0\right)}\left( \left \lvert  \left(d\widetilde{W}+\widetilde{F}\right)\left(x+h\right)\right \rvert+\left \lvert  \left(d\widetilde{W}+\widetilde{F}\right)\left(x\right)\right \rvert \right)^{\frac{p-2}{2}}  \left \lvert  \uptau_h  \left(d\widetilde{W}+\widetilde{F}\right)\right \rvert \eta \\
&\times \left( \left \lvert  \left(d\widetilde{W}+\widetilde{F}\right)\left(x+h\right)\right \rvert+\left \lvert  \left(d\widetilde{W}+\widetilde{F}\right)\left(x\right)\right \rvert \right)^{\frac{p-2}{2}}\left \lvert \uptau_h \widetilde{\widetilde{W}}\right \rvert \left[\text{By \eqref{G V upper}}\right]
\end{aligned}\\
&\begin{aligned}[t]\sim \frac{1}{r}&\int_{B_r\left(0\right)}   \left \lvert  \uptau_h V_p  \left(d\widetilde{W}+\widetilde{F}\right)\right \rvert \eta\\
& \times\left( \left \lvert  \left(d\widetilde{W}+\widetilde{F}\right)\left(x+h\right)\right \rvert+\left \lvert  \left(d\widetilde{W}+\widetilde{F}\right)\left(x\right)\right \rvert \right)^{\frac{p-2}{2}} \left \lvert \uptau_h \widetilde{\widetilde{W}}\right \rvert  \quad\left[\text{By \eqref{V sim}}\right]
\end{aligned}\\
&\quad\begin{aligned}[t]
\lesssim \epsilon\int_{B_r\left(0\right)} &\left \lvert  \uptau_hV_p  \left(d\widetilde{W}+\widetilde{F}\right)\right \rvert ^2\eta^2\\
&+\frac{c_\epsilon}{r^2}\int_{\operatorname{supp}\left(\eta\right)}\left[\left \lvert \left(d\widetilde{W}+\widetilde{F}\right)\left(x+h\right)\right \rvert+\left \lvert  \left(d\widetilde{W}+\widetilde{F}\right)\left(x\right)\right \rvert\right]^{p-2}\left \lvert \uptau_h \widetilde{\widetilde{W}}\right \rvert^2.
\end{aligned}\\
\end{align*}
In the last line we have used Young's Inequality with exponents $\left(2,2\right)$. On the other hand since, $p-2$ is non-negative, by suitable use of the Young's and H\"{o}lder's inequality with exponents $\left(\frac{p}{p-2}.\frac{p}{2}\right)$ and using \eqref{differencequotient} from  the Lemma \ref{differencelemma1} we get the following estimate:
\begin{align*}
&\int_{\operatorname{supp}\left(\eta\right)}\left[\left \lvert \left(d\widetilde{W}+\widetilde{F}\right)\left(x+h\right)\right \rvert+\left \lvert  \left(d\widetilde{W}+\widetilde{F}\right)\left(x\right)\right \rvert\right]^{p-2}\left \lvert \uptau_h \widetilde{\widetilde{W}}\right \rvert^2\\
&\quad\begin{aligned}[t]
\lesssim 
&\left \lvert h\right \rvert^2\left[\int_{B_r\left(0\right)}\left \lvert d\widetilde{W}+\widetilde{F}\right \rvert^p\right]^{\frac{p-2}{p}}\left[\int_{B_r\left(0\right)}\left \lvert \nabla \widetilde{\widetilde{W}}\right \rvert^p\right]^{\frac{2}{p}}
\end{aligned}\\
&\quad \begin{aligned}[t]\lesssim \left \lvert h\right \rvert^2\left[\int_{B_r\left(0\right)} \left \lvert \nabla \widetilde{\widetilde{W}}\right \rvert^p+\int_B\left \lvert \widetilde{F}\right \rvert^p\right]
\end{aligned}\\&\quad
\begin{aligned}[t]\lesssim \left \lvert h\right \rvert^2\left[\int_{B_r\left(0\right)} \left \lvert d \widetilde{\widetilde{W}}\right \rvert^p+\int_B\left \lvert \widetilde{F}\right \rvert^p\right] \left[\text{by Corollary \ref{equivnorm}; since $\widetilde{\widetilde{W}}\in W^{1,p}_{d^*,N}\left(B_r\left(0\right)\right)$ }\right]
\end{aligned}\\
&\quad\begin{aligned}[t]=\left \lvert h\right \rvert^2\left[\int_{B_r\left(0\right)} \left \lvert d \widetilde{{W}}\right \rvert^p+\int_B\left \lvert \widetilde{F}\right \rvert^p\right] \left[\text{since, $d\widetilde{\widetilde{W}}=d\widetilde{W}$ in $B_r\left(0\right)$}\right].
\end{aligned}
\end{align*}
Combining we get,
\begin{align}\label{A_2}
&\left \lvert A_2 \right \rvert \quad \begin{aligned}[t]
\lesssim \epsilon\int_{B_r\left(0\right)} &\left \lvert  \uptau_hV_p  \left(d\widetilde{W}+\widetilde{F}\right)\right \rvert ^2\eta^2+\frac{c_\epsilon}{r^2}\left \lvert h\right \rvert^2\left[\int_{B_r\left(0\right)} \left \lvert d \widetilde{{W}}\right \rvert^p+\int_B\left \lvert \widetilde{F}\right \rvert^p\right],
    \end{aligned}
\end{align}
    and similarly
    \begin{align}\label{A_3}
        \left \lvert A_3\right \rvert\lesssim \begin{aligned}[t]\epsilon \int_{B_r\left(0\right)} &\left \lvert \uptau_h V_p \left(d\widetilde{W}+\widetilde{F}\right)\right \rvert^2\eta^2+c_\epsilon\left \lvert h\right \rvert^2\int_{B_r\left(0\right)} \left [\left \lvert d\widetilde{W}+\widetilde{F}\right \rvert^p+\left \lvert \nabla \widetilde{F}\right \rvert^p\right].
        \end{aligned}
    \end{align}
Recall that $A_1 \leq \left \lvert A_2\right \rvert +\left \lvert A_3\right \rvert$, hence using \eqref{A_1}, \eqref{A_2}, \eqref{A_3} and 
choosing small $\epsilon>0$ we get
\begin{align*}\int_{B_{\frac{r}{4}}\left(0\right)} &\left \lvert \uptau_h V_p \left(d\widetilde{W}+\widetilde{F}\right)\right \rvert^2\leq c\left \lvert h\right \rvert^2\int_{B_r\left(0\right)} \left [\left \lvert d\widetilde{W}+\widetilde{F}\right \rvert^p+\left \lvert \nabla \widetilde{F}\right \rvert^p\right],
\end{align*} 
i.e.
\begin{align*}\int_{B_{\frac{r}{4}}\left(0\right)} &\left \lvert D_h V_p \left(d\widetilde{W}+\widetilde{F}\right)\right \rvert^2\leq c\int_{B_r\left(0\right)} \left [\left \lvert d\widetilde{W}+\widetilde{F}\right \rvert^p+\left \lvert \nabla \widetilde{F}\right \rvert^p\right].
\end{align*}
On the other hand, 
\begin{align*}
\int_{B_r\left(0\right)} &\left \lvert V_p \left(d\widetilde{W}+\widetilde{F}\right)\right \rvert^2= \int_{B_r\left(0\right)} \left \lvert d\widetilde{W}+\widetilde{F}\right \rvert^p<\infty.
\end{align*}
Thus, $V_p\left(d\widetilde{W}+\widetilde{F}\right)\in W^{1,2}_{\operatorname{loc}}\left(B_r\left(0\right)\right)$ and using Sobolev embedding $V_p\left(d\widetilde{W}+\widetilde{F}\right)\in L^{\frac{2n}{n-2}}_{\operatorname{loc}}B_r\left(0\right)$ for $n\geq 3$ and $V_p\left(d\widetilde{W}+\widetilde{F}\right)\in L^{\gamma}_{\operatorname{loc}}B_r\left(0\right)$ for all $\gamma$ in $\left[2,\infty\right)$. Therefore $d\widetilde{W}+\widetilde{F}\in L^{\widetilde{q}}_{\operatorname{loc}}B_r\left(0\right)$ where $\widetilde{q}$ is given in \eqref{qtilde}.\\

Now we deal with the case $1<p<2$. We split this case into  several Lemmas.\\
\begin{lemma}\label{Lemma 1<p<2 diff W diff VW}
    For $1<p<2, B_R \subset \subset B_r\left(0\right)$ and $\left \lvert h\right \rvert\leq \frac{R}{2}$  we have the following
    \begin{align}\label{diff to V diff appl}
        \int_{B_R}\left \lvert \uptau_h \left(d\widetilde{W}+\widetilde{F}\right)\right \rvert^p  \leq \delta \left \lvert h\right \rvert^{p-2} \int_{B_{R}}&\left \lvert \uptau_hV_p\left(d\widetilde{W}+\widetilde{F}\right)\right \rvert^2+c_\delta \left \lvert h \right \rvert^{p}\int_{B_{2R}} \left[\left \lvert d
        \widetilde{W} \right \rvert^p+\left \lvert  \widetilde{F} \right \rvert^p\right].
    \end{align}
\end{lemma}
\begin{proof}
Fix $B_{R}\subset \subset B_r\left(0\right)$ and using \eqref{V sim} with $\gamma=p$ we get the following:
    \begin{align*}
        &\int_{B_R}\left \lvert \uptau_h \left(d\widetilde{W}+\widetilde{F}\right)\right \rvert^p\notag\\
        &\quad \sim\int_{B_R}\left[\left(\left \lvert \left(d\widetilde{W}+\widetilde{F}\right) \left(x\right)\right \rvert+\left \lvert \left(d\widetilde{W}+\widetilde{F}\right)\left(x+h\right)\right \rvert\right)\right]^{\frac{2-p}{2}p}\left \lvert \uptau_h V_p\left(d\widetilde{W}+\widetilde{F}\right)\right \rvert^p.
    \end{align*}
Since $1<p<2$ by H\"{o}lder's Inequality with exponents $\left(\frac{2}{p},\frac{2}{2-p}\right)$, we have the following estimate,
\begin{align*}
&\int_{B_R}\left[\left \lvert \left(d\widetilde{W}+\widetilde{F}\right) \left(x\right)\right \rvert+\left \lvert \left(d\widetilde{W}+\widetilde{F}\right)\left(x+h\right)\right \rvert\right]^{\frac{2-p}{2}p}\left \lvert \uptau_h V_p\left(d\widetilde{W}+\widetilde{F}\right)\right \rvert^p \notag \\
&\quad= \begin{aligned}[t] \int_{B_R} &\left(h^{\frac{p-2}{2}p} \left \lvert \uptau_h V_p\left(d\widetilde{W}+\widetilde{F}\right)\right \rvert^p\notag\right) \left(\left \lvert h \right \rvert^{\frac{2-p}{2}p}\left[\left \lvert \left(d\widetilde{W}+\widetilde{F}\right) \left(x\right)\right \rvert+\left \lvert \left(d\widetilde{W}+\widetilde{F}\right)\left(x+h\right)\right \rvert\right]\right)  dx
\end{aligned}\notag\\
&\quad\leq\begin{aligned}[t] \delta &\left \lvert h\right \rvert^{p-2} \int_{B_{R}}\left \lvert \uptau_hV _p\left(\left(d\widetilde{W}+\widetilde{F}\right)\right)\right \rvert^2+c_\delta \left \lvert h \right \rvert^{p}\int_{B_R} \left(\left \lvert \left(d\widetilde{W}+\widetilde{F}\right) \left(x\right)\right \rvert+\left \lvert \left(d\widetilde{W}+\widetilde{F}\right)\left(x+h\right)\right \rvert\right)^{p}
        \end{aligned}\notag\\
        &\quad\lesssim\begin{aligned}[t] \delta \left \lvert h\right \rvert^{p-2} \int_{B_{R}}&\left \lvert \uptau_hV_p\left(d\widetilde{W}+\widetilde{F}\right)\right \rvert^2+c_\delta \left \lvert h \right \rvert^{p}\int_{B_{2R}} \left[\left \lvert d
        \widetilde{W} \right \rvert^p+\left \lvert  \widetilde{F} \right \rvert^p\right].
        \end{aligned}
    \end{align*}
    \end{proof}
\begin{lemma}
For $B_{10s}\subset \subset B_r\left(0\right)$ the following estimate holds:
\begin{align}\label{1<p<2 D_h Lemma}
\left \lbrace
\begin{aligned}
           \int_{B_{s}}&\left \lvert D_h V_p\left(d\widetilde{W}+\widetilde{F}\right)\right \rvert^2 \\
        &\leq \delta \int_{B_{4s}}\left \lvert D_hV_p\left(d\widetilde{W}+\widetilde{F}\right)\right \rvert^2\\
        &\qquad\qquad\qquad\qquad+c\left(\delta,p\right)s^{-2} \left[\left \lVert d\widetilde{W} \right \rVert_{L^{p}\left(B_{10s}\right)}^p+\left \lVert \widetilde{F} \right \rVert_{W^{2,p}\left(B_{10s}\right)}^p\right].
\end{aligned}
\right.
       \end{align}    
\end{lemma}
\begin{proof}
Fix a ball $B_{10s}\subset \subset B_r\left(0\right)$. Choose, non-negative smooth function $\eta$ with $\operatorname{supp}\left(\eta\right)\subseteq B_{2s}$ and $\eta \equiv 1$ in $B_{s}$.\\
In distributional  sense we have $d\left(d\widetilde{W}\right)=0$ in $B_{10s}$; so using Lemma \ref{1-div-curl-2} we have $\overline{W}\in W^{1,p}_{d^*,N}\left(B_{10s}\right)$ such that  $d\overline{{W}}=d\widetilde{W} \text{ in } B_{10s}$
Fix, $\left\lvert h\right \rvert\leq \frac{s}{1000}$ and taking test function $\phi=\uptau_{-h}\left(\eta^2\uptau_h \overline{W}\right)$ in \eqref{final weak formulation} we get the following
\begin{align*}
\int_{B_r\left(0\right)} &\left \langle \uptau_h\left[\left \langle G \left(d\widetilde{W}+\widetilde{F}\right), d\widetilde{W}+\widetilde{F}\right  \rangle^{\frac{p-2}{2}}G \left(d\widetilde{W}+\widetilde{F}\right)\right],\uptau_hd\overline{W}\right \rangle\eta^2\\
&\begin{aligned}[t]=-\int_{B_r\left(0\right)} \left \langle \uptau_h\left[\left \langle G \left(d\widetilde{W}+\widetilde{F}\right), d\widetilde{W}+\widetilde{F}\right  \rangle^{\frac{p-2}{2}}G \left(d\widetilde{W}+\widetilde{F}\right)\right],d \eta^2 \wedge\uptau_h \overline{W}\right \rangle.
\end{aligned}
\end{align*}
     Since, $d\overline{W}=d\widetilde{W}$ in $B_{10s}$ and $\operatorname{supp}\left(\eta\right)\subseteq B_{2s}$ we get the following
        \begin{align*}
        &\int_{B_r\left(0\right)} \left \langle \uptau_h\left[\left \langle G\left(d\widetilde{W}+\widetilde{F}\right), d\widetilde{W}+\widetilde{F}\right \rangle^{\frac{p-2}{2}}G \left(d\widetilde{W}+\widetilde{F}\right)\right],\uptau_hd\overline{W}\right \rangle\eta^2\\
        &=\begin{aligned}[t] &\int_{B_r\left(0\right)} \left \langle \uptau_h\left[\left \langle G\left(d\widetilde{W}+\widetilde{F}\right), d\widetilde{W}+\widetilde{F}\right \rangle^{\frac{p-2}{2}}G \left(d\widetilde{W}+\widetilde{F}\right)\right],\uptau_h\left(d\widetilde{W}+\widetilde{F}\right)\right \rangle\eta^2\\
        &-\int_{B_r\left(0\right)} \left \langle \uptau_h\left[\left \langle G\left(d\widetilde{W}+\widetilde{F}\right), d\widetilde{W}+\widetilde{F}\right \rangle^{\frac{p-2}{2}}G \left(d\widetilde{W}+\widetilde{F}\right)\right],\uptau_h\widetilde{F}\right \rangle\eta^2.
        \end{aligned}
    \end{align*}
    Define
    \begin{align*}
    &B_1:=\int_{B_r\left(0\right)} \left \langle \uptau_h\left[\left \langle G\left(d\widetilde{W}+\widetilde{F}\right), d\widetilde{W}+\widetilde{F}\right \rangle^{\frac{p-2}{2}}G \left(d\widetilde{W}+\widetilde{F}\right)\right],\uptau_h\left(d\widetilde{W}+\widetilde{F}\right)\right \rangle\eta^2,\\
    &B_2:=\int_{B_r\left(0\right)} \left \langle \uptau_h\left[\left \langle G\left(d\widetilde{W}+\widetilde{F}\right), d\widetilde{W}+\widetilde{F}\right \rangle^{\frac{p-2}{2}}G \left(d\widetilde{W}+\widetilde{F}\right)\right],\uptau_h\widetilde{F}\right \rangle\eta^2,\\
        &B_3:= \int_{B_r\left(0\right)} \left \langle \uptau_h\left[\left \langle G\left(d\widetilde{W}+\widetilde{F}\right), d\widetilde{W}+\widetilde{F}\right \rangle^{\frac{p-2}{2}}G \left(d\widetilde{W}+\widetilde{F}\right)\right],\uptau_hd\overline{W}\right \rangle\eta^2.
    \end{align*}
So, $B_3=B_1-B_2$ which implies that $B_1\leq \left \lvert B_2\right \rvert+\left \lvert B_3\right \rvert$.
By \eqref{c1p} we have
    \begin{align}\label{B_1}
         c_1\left(p\right)\int_{B_r\left(0\right)}\left \lvert \uptau_h V_p\left(d\widetilde{W}+\widetilde{F}\right)\right \rvert^2 \eta^2 \leq B_1.
    \end{align}
    We have assumed that $\widetilde{F}\in W^{2,\widetilde{q}}\left({B_r\left(0\right)}\right)$, on the other hand $\eta$ is smooth with $\operatorname{supp}\left(\eta\right)\subset B_s\subset B_{10s}\subset \subset B_r\left(0\right)$ thus
    $\eta^2\uptau_h\widetilde{F}\in W_0^{1,p}\left(B_r\left(0\right)\right)$. Also 
    \begin{align*}
        \left \lvert \left \langle G\left(d\widetilde{W}+\widetilde{F}\right), d\widetilde{W}+\widetilde{F}\right \rangle^{\frac{p-2}{2}}G_{\Phi}\left(d\widetilde{W}+\widetilde{F}\right)\right \rvert\lesssim \left \lvert d\widetilde{W}+\widetilde{F}\right \rvert^{p-1},
    \end{align*}
    which implies $\left \langle G\left(d\widetilde{W}+\widetilde{F}\right), d\widetilde{W}+\widetilde{F}\right \rangle^{\frac{p-2}{2}}G_{\Phi}\left(d\widetilde{W}+\widetilde{F}\right) \in L^{\frac{p}{p-1}}\left(B_r\left(0\right)\right)$. On the otherhand  $\eta^2 \uptau_h \widetilde{F}\in L^p\left(B_r\left(0\right)\right)$. Hence using \eqref{<tau h G, H>} from Lemma \ref{Lemma <G,H>} we get,
    \begin{align*}
        \left \lvert B_2 \right \rvert
        &\leq c\left(s\right)\left \lvert h\right \rvert \int_{B_{2s}} \left[\int_0^1\left \lvert \left(d\widetilde{W}+\widetilde{F}\right)\left(x+th\right)\right \rvert^{{p-1}}dt\right]\left(\left\lvert \uptau_h \nabla \widetilde{F}\right\rvert+\left\lvert \uptau_h  \widetilde{F}\right\rvert\right)dx\\
        &\leq c\begin{aligned}[t]&\left(s\right)\left \lvert h\right \rvert \left(\int_{B_{2s}} \left[\int_0^1\left \lvert \left(d\widetilde{W}+\widetilde{F}\right)\left(x+th\right)\right \rvert^{{p-1}}dt\right]^{\frac{p}{p-1}}dx\right)^{\frac{p-1}{p}}\\
        &\qquad\qquad\times\left(\int_{B_{2s}}\left(\left\lvert \uptau_h \nabla \widetilde{F}\right\rvert+\left\lvert \uptau_h  \widetilde{F}\right\rvert\right)^pdx\right)^{\frac{1}{p}}\left[\text{by H\"{o}lder's Inequality}\right]
        \end{aligned}\\
        &\leq c\begin{aligned}[t]&\left(s\right)\left \lvert h\right \rvert \left(\int_{B_{2s}} \left[\int_0^1\left \lvert \left(d\widetilde{W}+\widetilde{F}\right)\left(x+th\right)\right \rvert^{{p}}dt\right]dx\right)^{\frac{p-1}{p}}\\
        &\times\left(\int_{B_{2s}}\left\lvert \uptau_h \nabla \widetilde{F}\right\rvert^p+\int_{B_{2s}}\left\lvert \uptau_h  \widetilde{F}\right\rvert^p\right)^{\frac{1}{p}}\left[\text{by Jensen's Inequality for $t\mapsto t^{\frac{p}{p-1}}$}\right].
        \end{aligned}
        \end{align*}
        By \eqref{differencequotient} we get the following
        \begin{align*}
        \int_{B_{2s}}\left\lvert \uptau_h \nabla \widetilde{F}\right\rvert^p+\int_{B_{2s}}\left\lvert \uptau_h  \widetilde{F}\right\rvert^p
        &\quad\leq \begin{aligned}[t]&\left \lvert h\right \rvert^p\int_{B_{3s}}\left(\left\lvert D^2 \widetilde{F}\right\rvert^p+\left\lvert \nabla \widetilde{F}\right\rvert^p\right),
        \end{aligned}
    \end{align*}
    and by Fubini's theorem
    \begin{align*}
        \int_{B_{2s}} \left[\int_0^1\left \lvert \left(d\widetilde{W}+\widetilde{F}\right)\left(x+th\right)\right \rvert^{{p}}dt\right]dx&\leq  \int_0^1\left[\int_{B_{2s}}\left \lvert \left(d\widetilde{W}+\widetilde{F}\right)\left(x+th\right)\right \rvert^{{p}}dx\right]dt\\
        &\leq \int_{B_{3s}}\left \lvert d\widetilde{W}+\widetilde{F}\right \rvert^{{p}}.
    \end{align*}
    Combining we get that
    \begin{align*}
        &\left \lvert B_2 \right \rvert \leq \begin{aligned}[t]&c\left(s\right)\left \lvert h\right \rvert^2\left(\int_{B_{3s}}\left \lvert d\widetilde{W}+\widetilde{F}\right \rvert^{{p}}\right)^{\frac{p-1}{p}}\left(\int_{B_{3s}}\left(\left\lvert D^2 \widetilde{F}\right\rvert^p+\left\lvert \nabla \widetilde{F}\right\rvert^p\right)\right)^{\frac{1}{p}}.
        \end{aligned}
    \end{align*}
    Using Young's Inequality with exponents $\left(\frac{p}{p-1},p\right)$ we get
    \begin{align}\label{B_2}
        \left \lvert B_2\right \rvert
        &\leq c\begin{aligned}[t]&\left(s\right)\left \lvert h\right \rvert^2  \int_{B_{2s}}\left(\left \lvert \left(d\widetilde{W}+\widetilde{F}\right)\right \rvert^{{p}}+\left\lvert D^2 \widetilde{F}\right\rvert^p+\left\lvert \nabla \widetilde{F}\right\rvert^p\right)
        \end{aligned}\notag \\
        &\leq c\begin{aligned}[t]&\left(s\right)\left \lvert h\right \rvert^2  \left(\int_{B_{10s}}\left \lvert d\overline{W}\right \rvert^p+\left\lVert  \widetilde{F}\right\rVert^p_{W^{2,p}\left(B_{10s}\right)}\right)\left[\text{since, $d\widetilde{W}=d\overline{W}$ in $B_{10s}$}\right].
        \end{aligned}
    \end{align}
    Similarly using Lemma \ref{Lemma <G,H>} we have
        \begin{align*}
        \left \lvert B_3\right \rvert
        &\leq c\left(s\right)\left \lvert h\right \rvert \int_{B_{2s}} \int_0^1\left \lvert \left(d\widetilde{W}+\widetilde{F}\right)\left(x+th\right)\right \rvert^{{p-1}}\left(\left\lvert \uptau_h \overline{W}\right\rvert+\left\lvert \uptau_h \nabla \overline{W}\right\rvert\right).
    \end{align*}
    We define
    \begin{align*}
        &J_1:=c\left(s\right)\left \lvert h\right \rvert \int_{B_{2s}} \left[\int_0^1\left(\left \lvert \left(d\widetilde{W}+\widetilde{F}\right)\left(x+th\right)\right \rvert\right)^{p-1}dt\right]\left\lvert \uptau_h \overline{W}\right\rvert dx,\\
        &J_2:= c\left(s\right)\left \lvert h\right \rvert \int_{B_{2s}} \left[\int_0^1\left(\left \lvert \left(d\widetilde{W}+\widetilde{F}\right)\left(x+th\right)\right \rvert\right)^{{p-1}}dt\right]\left\lvert \uptau_h \nabla \overline{W}\right\rvert dx;
    \end{align*}
hence, $\left \lvert B_3\right \rvert\leq J_1+J_2$.\\
Using H\"{o}lder's inequality with exponents $\left(\frac{p}{p-1},p\right)$, first w.r.t $dx $ and then w.r.t. $dt$, Fubini's theorem and \eqref{differencequotient} we derive the following estimate:
\begin{align}\label{1<p<2 I1 estimae}
J_1&\leq  c\left(s\right)\begin{aligned}[t]&\left \lvert h\right \rvert \left[\int_{B_{2s}} \left(\int_0^1\left \lvert \left(d\widetilde{W}+\widetilde{F}\right)\left(x+th\right)\right \rvert^{{p-1}}dt\right)^{\frac{p}{p-1}}dx\right]^{\frac{p-1}{p}}\left[\int_{B_{2s}} \left\lvert \uptau_h \overline{W}\right\rvert^p dx\right]^{\frac{1}{p}}
\end{aligned}\notag \\
&\begin{aligned}[t] \leq c\left(s\right)&\left \lvert h\right \rvert \left[\int_{B_{2s}} \left(\int_0^1\left \lvert \left(d\widetilde{W}+\widetilde{F}\right)\left(x+th\right)\right \rvert^{{p}}dt\right)dx\right]^{\frac{p-1}{p}}\left[\int_{B_{2s}} \left\lvert \uptau_h \overline{W}\right\rvert^p dx\right]^{\frac{1}{p}}\end{aligned}\notag \\
&\begin{aligned}[t] \leq c\left(s\right)&\left \lvert h\right \rvert \left[\int_0^1\left(\int_{B_{2s}}\left \lvert \left(d\widetilde{W}+\widetilde{F}\right)\left(x+th\right)\right \rvert^{{p}}dx\right)dt\right]^{\frac{p-1}{p}}\left[\left \lvert h\right \rvert^p\int_{B_{3s}} \left\lvert \nabla \overline{W}\right\rvert^p dx\right]^{\frac{1}{p}}
        \end{aligned}\notag \\
        &\begin{aligned}[t] \leq c\left(s\right)&\left \lvert h\right \rvert^2 \left[\int_{B_{3s}}\left \lvert d\widetilde{W}+\widetilde{F}\right \rvert^{{p}}\right]^{\frac{p-1}{p}}\left[\int_{B_{3s}} \left\lvert \nabla \overline{W}\right\rvert^p \right]^{\frac{1}{p}}
        \end{aligned} \leq c\left(s\right) \left \lvert h\right \rvert^2\int_{3s}\left[\left \lvert \nabla \overline{W}\right\rvert^{p}+\left \lvert \widetilde{F}\right\rvert^{p}\right].
    \end{align}
    In the last line while estimating $J_1$ we used the fact that $\left \lvert d\overline{W}\right \rvert\leq c\left(n,k\right)\left \lvert \nabla \overline{W}\right \rvert$ point-wise in $B_{10s}$
    and in a similar way we get the following:
    \begin{align}\label{I2 inequality}
    \left \lbrace
    \begin{aligned}
        J_2
        &\leq \begin{aligned}[t] c\left(s\right)&\left \lvert h\right \rvert \left[\int_{B_{3s}}\left[\left \lvert d\widetilde{W}\right\rvert^{p}+\left \lvert \widetilde{F}\right\rvert^{p}\right]dx\right]^{\frac{p-1}{p}}\left[\int_{B_{2s}} \left\lvert  \uptau_h \nabla \overline{W}\right\rvert^p dx\right]^{\frac{1}{p}}\\
        \end{aligned}\\
        &=\begin{aligned}[t]  c\left(s\right)\left \lvert h\right \rvert \left[\int_{B_{3s}}\left[\left \lvert d\overline{W}\right\rvert^{p}+\left \lvert \widetilde{F}\right\rvert^{p}\right]dx\right]^{\frac{p-1}{p}}&\left[\int_{B_{2s}} \left\lvert  \uptau_h \nabla \overline{W}\right\rvert^p dx\right]^{\frac{1}{p}}\\
        &\text{[since, $d\widetilde{W}=d\overline{W}$ in $B_{10s}$]}
        \end{aligned}\\
        &\leq \begin{aligned}[t] c\left(s\right)&\left[\left \lvert h\right \rvert^2\int_{B_{3s}}\left[\left \lvert d\overline{W}\right\rvert^{p}+\left \lvert \widetilde{F}\right\rvert^{p}\right]dx\right]^{\frac{p-1}{p}}\times \left \lvert h\right \rvert ^{\frac{2-p}{p}}\left[\int_{B_{2s}} \left\lvert  \uptau_h \nabla \overline{W}\right\rvert^p dx\right]^{\frac{1}{p}}
        \end{aligned}\\
        &\leq \begin{aligned}[t] c\left(s\right)&\left \lvert h\right \rvert^2\int_{B_{3s}} \left[\left \lvert d\overline{W}\right\rvert^{p}+\left \lvert \widetilde{F}\right\rvert^{p}\right]+ c\left(s\right)\left \lvert h\right \rvert ^{2-p}\int_{B_{2s}} \left\lvert  \uptau_h \nabla \overline{W}\right\rvert^p .
        \end{aligned}
        \end{aligned}
        \right.
    \end{align}In the last estimate we used Young's inequality with exponents $\left(\frac{p}{p-1},p\right)$.
    Choose, smooth non-negative function $\xi$ such that $\xi\equiv 1$ in $B_{2s}$ and $\operatorname{supp}\left(\eta\right)\subseteq B_{3s}$. Now, 
    \begin{align*}
    \uptau_hd\left(\xi \overline{W}\right)
    &=\left(\uptau_hd\xi\right)\wedge \overline{W}+d\xi\left(x+h\right)\wedge \uptau_h \overline{W}+\left(\uptau_h\xi \right) d\overline{W}+\xi \left(x+h\right)\left(\uptau_h d\overline{W}\right).
    \end{align*}
    Hence,
        $\left \lvert \uptau_hd\left(\xi \overline{W}\right)\right \rvert
    \leq c\left(s\right)\left \lvert h\right \rvert \left[\left \lvert \overline{W} \right \rvert+\left \lvert d\overline{W} \right \rvert\right]+c\left(s\right)\left \lvert  \uptau_h \overline{W}\right \rvert+\left \lvert \uptau_h d\overline{W}\right \rvert$.
    Thus using \eqref{differencequotient} we get the following
    \begin{align}\label{d inequality}
        \int_{B_{4s}}&\left \lvert \uptau_hd\left(\xi \overline{W}\right)\right \rvert^p\leq c\left(s\right)\left \lvert h\right \rvert^p\int_{B_{10s}}\left[\left \lvert \overline{W} \right \rvert^p+\left \lvert \nabla \overline{W}\right \rvert^p\right]+c\left(p\right)\int_{B_{4s}}\left \lvert \uptau_h d\overline{W} \right \rvert^p.
    \end{align}
Similarly using the fact that $d^*\overline{W}=0$ in $B_{10s}$ we get, $\left \lvert \uptau_h d^*\left(\xi \overline{W}\right)\right \rvert \leq c\left(s\right)\left \lvert h\right \rvert\left \lvert \overline{W} \right \rvert+c\left(s\right)\left \lvert  \uptau_h \overline{W}\right \rvert$ and using \eqref{differencequotient} again we arrive at the following
    \begin{align}\label{d* inequality}
        \int_{B_{4s}}\left \lvert \uptau_hd^*\left(\xi \overline{W}\right)\right \rvert^p
        \leq c\left(s\right)\left \lvert h\right \rvert^p\int_{B_{10s}}\left[\left \lvert \overline{W} \right \rvert^p+\left \lvert \nabla \overline{W}\right \rvert^p\right].
    \end{align}
    Now $\operatorname{supp}\left(\xi\right)\subset\subset B_{3s}$  implies  $\uptau_h\left(\xi\overline{W}\right) \in W^{1,p}_T\left(B_{4s}\right)$ for $\left \lvert h\right \rvert\leq 10^{-3}s$, hence using Theorem \ref{nabla d d*} we get the following estimate
    \begin{align*}
        \int_{B_{2s}}\left \lvert \uptau_h\nabla \overline{W} \right \rvert^p
        &\leq \int_{B_{4s}}\left \lvert \uptau_h\nabla\left( \xi \overline{W}\right) \right \rvert^p\left[\text{Since, $\xi\equiv 1$ in $B_{2s}$}\right]\\
        &\leq c\left(s\right)\begin{aligned}[t]\int_{B_{4s}}&\left \lvert \uptau_h d\left(\xi \overline{W}\right) \right \rvert^p+c\left(s\right)\int_{B_{4s}}\left \lvert \uptau_h d^*\left(\xi \overline{W}\right) \right \rvert^p
        \end{aligned}\\
        &\leq c\left(p\right)\begin{aligned}[t] \int_{B_{4s}} &\left \lvert \uptau_h d\widetilde{W} \right \rvert^p+c\left(s\right)\left\lvert h\right\rvert^p \int_{B_{10s}}\left[\left \lvert \overline{W} \right \rvert^p+\left \lvert \nabla \overline{W}\right \rvert^p\right]\\
        &\left[\text{By \eqref{d inequality} and \eqref{d* inequality} with the fact, $d\overline{W}=d\widetilde{W}$ in $B_{10s}$}\right]
        \end{aligned}\\
        &\leq c\left(p\right)\begin{aligned}[t] &\int_{B_{4s}} \left \lvert \uptau_h \left(d\widetilde{W}+\widetilde{F}\right) \right \rvert^p+c\left(s\right)\left\lvert h\right\rvert^p \int_{B_{10s}}\left[\left \lvert \overline{W} \right \rvert^p+\left \lvert \nabla \widetilde{F} \right \rvert^p+\left \lvert \nabla \overline{W}\right \rvert^p\right]\\
        &\left[\text{By triangle inequality and \eqref{differencequotient} for $\widetilde{F}$}\right]
        \end{aligned}\\
        &\leq \begin{aligned}[t]&\delta \left \lvert h\right \rvert^{p-2} \int_{B_{4s}}\left \lvert \uptau_hV_p\left(d\widetilde{W}+\widetilde{F}\right)\right \rvert^2\\
        &+c\left(s,\delta,p\right) \left \lvert h \right \rvert^{p}\left[\left \lVert \overline{W} \right \rVert_{W^{1,p}\left(B_{10s}\right)}^p+\left \lVert \widetilde{F} \right \rVert_{W^{1,p}\left(B_{10s}\right)}^p\right]\left[\text{by \eqref{diff to V diff appl} for $R=4s$}\right].
        \end{aligned}
    \end{align*}
    Combining the last estimate with \eqref{I2 inequality} we arrive at the following:
    \begin{align}\label{J_2}
        J_2
        &\leq \delta \int_{B_{4s}}\left \lvert \uptau_hV_p\left(d\widetilde{W}+\widetilde{F}\right)\right \rvert^2+c\left(s,\delta,p\right) \left \lvert h \right \rvert^{2}\left[\left \lVert \overline{W} \right \rVert_{W^{1,p}\left(B_{10s}\right)}^p+\left \lVert \widetilde{F} \right \rVert_{W^{1,p}\left(B_{10s}\right)}^p\right].
    \end{align}
    We  combine \eqref{B_1},\eqref{B_2}, \eqref{1<p<2 I1 estimae} and \eqref{J_2} with the fact that, $B_1\leq \left \lvert B_2\right \rvert+\left \lvert B_3\right \rvert\leq \left \lvert B_2\right \rvert+J_1+J_2$ and we get that, if $\left \lvert h\right \rvert\leq 10^{-3}s$ then
    \begin{align*}
        &\int_{B_{s}}\left \lvert \uptau_h V_p\left(d\widetilde{W}+\widetilde{F}\right)\right \rvert^2 \\
        &\quad\leq \delta \int_{B_{4s}}\left \lvert \uptau_hV_p\left(d\widetilde{W}+\widetilde{F}\right)\right \rvert^2+c\left(s,\delta,p\right) \left \lvert h \right \rvert^{2}\left[\left \lVert \overline{W} \right \rVert_{W^{1,p}\left(B_{10s}\right)}^p+\left \lVert \widetilde{F} \right \rVert_{W^{2,p}\left(B_{10s}\right)}^p\right],
    \end{align*}
    i.e. for $\left \lvert h\right \rvert\leq \frac{s}{1000}$ we have
    \begin{align*}
        \int_{B_{2s}}&\left \lvert D_h V_p\left(d\widetilde{W}+\widetilde{F}\right)\right \rvert^2 \eta^2\\
        &\leq \delta \int_{B_{4s}}\left \lvert D_hV_p\left(d\widetilde{W}+\widetilde{F}\right)\right \rvert^2+c\left(s,\delta,p\right) \left[\left \lVert \overline{W} \right \rVert_{W^{1,p}\left(B_{10s}\right)}^p+\left \lVert \widetilde{F} \right \rVert_{W^{2,p}\left(B_{10s}\right)}^p\right].
    \end{align*}
    But, $\overline{W}\in W^{1,p}_{d^*,N}\left(B_{10s}\right)$ thus using Theorem \ref{nabla d d*} and the fact $d\overline{W}=d\widetilde{W}$  in $B_{10s}$ we get
    \begin{align*}
        &\int_{B_{s}}\left \lvert D_h V_p\left(d\widetilde{W}+\widetilde{F}\right)\right \rvert^2 \notag \\
        &\quad\leq \delta \int_{B_{4s}}\left \lvert D_hV_p\left(d\widetilde{W}+\widetilde{F}\right)\right \rvert^2+c\left(s,\delta,p\right) \left[\left \lVert d\overline{W} \right \rVert_{L^{p}\left(B_{10s}\right)}^p+\left \lVert \widetilde{F} \right \rVert_{W^{2,p}\left(B_{10s}\right)}^p\right]\notag \\
        &\quad\leq \delta \int_{B_{4s}}\left \lvert D_hV_p\left(d\widetilde{W}+\widetilde{F}\right)\right \rvert^2+c\left(s,\delta,p\right) \left[\left \lVert d\widetilde{W} \right \rVert_{L^{p}\left(B_{10s}\right)}^p+\left \lVert \widetilde{F} \right \rVert_{W^{2,p}\left(B_{10s}\right)}^p\right].
    \end{align*}
       Using standard scaling argument we get the estimate \eqref{1<p<2 D_h Lemma}.
       \end{proof}
Finally we are in a position to prove that for $1<p<2$, $d\widetilde{W}+\widetilde{F}\in L^{\widetilde{q}}_{\operatorname{loc}}B_r\left(0\right)$ where $\widetilde{q}$ satisfies \eqref{qtilde} for which take the help of following Lemma:
\begin{lemma}\label{h_j Lemma}
Fix $i\in \left \{1,2,3,\cdots,n\right\}$ and $B_{10r_1}\subset \subset B_r\left(0\right)$. For all $j\in \mathbb{N}$, define $h_j=10^{-6}\frac{r_1}{2^{j+1}}e_i$ where $e_i$ is the $i$-th basis vector for $\mathbb{R}^n$, then 
\begin{align}\label{A0 and difference quotient}
        \int_{B_{r_1}}\left \lvert D_{h_j}V_p\left(d\widetilde{W}+\widetilde{F}\right)\right \rvert^2 \leq A_0;
    \end{align}    
where $A_0=c\left(n,p,r_1\right)\left[\left \lVert \widetilde{W} \right \rVert^p_{W^{1,p}\left(B_{r_1}\right)}+\left \lVert \widetilde{F} \right \rVert^p_{W^{2,p}\left(B_{r_1}\right)}\right]$.
\end{lemma}
\begin{proof}  
Choose, $r_1\leq t<T\leq 2r_1$ where $B_{10r_1}\subset\subset B_{r}\left(0\right)$.
Take, $s=\frac{T-t}{1000}$ and we can cover $B_{t}$ by finitely many balls with centers say $y_1,y_2,...y_{M}$ in $B_{r_1}$ and each radius $s$ such that each ball intersects at most $c\left(n\right)$ (independent of $s$) many balls from this family of $M$ balls. Thus, for $\left \lvert h\right\rvert\leq \frac{s}{10^3} $ we use  \eqref{1<p<2 D_h Lemma} and get the following estimate:
\begin{align*}
&\int_{B_{t}}\left \lvert D_h V_p\left(d\widetilde{W}+\widetilde{F}\right)\right \rvert^2 \eta^2\\
&\quad\leq \sum_{j=1}^{M}\int_{B_{s}\left(y_j\right)}\left \lvert D_h V_p\left(d\widetilde{W}+\widetilde{F}\right)\right \rvert^2 \eta^2\\
&\quad\leq \delta \begin{aligned}[t]\sum_{j=1}^{M}&\int_{B_{4s}\left(y_j\right)}\left \lvert D_hV_p\left(d\widetilde{W}+\widetilde{F}\right)\right \rvert^2+c\left(p,\delta\right)s^{-2}\sum_{j=1}^{M} \left[\left \lVert d\widetilde{W} \right \rVert_{L^{p}\left(B_{10s}\left(y_j\right)\right)}^p+\left \lVert \widetilde{F} \right \rVert_{W^{2,p}\left(B_{10s}\left(y_j\right)\right)}^p\right]
\end{aligned}\\
&\quad\leq \delta \begin{aligned}[t]c^*\left(n\right)&\int_{B_{T}}\left \lvert D_hV_p\left(d\widetilde{W}+\widetilde{F}\right)\right \rvert^2+c\left(n,p,\delta\right)\left(T-t\right)^{-2} \left[\left \lVert d\widetilde{W} \right \rVert_{L^{p}\left(B_{T}\right)}^p+\left \lVert \widetilde{F} \right \rVert_{W^{2,p}\left(B_{T}\right)}^p\right].
\end{aligned}
\end{align*}
    Taking $\delta=\frac{1}{8c^*n}$
    \begin{align*}
        \int_{B_{t}}&\left \lvert D_h V_p\left(d\widetilde{W}+\widetilde{F}\right)\right \rvert^2 \eta^2\\
        &\leq \frac{1}{8} \begin{aligned}[t]&\int_{B_{T}}\left \lvert D_hV_p\left(d\widetilde{W}+\widetilde{F}\right)\right \rvert^2+c\left(n,p\right)\left(T-t\right)^{-2}\left[\left \lVert d\widetilde{W} \right \rVert_{L^{p}\left(B_{T}\right)}^p+\left \lVert \widetilde{F} \right \rVert_{W^{2,p}\left(B_{T}\right)}^p\right].
        \end{aligned}
    \end{align*}
    For $r_1\leq t<T\leq 2r_1$ and $\left \lvert h\right \rvert\leq \frac{T-t}{10^6}$ 
    \begin{align}\label{iteration set up}
        \int_{B_t}\left \lvert D_hV_p\left(d\widetilde{W}+\widetilde{F}\right)\right \rvert^2\leq \frac{1}{8} \int_{B_T}\left \lvert D_hV_p\left(d\widetilde{W}+\widetilde{F}\right)\right \rvert^2+\frac{A}{\left(T-t\right)^2}
    \end{align}
    where $A=c\left(n,p\right)\left[\left \lVert \widetilde{W} \right \rVert^p_{W^{1,p}\left(B_{T}\right)}+\left \lVert \widetilde{F} \right \rVert^p_{W^{2,p}\left(B_{T}\right)}\right]$.\\
    Fix, $i\in \left\{1,2,\cdots,n\right\}$; choose arbitrary $j\in \mathbb{N}$ and define \begin{align*}\delta_j=\frac{r_1}{2^{j+1}},t_0=r_1,t_{j+1}=t_j+\delta_j, h_j=10^{-6}\delta_je_i;
    \end{align*}thus for $l=0,1,2,\cdots,j$ 
    \begin{align*}
        \left \lvert h_j\right \rvert=\frac{10^{-6}r_1}{2^{j+1}}\leq \frac{10^{-6}r_1}{2^{l+1}}=10^{-6}\delta_l=\frac{t_{l+1}-t_l}{10^6}.
    \end{align*}
    Now
    \begin{align*}
        \int_{B_{r_1}}&\left \lvert D_{h_j}V_p\left(d\widetilde{W}+\widetilde{F}\right)\right \rvert^2\\&\leq \frac{1}{8} \int_{B_{t_1}}\left \lvert D_{h_j}V_p\left(d\widetilde{W}+\widetilde{F}\right)\right \rvert^2+\delta_{0}^{-2}A\text{ [Taking $h=h_j,T=t_1,t=t_0$ in \eqref{iteration set up}]}\\
        &\leq \begin{aligned}[t]\frac{1}{8^2} \int_{B_{t_2}}\left \lvert D_{h_j}V_p\left(d\widetilde{W}+\widetilde{F}\right)\right \rvert^2&+A\left(\delta_0^{-2}+\frac{1}{8}\delta_1^{-2}\right)\\
        &\text{ [We chose $h=h_j,T=t_2,t=t_1$ in \eqref{iteration set up}]}
        \end{aligned}\\
        &\leq \begin{aligned}[t]\frac{1}{8^3} \int_{B_{t_3}}\left \lvert D_{h_j}V_p\left(d\widetilde{W}+\widetilde{F}\right)\right \rvert^2&+A\left(\delta_0^{-2}+\frac{1}{8}\delta_1^{-2}+\frac{1}{8^2}\delta_2^{-2}\right)\\
        &\text{ [Putting $h=h_j,T=t_1,t=t_0$ in \eqref{iteration set up}]}
        \end{aligned}\\
        &\qquad \qquad\qquad\qquad\qquad\qquad\vdots\\
        &\leq \frac{1}{8^j} \int_{B_{t_j}}\left \lvert D_{h_j}V_p\left(d\widetilde{W}+\widetilde{F}\right)\right \rvert^2+A\sum_{l=0}^j \left(\frac{1}{8^l\delta_l^2}\right).\\
    \end{align*}
    Now, 
    \begin{align*}
        \left \lvert D_{h_j}V_p\left(d\widetilde{W}+\widetilde{F}\right)\right \rvert^2
        &=\frac{\left \lvert V_p\left(d\widetilde{W}+\widetilde{F}\right)\left(x+h_j\right)-V_p\left(d\widetilde{W}+\widetilde{F}\right)\left(x\right)\right \rvert^2}{\left\lvert h_j\right \rvert^2}\\
        &= 10^{12}\times 2^{2j+2}\frac{\left \lvert V_p\left(d\widetilde{W}+\widetilde{F}\right)\left(x+h_j\right)+ V_p\left(d\widetilde{W}+\widetilde{F}\right)\left(x\right)\right \rvert^2}{r_1^2}\\
        &= 4 \times 10^{12}\times 4^{j}\frac{\left \lvert V_p\left(d\widetilde{W}+\widetilde{F}\right)\left(x+h_j\right)\right \rvert^2+\left \lvert V_p\left(d\widetilde{W}+\widetilde{F}\right)\left(x\right)\right \rvert^2}{r_1^2}\\
        &= 4 \times 10^{12}\times 4^{j}\frac{\left\lvert \left(d\widetilde{W}+\widetilde{F}\right)\left(x+h_j\right)\right \rvert^p+\left \lvert \left(d\widetilde{W}+\widetilde{F}\right)\left(x\right)\right \rvert^p}{r_1^2}.
    \end{align*}
    On the other hand, for $l=1,2,\cdots,j$ we have 
        $8^l\delta_l^2=\frac{1}{4}{r_1^2} \times 2^l$.\\
    Define $A_0=c\left(n,p,r_1\right)\left[\left \lVert \widetilde{W} \right \rVert^p_{W^{1,p}\left(B_{r_1}\right)}+\left \lVert \widetilde{F} \right \rVert^p_{W^{2,p}\left(B_{r_1}\right)}\right]$ and thus for $i=1,2,\cdots n$ and $h_j=10^{-6}\frac{r_1}{2^{j+1}}e_i$, we have \eqref{A0 and difference quotient}.
    \end{proof}
By Lemma \ref{h_j Lemma}, we have the estimate in \eqref{A0 and difference quotient} where $h_j=10^{-6}\frac{r_1}{2^{j+1}}e_i\rightarrow 0$; on the other hand
    Also
    \begin{align*}
        \int_{B_{r_1}} \left \lvert V_p\left(d\widetilde{W}+\widetilde{F}\right)\right \rvert^2= \int_{B_{r_1}} \left \lvert d\widetilde{W}+\widetilde{F}\right \rvert^p.
    \end{align*}
    Thus by \eqref{A0 and difference quotient} and Banach-Alaoglu Theorem $V_p\left(d\widetilde{W}+\widetilde{F}\right) $ is in $ W^{1,2}\left(B_{r_1};\Lambda^{k+1}\right)$.\\ Finally using Sobolev Embedding we get the following cases (see \cite[Corollary~$9.14$]{brezis2010functional}):\\
    If $n=2$, then for all $\gamma\in \left[2,\infty\right)$ we have
        $V_p\left(d\widetilde{W}+\widetilde{F}\right)\in L^{\gamma}\left(B_{r_1}\right)$. On the other hand
    if $n\geq 3$
    then for all $\gamma\in \left[2,\frac{2n}{n-2}\right]$,
     $V_p\left(d\widetilde{W}+\widetilde{F}\right)\in L^{\gamma}\left(B_{r_1}\right)$. 
    Also
    \begin{align*}
        \int_{B_{r_1}}\left \lvert V_p\left(d\widetilde{W}+\widetilde{F}\right)\right \rvert^{\gamma}&= \int_{B_{r_1}}\left \lvert V_p\left(d\widetilde{W}+\widetilde{F}\right)\right \rvert^{2.\frac{\gamma}{2}}=\int_{B_{r_1}}\left \lvert d\widetilde{W}+\widetilde{F}\right \rvert^{\frac{\gamma p}{2}}.
    \end{align*}
    As a consequence $\left(d\widetilde{W}+\widetilde{F}\right)\in L^{{\tilde{q}}}\left(B_{r_1}\right)$ where for $n=2$, $\tilde{q}\in \left[p,\infty\right)$ and  for $n\geq 3$, $\tilde{q}$ is any number in $\left[p,\frac{np}{n-2}\right]$. Thus for both of the cases, $2\leq p$ and $1<p<2$ Proposition \ref{Proposition with F in W2,p} is proved.
    \end{proof}
        
Proof of Proposition \ref{Application of HL FS KU Proposition}:
\begin{proof}
    For $B_R\subset \subset B_{\frac{r}{4}}$
we define ${v}_R$ as the minimizer of the following functional
\begin{align*}
    \overline{w}\mapsto \int_{B_R}\left \lvert d\overline{w}\right \rvert^p; \overline{w}\in \widetilde{W}+W^{1,p}_{d^*,T}\left(B_R;\Lambda^k\right)
\end{align*}
By definition
\begin{align}\label{in particular minimization}
    \int_{B_R}\left \lvert dv_R\right \rvert^p\leq \int_{B_R}\left \lvert d\widetilde{W}\right \rvert^p.
\end{align}
and we have
\begin{align}\label{vR el}
        \int_{B_R}\left\langle \left \lvert d{v}_R\right \rvert^{p-2}d{v}_R,d\phi\right \rangle=0, \text{ for $\phi \in W^{1,p}_{d^*,T}\left(B_R;\Lambda^k\right)$}.
        \end{align}
        On the other hand by K.Ulhenbeck's Theorem for $B_\rho \subset \subset B_R$ we have the following decay estimate (see \ref{decay estimate}):
        \begin{align}\label{Ulhenbeck Application}
            \fint_{B_\rho}& \left \lvert V_p\left(dv_R\right)-\left(V_p\left(dv_R\right)\right)_{B_\rho} \right \rvert^2 \leq c_U\left(\frac{\rho}{R}\right)^{2\alpha_1}\fint_{B_R} \left \lvert V_p\left(dv_R\right)-\left(V_p\left(dv_R\right)\right)_{B_R} \right \rvert^2.
        \end{align}
One weak formulation of \eqref{2nd Proposition eqn. mention} (also see \eqref{G weak formulation in theorem} or \eqref{el fo W}) is \eqref{d* T weal formulation} i.e. for all $\chi\in W^{1,p}_{d^*,T}\left(B_{r}\left(0\right)\right)$ we have the following
    \begin{align}\label{d* T weal formulation in proof of proposition}
        \int_{B_r\left(0\right)}\left \langle \left \langle G \left(d\widetilde{W}+\widetilde{F}\right), d\widetilde{W}+\widetilde{F}\right  \rangle^{\frac{p-2}{2}}G \left(d \widetilde{W}+\widetilde{F}\right),d\chi\right \rangle=0.
    \end{align}
    Hence by \eqref{vR el}, \eqref{d* T weal formulation in proof of proposition} along with the fact that $\widetilde{W}-v_R\in W^{1,p}_{d^*,T}\left(B_R;\Lambda^k\right)$ we get
    \begin{align*}
        0=\int_{B_R} &\left \langle G \left(d\widetilde{W}+\widetilde{F}\right),d\widetilde{W}+\widetilde{F}\right  \rangle^{\frac{p-2}{2}}\\ 
        &\qquad \qquad \times\left \langle  G \left(d\widetilde{W}+\widetilde{F}\right)- \left \lvert dv_R\right \rvert^{p-2}dv_R, d\widetilde{W}-dv_R \right \rangle  ,
    \end{align*}
    and we use it to prove the following Lemma:
    \begin{lemma}\label{First Lemma in Fefferman-Stein estimate}
        For, $x\in B_r\left(0\right)$ and $0<\rho<R$ such that  $B_{\rho}\left(x\right)\subset \subset B_R\left(x\right)\subset \subset B_r\left(0\right)$ we have the following estimate:
        \begin{align}\label{First Lemma Inequality in FS Application}
        \left \lbrace
        \begin{aligned}
         \fint_{B_\rho\left(x\right)}& \left \lvert V_p\left(d\widetilde{W}+\widetilde{F}\right)-\left(V_p\left(d\widetilde{W}+\widetilde{F}\right)\right)_{\rho} \right \rvert^2\\
        &\quad\begin{aligned}[t]
            \leq 4&c_U \left(\frac{\rho}{R}\right)^{2\alpha_1}\fint_{B_R\left(x\right)} \left \lvert d\widetilde{W}\right \rvert^p
            +c\left(p,\epsilon\right) \left(\frac{R}{\rho}\right)^n\fint_{B_R} \left \lvert \widetilde{F}\right \rvert^{p}\\
            &+\left(\frac{R}{\rho}\right)^n\left[\frac{2^{p+1}}{c_0\left(p\right)} K_p\left \lVert G-I\right \rVert_{L^\infty\left(B_r\left(0\right)\right)}+\epsilon\right]\fint_{B_R\left(x\right)}\left \lvert d\widetilde{W}\right \rvert^p.
        \end{aligned}
        \end{aligned}
        \right.
    \end{align}
    \end{lemma}
    \begin{proof}
    Define
    \begin{align*}
         &I_1:=\int_{B_R} \left \langle \left \lvert d\widetilde{W}+\widetilde{F}\right  \rvert^{p-2} \left(d\widetilde{W}+\widetilde{F}\right)- \left \lvert dv_R\right \rvert^{p-2}dv_R, \left(d\widetilde{W}+\widetilde{F}\right)-dv_R\right \rangle,\\
         &I_2:=\int_{B_R} \left \langle \left \lvert d\widetilde{W}+\widetilde{F}\right  \rvert^{p-2} \left(d\widetilde{W}+\widetilde{F}\right)- \left \lvert dv_R\right \rvert^{p-2}dv_R, \widetilde{F}\right \rangle,\\
         &I_3
        :=\begin{aligned}[t]
        \int_{B_R} &\left \langle \left \langle G \left(d\widetilde{W}+\widetilde{F}\right), d\widetilde{W}+\widetilde{F}\right  \rangle^{\frac{p-2}{2}}\right.\\ \times &G \left(d\widetilde{W}+\widetilde{F}\right)\left.- \left \lvert d\widetilde{W}+\widetilde{F}\right  \rvert^{p-2}\left(d\widetilde{W}+\widetilde{F}\right)
        , d\widetilde{W}-dv_R\right \rangle.
        \end{aligned}
    \end{align*}
Notice, $I_1-I_2+I_3=0$ and  $I_1\leq \left \lvert I_2\right \rvert+ \left \lvert I_3\right \rvert$ and by \eqref{V upsim} we have
    \begin{align}\label{I_1}
        I_1 \geq c_0\left(p\right) \int_{B_R}\left \lvert V_p\left(d\widetilde{W}+\widetilde{F}\right)-V_p\left(dv_R\right) \right \rvert^2
    \end{align}
    On the other hand using Young's Inequality with exponents $\left(\frac{p}{p-1},p\right)$ and \eqref{in particular minimization} we estimate $I_2$:
    \begin{align*}
        I_2
        &\leq c\left(p\right)  \int_{B_R}\left[ \left \lvert d\widetilde{W}\right\rvert +\left \lvert \widetilde{F}\right\rvert+\left \lvert dv_R\right\rvert \right ]^{p-1} \left \lvert\widetilde{F}\right \rvert\\
        &\leq \begin{aligned}[t]\epsilon  &\int_{B_R}\left[ \left \lvert d\widetilde{W}\right\rvert+\left \lvert \widetilde{F}\right\rvert+\left \lvert dv_R\right\rvert \right ]^{p}+c_\epsilon \int_{B_R}\left \lvert\widetilde{F}\right \rvert^p
        \end{aligned}\sim \epsilon\int_{B_R}\left \lvert d\widetilde{W}\right\rvert ^{p}+c_\epsilon\int_{B_R} \left \lvert\widetilde{F}\right \rvert^p.
    \end{align*}
    Similarly with Young's inequality with exponents $\left(\frac{p}{p-1},p\right)$, using \eqref{Kp} we estimate $I_3$:
    \begin{align*}
        I_3
        &\begin{aligned}[t]\leq 
        \int_{B_R} &\left \lvert \left \langle G \left(d\widetilde{W}+\widetilde{F}\right), d\widetilde{W}+\widetilde{F}\right  \rangle^{\frac{p-2}{2}}G \left(d\widetilde{W}+\widetilde{F}\right)- \left \lvert d\widetilde{W}+\widetilde{F}\right  \rvert^{p-2}\left(d\widetilde{W}+\widetilde{F}\right)
        \right \rvert \left \lvert d\widetilde{W}-dv_R\right \rvert
        \end{aligned}\\
        &\leq K_p\left \lVert G-I\right \rVert_{L^\infty\left(B_r\left(0\right)\right)}\int_{B_R}\left[\left \lvert d\widetilde{W}+ \widetilde{F}\right \rvert\right]^{p-1}\left \lvert d\widetilde{W}-{dv_R} \right \rvert\\
        &\leq K_p\begin{aligned}[t]&\left \lVert G-I\right \rVert_{L^\infty\left(B_r\left(0\right)\right)}\left[\frac{p-1}{p}\int_{B_R}\left[\left \lvert d\widetilde{W}+ \widetilde{F}\right \rvert\right]^{p}+\frac{1}{p}\int_{B_R}\left \lvert d\widetilde{W}-{dv_R} \right \rvert^p\right].
        \end{aligned}
        \end{align*}
    Now, $\left \lvert d\widetilde{W}+\widetilde{F} \right \rvert^p\leq 2^{p-1}\left(\left \lvert d\widetilde{W}\right \rvert^p+\left \lvert \widetilde{F}\right \rvert^p\right)$ and $\left \lvert d\widetilde{W}- dv_R \right \rvert^p\leq 2^{p-1}\left(\left \lvert d\widetilde{W}\right \rvert^p+\left \lvert dv_R\right \rvert^p\right)$ along with \eqref{in particular minimization} gives 
    \begin{align}\label{I_3}
        I_3 \leq 2^pK_p\left \lVert G-I\right \rVert_{L^\infty\left(B_r\left(0\right)\right)}\int_{B_R}\left[\left \lvert d\widetilde{W}\right \rvert^p+ \left \lvert \widetilde{F}\right \rvert^{p}\right].
    \end{align}
    Combining \eqref{I_1},\eqref{I2 inequality} and \eqref{I_3} with the fact $I_1\leq \left \lvert I_2 \right \rvert+\left \lvert I_3\right \rvert$ we derive the following:
    \begin{align}\label{w+F-v to Kp inq}
        \int_{B_R}& \left \lvert V_p\left(d\widetilde{W}+\widetilde{F}\right)-V_p\left(dv_R\right)\right \rvert^2\notag\\
        &\leq \left[\frac{2^p}{c_0\left(p\right)} K_p\left \lVert G-I\right \rVert_{L^\infty\left(B_r\left(0\right)\right)}+\epsilon\right]\int_{B_R}\left \lvert d\widetilde{W}\right \rvert^p+ c\left(p,\epsilon\right) \int_{B_R} \left \lvert \widetilde{F}\right \rvert^{p}.
    \end{align}
    Fix $x\in B_r\left(0\right)$ and $0<\rho<R$ such that  $B_{\rho}\left(x\right)\subset \subset B_R\left(x\right)\subset \subset B_r\left(0\right)$. Then 
    \begin{align*}
        &\fint_{B_\rho\left(x\right)} \left \lvert V_p\left(d\widetilde{W}+\widetilde{F}\right)-\left(V_p\left(d\widetilde{W}+\widetilde{F}\right)\right)_{\rho} \right \rvert^2\\
        &\quad \begin{aligned}[t]\leq \fint_{B_\rho\left(x\right)} \left \lvert V_p\left(d\widetilde{W}+\widetilde{F}\right)-\left(V_p\left(dv_R\right)\right)_{\rho} \right \rvert^2\left[\text{By minimality of the mean}\right]
        \end{aligned}\\
         &\quad\begin{aligned}[t]
            \leq 2 \left[\fint_{B_\rho\left(x\right)} \left \lvert V_p\left(dv_R\right)-\left(V_p\left(dv_R\right)\right)_{\rho}\right \rvert^2+\fint_{B_\rho\left(x\right)} \left \lvert V_p\left(d\widetilde{W}+\widetilde{F}\right)-V_p\left(dv_R\right)\right \rvert^2\right].
        \end{aligned}
        \end{align*}
        By \eqref{Ulhenbeck Application} and minimality of the mean we have
        \begin{align*}
            &\fint_{B_\rho\left(x\right)} \left \lvert V_p\left(dv_R\right)-\left(V_p\left(dv_R\right)\right)_{\rho}\right \rvert^2\\
            &\quad\leq c_{U} \left(\frac{\rho}{R}\right)^{2\alpha_1}\fint_{B_R\left(x\right)} \left \lvert V_p\left(dv_R\right)-\left(V_p\left(dv_R\right)\right)_{R}\right \rvert^2\\
            &\quad\leq c_{U} \left(\frac{\rho}{R}\right)^{2\alpha_1}\fint_{B_R\left(x\right)} \left \lvert V_p\left(dv_R\right)\right \rvert^2=c_U \left(\frac{\rho}{R}\right)^{2\alpha_1}\fint_{B_R\left(x\right)} \left \lvert dv_R\right \rvert^p\leq c_U \left(\frac{\rho}{R}\right)^{2\alpha_1}\fint_{B_R\left(x\right)} \left \lvert d\widetilde{W}\right \rvert^p.
        \end{align*}
        In the last line we have used \eqref{in particular minimization}.
        On the other hand
        \begin{align*}
            \fint_{B_\rho\left(x\right)} \left \lvert V_p\left(d\widetilde{W}+\widetilde{F}\right)-V_p\left(dv_R\right)\right \rvert^2& \leq \left(\frac{R}{r}\right)^n\fint_{B_R\left(x\right)} \left \lvert V_p\left(d\widetilde{W}+\widetilde{F}\right)-V_p\left(dv_R\right)\right \rvert^2.
        \end{align*}
        Combining the previous estimate along with \eqref{w+F-v to Kp inq} we prove \eqref{First Lemma Inequality in FS Application}.
    \end{proof}
    Fix, $x\in B_\frac{r}{16}\left(0\right)$ and $0<\rho_0<\frac{r}{16}$ such that 
    \begin{align}\label{rho0}\delta_*^{-1}\rho_0<\frac{r}{2},
    \end{align}
    where $\delta_*$ is chosen in \eqref{delta}; thus $\rho_0$ depends on $n,k,p,\widetilde{q}, \alpha_1,c_U$ and $r$. On the other hand $c_U$ and $\alpha_1$ depends on $n,k,p$. Therefore, $\rho_0$ depends only on $n,k,p,\widetilde{q}$ and $r$. Choose, $\rho\in \left(0,\rho_0\right]$ and taking $R=\delta_*^{-1} \rho$ we arrive at the following:
    
    \begin{align}\label{First use of HL}
        \fint_{B_\rho\left(x\right)}&\left \lvert V_p\left(d\widetilde{W}+\widetilde{F}\right)-\left(V_p\left( d\widetilde{W}+\widetilde{F}\right)\right)_{\rho}\right \rvert^2
        \notag \\ 
         &\begin{aligned}[t]
            \quad\leq 4&c_U \delta_*^{2\alpha_1}\fint_{B_R\left(x\right)} \left \lvert d\widetilde{W}\right \rvert^p+c\left(p,\epsilon\right) \delta_*^{-n}\fint_{B_R\left(x\right)} \left \lvert \widetilde{F}\right \rvert^{p}\notag \\
            &+\delta_*^{-n}\left[\frac{2^{p+1}}{c_0\left(p\right)} K_p\left \lVert G-I\right \rVert_{L^\infty\left(B_r\left(0\right)\right)}+\epsilon\right]\fint_{B_R\left(x\right)}\left \lvert d\widetilde{W}\right \rvert^p
        \end{aligned}
        \notag \\
        &\begin{aligned}[t]\quad\leq 4c_U&\delta_*^{2\alpha_1}M\left(\left \lvert d\widetilde{W}\right \rvert^p\right) \left(x\right)+c\left(\epsilon,p\right)\delta_*^{-n}M\left(\left \lvert F\right\rvert^p\right)\left(x\right) \\&+\delta_*^{-n}\left[\frac{2^{p+1}}{c_0\left(p\right)} K_p\left \lVert G-I\right \rVert_{L^\infty\left(B_r\left(0\right)\right)}+\epsilon\right] M\left(\left \lvert d\widetilde{W}\right \rvert^p\right)\left(x\right).
        \end{aligned}
    \end{align}
    Let $Q_0$ be the cube inscribed in $B_{\rho_0}\left(0\right)$ with sides parallel to the co-ordinate axes, choose $x\in Q\left(x\right)\subseteq Q_0$. There exists $\rho \in \left(0,\rho_0\right]$ such that $Q\left(x\right)$ is inscribed inside $B_\rho\left(x\right)$ hence $\left \lvert B_\rho\left(x\right)\right \rvert=c_0^*\left(n\right)\left \lvert Q\left(x\right)\right \rvert$. Again using minimality of mean we have the following estimate:
    \begin{align*}
    \fint_{Q\left(x\right)}&\left \lvert V_p\left(d\widetilde{W}+\widetilde{F}\right)\left(y\right)-\left(V_p\left( d\widetilde{W}+\widetilde{F}\right)\right)_{Q\left(x\right)}\right \rvert^2dy\\
    &\leq \fint_{Q\left(x\right)}\left \lvert V_p\left(d\widetilde{W}+\widetilde{F}\right)\left(y\right)-\left(V_p\left( d\widetilde{W}+\widetilde{F}\right)\right)_{B_{\rho}\left(x\right)}\right \rvert^2dy\\
        &\leq c_0^*\left(n\right)\fint_{B_\rho\left(x\right)}\left \lvert V_p\left(d\widetilde{W}+\widetilde{F}\right)\left(y\right)-\left(V_p\left( d\widetilde{W}+\widetilde{F}\right)\right)_{B_\rho\left(x\right)}\right \rvert^2dy.
    \end{align*}
    Using Jensen inequality we get that
    \begin{align*}
    &\left[\fint_{Q\left(x\right)}\left \lvert V_p\left(d\widetilde{W}+\widetilde{F}\right)\left(y\right)-\left(V_p\left( d\widetilde{W}+\widetilde{F}\right)\right)_{Q\left(x\right)}\right \rvert dy\right]^2\\
    &\qquad\leq \fint_{Q\left(x\right)}\left \lvert V_p\left(d\widetilde{W}+\widetilde{F}\right)\left(y\right)-\left(V_p\left( d\widetilde{W}+\widetilde{F}\right)\right)_{Q\left(x\right)}\right \rvert^2dy\\
        &\qquad\leq c_0^*\left(n\right)\fint_{B_\rho\left(x\right)}\left \lvert V_p\left(d\widetilde{W}+\widetilde{F}\right)\left(y\right)-\left(V_p\left( d\widetilde{W}+\widetilde{F}\right)\right)_{B_\rho\left(x\right)}\right \rvert^2dy\\
        &\qquad\qquad\qquad\qquad\qquad\qquad\text{ [$c^*_0\left(n\right)$ was  defined in \eqref{c*0}]}.
    \end{align*}
    Taking the supremum over all cubes $Q\left(x\right)(\subseteq Q_0)$ centered at $x$ and using \eqref{First use of HL} we get, 
    \begin{align*}
        &\left[ V_p\left( d\widetilde{W}+\widetilde{F}\right) ^{\#}\left(x\right)\right]^2\\ 
         &\quad\leq \begin{aligned}[t] 4c^*_0&\left(n\right)c_U\delta_*^{2\alpha_1}M\left(\left \lvert d\widetilde{W}\right \rvert^p\right) \left(x\right)+c_0^*\left(n\right)c\left(\epsilon,p\right)\delta_*^{-n}M\left(\left \lvert F\right\rvert^p\right)\left(x\right)\\&+c_0^*\left(n\right)\delta_*^{-n}\left[\frac{2^{p+1}}{c_0\left(p\right)} K_p\left \lVert G-I\right \rVert_{L^\infty\left(B_r\left(0\right)\right)}+\epsilon\right] M\left(\left \lvert d\widetilde{W}\right \rvert^p\right)\left(x\right).
        \end{aligned}
    \end{align*}
    Integrating in $Q_0$ gives,
    \begin{align*}
        \int_{Q_0}&\left[ V_p\left( d\widetilde{W}+\widetilde{F}\right)^{\#}\right]^{\frac{2\tilde{q}}{p}}\\ 
         &\begin{aligned}[t]\leq &\left[8c^*_0\left(n\right)c_U\right]^{\frac{\tilde{q}}{p}}\delta_*^{\frac{2\alpha_1\tilde{q}}{p}}\int_{Q_0}M\left(\left \lvert d\widetilde{W}\right \rvert^p\right) ^{\frac{\tilde{q}}{p}}+\left[2c_0^*\left(n\right)c\left(\epsilon,p\right)\right]^\frac{\tilde{q}}{p}\delta_*^{-\frac{n\tilde{q}}{p}}\int_{Q_0}M\left(\left \lvert F\right\rvert^p\right)^{\frac{\tilde{q}}{p}}\\&+\left[2c_0^*\left(n\right)\right]^{\frac{\tilde{q}}{p}}\delta_*^{-\frac{n\tilde{q}}{p}}\left[\frac{2^{p+1}}{c_0\left(p\right)} K_p\left \lVert G-I\right \rVert_{L^\infty\left(B_r\left(0\right)\right)}+\epsilon\right]^\frac{\tilde{q}}{p} \int_{Q_0}M\left(\left \lvert d\widetilde{W}\right \rvert^p\right)^{\frac{\tilde{q}}{p}}.
         \end{aligned}\\
    \end{align*}
 Now both $\left \lvert d\widetilde{W}\right \rvert^p, \left \lvert \widetilde{F}\right \rvert^p \in L^{\frac{\widetilde{q}}{p}}\left(B_{\frac{r}{4}}\left(0\right)\right)$. Therefore using \eqref{Mf <f} with $\gamma=\frac{\tilde{q}}{p}>1$ we get the following:
    \begin{align*}
    \int_{Q_0}&\left[ V_p\left( d\widetilde{W}+\widetilde{F}\right)^{\#}\right]^{\frac{2\tilde{q}}{p}}\\
        &\begin{aligned}[t]\leq &\left[8c^*_0\left(n\right)c_U\right]^{\frac{\tilde{q}}{p}}\delta_*^{\frac{2\alpha_1\tilde{q}}{p}}A_0\left(n,\frac{\widetilde{q}}{p}\right)\int_{Q_0}\left \lvert d\widetilde{W}\right \rvert^{{\tilde{q}}}\\
        &+\left[2c_0^*\left(n\right)\delta_*^{-n}c\left(\epsilon,p\right)\right]^\frac{\tilde{q}}{p}A_0\left(n,\frac{\widetilde{q}}{p}\right)\int_{Q_0}\left \lvert F\right\rvert^{\tilde{q}}\\
        &+\left[2c_0^*\left(n\right) \delta_*^{-n}\right]^{\frac{\tilde{q}}{p}}\begin{aligned}[t]&\left[\frac{2^{p+1}}{c_0\left(p\right)} K_p\left \lVert G-I\right \rVert_{L^\infty\left(B_r\left(0\right)\right)}+\epsilon\right]^\frac{\tilde{q}}{p}A_0\left(n,\frac{\widetilde{q}}{p}\right) \int_{Q_0}\left \lvert d\widetilde{W}\right \rvert^{\widetilde{q}}.
        \end{aligned}
        \end{aligned}
    \end{align*}
    We use the Fefferman-Stein Inequality (see \eqref{Fefferman Stein Cube} with $\gamma=\frac{2\tilde{q}}{p}$) and the following:
    \begin{align}
        \int_{Q_0}&\left \lvert V_p\left(d\widetilde{W}+\widetilde{F}\right) \right \rvert^{\frac{2\widetilde{q}}{p}}\notag\\
        &\leq A\left(n,\frac{2\tilde{q}}{p}\right)\left[\left \lvert Q_0\right \rvert\left[\fint_{Q_0}\left \lvert V_p\left(d\widetilde{W}+\widetilde{F}\right)\right \rvert\right]^{\frac{2\widetilde{q}}{p}}+\int_{Q_0}\left \lvert V_p\left(d\widetilde{W}+\widetilde{F}\right)^{\#} \right \rvert^{\frac{2\widetilde{q}}{p}} \right]
        \end{align}\\
        Using Jensen's inequality and the fact that $\left \lvert V_p\left(z\right)\right \rvert^2=\left \lvert z\right \rvert^p$ implies, 
        \begin{align*}
            \left[\fint_{Q_0}\left \lvert V_p\left(d\widetilde{W}+\widetilde{F}\right)\right \rvert\right]^{\frac{2\widetilde{q}}{p}}& \leq \left[\fint_{Q_0}\left \lvert d\widetilde{W}+\widetilde{F} \right \rvert^p\right]^{\frac{\tilde{q}}{p}}.
        \end{align*}
        Combining we get that
    \begin{align}\label{Fefferman-Stein Application First}
       \left \lbrace
       \begin{aligned}
        \int_{Q_0}&\left \lvert V_p\left(d\widetilde{W}+\widetilde{F}\right) \right \rvert^{\frac{2\widetilde{q}}{p}}\\
        &\leq \begin{aligned}[t]
        &c\left(n,\widetilde{q},{p},\rho_0\right) \left[\int_{Q_0}\left \lvert d\widetilde{W}+\widetilde{F} \right \rvert^p\right]^{\frac{\tilde{q}}{p}}+c\left(n,\tilde{q},p,\delta_*,\alpha_1,\epsilon\right)\int_{Q_0}\left \lvert F\right\rvert^{\tilde{q}} \\&
        \quad +\left[8c^*_0\left(n\right)c_U\right]^{\frac{\tilde{q}}{p}}\delta_*^{\frac{2\alpha_1\tilde{q}}{p}}A_0\left(n,\frac{\widetilde{q}}{p}\right)A\left(n,\frac{2\widetilde{q}}{p}\right)\int_{Q_0}\left \lvert d\widetilde{W}\right \rvert^{{\tilde{q}}}\\
        &\qquad+\left[4c_0^*\left(n\right)\right]^{\frac{\tilde{q}}{p}}\begin{aligned}[t]\delta_*^{-\frac{n\tilde{q}}{p}}&\left[\frac{2^{p+1}}{c_0\left(p\right)} K_p\left \lVert G-I\right \rVert_{L^\infty\left(B_r\left(0\right)\right)}+\epsilon\right]^\frac{\tilde{q}}{p} \\
        &\times A_0\left(n,\frac{\widetilde{q}}{p}\right)A\left(n,\frac{2\widetilde{q}}{p}\right) \int_{Q_0}\left \lvert d\widetilde{W}\right \rvert^{\tilde{q}}.
        \end{aligned}
        \end{aligned}
        \end{aligned}
        \right.
    \end{align} 
    By our choice of $\delta_*$ in \eqref{delta} (which depends only on $n,p,\tilde{q}$ and $k$) we get the following
    \begin{align}\label{justi delta}
        \left[8c^*_0\left(n\right)c_U\right]^{\frac{\tilde{q}}{p}}\delta_*^{\frac{2\alpha_1\tilde{q}}{p}}A_0\left(n,\frac{\widetilde{q}}{p}\right)A\left(n,\frac{2\widetilde{q}}{p}\right)\leq \frac{1}{1000. 2^{\tilde{q}}}.
    \end{align}
    Depending on this $\delta_*$ fix $\epsilon>0$ such that 
    \begin{align}\label{epsilon}
        \left[8c_0^*\left(n\right)\delta_*^{-n}\epsilon\right]^\frac{\tilde{q}}{p}A_0\left(n,\frac{\widetilde{q}}{p}\right)A\left(n,\frac{2\widetilde{q}}{p}\right)\leq \frac{1}{1000.2^{\tilde{q}}}.
    \end{align}
     Inequality \eqref{r} of $*$-condition and \eqref{epsilon} give the following estimate:
    \begin{align*}
        \left[4c_0^*\left(n\right)\right]^{\frac{\tilde{q}}{p}}\delta_*^{-\frac{n\tilde{q}}{p}}\left[\frac{2^{p+1}}{c_0\left(p\right)} K_p\left \lVert G-I\right \rVert_{L^\infty\left(B_r\left(0\right)\right)}+\epsilon\right]^\frac{\tilde{q}}{p}  A_0\left(n,\frac{\widetilde{q}}{p}\right)A&\left(n,\frac{2\widetilde{q}}{p}\right)\leq \frac{1}{500 2^{\widetilde{q}}}.
    \end{align*}
    Applying that in \eqref{Fefferman-Stein Application First} we get
    \begin{align*}
        \int_{Q_0}\left \lvert d\widetilde{W}+\widetilde{F} \right \rvert^{\tilde{q}}
        &= \int_{Q_0}\left \lvert V_p\left(d\widetilde{W}+\widetilde{F}\right) \right \rvert^{\frac{2\widetilde{q}}{p}}\\
        &\leq \begin{aligned}[t]
        &c\left(n,\widetilde{q},{p},\rho_0\right) \left[\int_{Q_0}\left \lvert d\widetilde{W}+\widetilde{F} \right \rvert^p\right]^{\frac{\tilde{q}}{p}}
        +\frac{1}{1000. 2^{\tilde{q}}}\int_{Q_0}\left \lvert d\widetilde{W}\right \rvert^{{\tilde{q}}}\\
        &\qquad \qquad+c\left(n,\tilde{q},p,k\right)\int_{Q_0}\left \lvert F\right\rvert^{\tilde{q}}+ \frac{1}{500.2^{\tilde{q}}}\int_{Q_0}\left \lvert d\widetilde{W}\right \rvert^{\tilde{q}}
        \end{aligned}\\
        &\leq \begin{aligned}[t]
        &c\left(n,\widetilde{q},{p},k,\rho_0\right) \left[\int_{Q_0}\left \lvert d\widetilde{W}+\widetilde{F} \right \rvert^p\right]^{\frac{\tilde{q}}{p}}
        +\frac{1}{10}\int_{Q_0}\left \lvert d\widetilde{W}+\widetilde{F}\right \rvert^{{\tilde{q}}}\\
        &\qquad\qquad+c\left(n,\tilde{q},p,k,\rho_0\right)\int_{Q_0}\left \lvert \widetilde{F}\right\rvert^{\tilde{q}}.\\
        \end{aligned}
    \end{align*}
Since, $Q_0$ is a (maximal) cube inscribed in $B_{\rho_0}\left(0\right)$;  $\left \lvert Q_0\right \rvert$ is a function on $\rho_0$; which we used in the last estimate. On on the other hand $\rho_0$ depends on $\delta_*$ and $r$ hence on $n,k,p,\tilde{q}$ and $r$ only. Hence $\left \lvert Q_0\right \rvert$  depends only on $n,p,\widetilde{q},k$ and $r$. Finally after absorption of $\frac{1}{10}\int_{Q_0}\left \lvert d\widetilde{W}+\widetilde{F} \right \rvert^{\widetilde{q}}$ on the right hand side we get the following 
\begin{align*}
    \int_{Q_0}&\left \lvert d\widetilde{W}+\widetilde{F} \right \rvert^{\widetilde{q}}  \leq \begin{aligned}[t]
    &c\left(n,\tilde{q},p,k,r\right) \left[\left[\int_{Q_0}\left \lvert d\widetilde{W}+\widetilde{F} \right \rvert^p\right]^{\frac{\tilde{q}}{p}}
    +\int_{Q_0}\left \lvert F\right\rvert^{\tilde{q}}\right].
\end{aligned}
\end{align*}
\end{proof}
  We use Proposition \ref{Proposition with F in W2,p} and Proposition \ref{Application of HL FS KU Proposition} and finally finish the proof of Theorem \ref{G Theorem} below:
     \begin{proof}
     Choose, a sequence $\left\{\widetilde{F}_j\right\}_{j\in \mathbb{N}}$ in $ C^\infty_c\left(B_r\left(0\right);\Lambda^{k+1}\right)$ converging to $\widetilde{F}$ in  $L^{\widetilde{q}}\left(B_r\left(0\right);\Lambda^{k+1}\right)$ and for each $j \in  \mathbb{N}$ 
     let us consider the minimization problem:
     \begin{align}\label{minimization prblm Ij}
         \zeta\mapsto I_j\left[\zeta\right]:=\int_{B_r\left(0\right)}\left \langle G\left(d\zeta+\widetilde{F}_j\right),d\zeta+\widetilde{F}_j\right \rangle^{\frac{p}{2}}; \text{ where $\zeta \in\widetilde{W}+ W^{1,p}_{d^*T}\left(B_r\left(0\right);\Lambda^k\right)$.}
     \end{align}
     \begin{lemma}\label{existenceof the minimizer I_j Lemma}
         For each $j\in \mathbb{N}$ the minimizer for the functional defined in \eqref{minimization prblm Ij} exists uniquely.
     \end{lemma}
     \begin{proof}
         Fix, $j$ in $\mathbb{N}$ and by \eqref{ellipticity} for each $\zeta \in\widetilde{W}+ W^{1,p}_{d^*T}\left(B_r\left(0\right);\Lambda^k\right)$, $I_j\left[\zeta\right] \in \left[0,\infty \right)$.
         Define, $M_0=\inf_{\zeta \in \mathcal{A}}I_j\left[\zeta\right]$ where $\mathcal{A}=\widetilde{W}+ W^{1,p}_{d^*T}\left(B_r\left(0\right);\Lambda^k\right)$. So there exists a sequence, $\left\{\zeta_m\right\}_m$ in $\mathcal{A}$ such that, $I_j\left[\zeta_m\right]\rightarrow M_0$ as $m\rightarrow \infty$. Thus, $I_j\left[\zeta_m\right] \leq M_0+1$ eventually. Hence by \eqref{ellipticity} we have,
         \begin{align*}
             \int_{B_r\left(0\right)}&\left \lvert d\zeta_m-d\widetilde{W}\right \rvert^p \\
             &\lesssim \int_{B_r\left(0\right)}\left \lvert d\zeta_m+\widetilde{F}\right \rvert^p
             +\int_{B_r\left(0\right)}\left \lvert d\widetilde{W}-\widetilde{F}\right \rvert^p\\
             &\lesssim \begin{aligned}[t]I_j\left[\zeta_m\right]+\int_{B_r\left(0\right)}\left \lvert d\widetilde{W}-\widetilde{F}\right \rvert^p=\int_{B_r\left(0\right)}&\left \langle G\left(d\zeta_m+\widetilde{F}_j\right),d\zeta_m+\widetilde{F}_j\right \rangle^{\frac{p}{2}}+\int_{B_r\left(0\right)}\left \lvert d\widetilde{W}-\widetilde{F}\right \rvert^p\\
             &\leq  M_0+1+\int_{B_r\left(0\right)}\left \lvert d\widetilde{W}-\widetilde{F}\right \rvert^p \text{ eventually}.
             \end{aligned}
         \end{align*}
         But, $\zeta_m-\widetilde{W}$ is in $W^{1,p}_{d^*,T}\left(B_r\left(0\right);\Lambda^k\right)$ and since $\zeta \rightarrow \left \lVert d\zeta \right \rVert_{L^p\left(B_r\left(0\right)\right)}$ is (an equivalent) norm in $W^{1,p}_{d^*,T}\left(B_r\left(0\right);\Lambda^k\right)$ (see Corollary \ref{equivnorm}). Hence using Banach-Alaoglu and separation form of Hahn-Banach theorem we have $\xi_j \in W^{1,p}_{d^*,T}\left(B_r\left(0\right);\Lambda^k\right)$ such that
         \begin{align*}
             \zeta_m-\widetilde{W} \rightharpoonup \xi_j \text{ in } W^{1,p}\left(B_r\left(0\right);\Lambda^k\right); \text{ hence in particular } \nabla \left( \zeta_m-\widetilde{W}\right) \rightharpoonup \xi_j \text{ in } L^p\left(B_r\left(0\right);\Lambda^k\right).
         \end{align*}
         Define, $\widetilde{W}_j=\widetilde{W}+\xi_j\in \mathcal{A}$ and we get that
         \begin{align*}
            d\zeta_m+\widetilde{F}_j \rightharpoonup d\widetilde{W}_j+F_j \text{ in } L^p\left(B_r\left(0\right);\Lambda^k\right).
         \end{align*}
         Now we use $\mathcal{H}\left(x,z\right)=\left \langle G\left(x\right)z,z\right \rangle^{\frac{p}{2}}$ in Lemma \ref{Convex Functional Lemma from Daco's Book} and get the following:
         \begin{align*}
             M_0\leq I_j\left[\widetilde{W}_j\right]\leq \liminf_{m\rightarrow \infty}I_j\left[\zeta_m\right]=\lim_{m\rightarrow \infty}I_j\left[\zeta_m\right]=M_0.
         \end{align*}
         Therefore, $\widetilde{W}_j$ is the minimizer of the functional mentioned in \eqref{minimization prblm Ij} and the uniqueness of $\widetilde{W}_j$ follows from the strict-convexity of $z\mapsto \left \langle G\left(x\right)z,z\right \rangle^{\frac{p}{2}}$ for $x\in B_r\left(0\right)$.
     \end{proof}
    Since, $\widetilde{W}_j$ is the unique minimizer for the functional defined in \eqref{minimization prblm Ij} we have the following:
        \begin{align*}
            \int_{B_r\left(0\right)}\left \langle G\left(d\widetilde{W}_j+\widetilde{F}_j\right),d\widetilde{W}_j+\widetilde{F}_j\right \rangle^{\frac{p}{2}} \leq \int_{B_r\left(0\right)}\left \langle G\left(d\widetilde{W}+\widetilde{F}_j\right),d\widetilde{W}+\widetilde{F}_j\right \rangle^{\frac{p}{2}} 
        \end{align*}
        and using \eqref{ellipticity} we get
        \begin{align}\label{Ij min ineq}
            \int_{B_r\left(0\right)}\left \lvert d\widetilde{W}_j+\widetilde{F}_j \right \rvert^{p} \lesssim \int_{B_r\left(0\right)}\left \lvert d\widetilde{W}+\widetilde{F}_j \right \rvert^{p}.
        \end{align}
       Thus for every $\eta \in W^{1,p}_{d^*,T}\left(B_r\left(0\right);\Lambda^k\right)$,
    \begin{align}\label{d*T min problem for Wj}
        \int_{B_r\left(0\right)}\left \langle \left \langle G \left(d\widetilde{W}_j+\widetilde{F}_j\right), d\widetilde{W}_j+\widetilde{F}_j\right  \rangle^{\frac{p-2}{2}}G \left(d \widetilde{W}_j+\widetilde{F}_j\right),d\eta \right \rangle =0.
    \end{align}
    Moreover by Proposition \ref{Application of HL FS KU Proposition} for each $j=1,2,..,$ we have cube (with sides parallel to the co-ordinate axes) such that  
    \begin{align}\label{smooth estimate j}
        \int_{Q_0}&\left \lvert d\widetilde{W}_j+\widetilde{F}_j \right \rvert^{\widetilde{q}}  \leq \begin{aligned}[t]
        &c\left(n,\tilde{q},p,k,r\right) \left[\left[\int_{Q_0}\left \lvert d\widetilde{W}_j+\widetilde{F}_j \right \rvert^p\right]^{\frac{\tilde{q}}{p}}
        +\int_{Q_0}\left \lvert \widetilde{F}_j\right\rvert^{\tilde{q}}\right].
    \end{aligned}
    \end{align}
    Recall the weak formulation \eqref{d* T weal formulation}:
    \begin{align}\label{d*T  formulation in  step for 3 W}
        \int_{B_r\left(0\right)}\left \langle \left \langle G \left(d\widetilde{W}+\widetilde{F}\right), d\widetilde{W}+\widetilde{F}\right  \rangle^{\frac{p-2}{2}}G \left(d \widetilde{W}+\widetilde{F}\right),d\eta \right \rangle =0,
    \end{align}
    for  $\eta \in W^{1,p}_{d^*,T}\left(B_r\left(0\right);\Lambda^k\right)$.
    By \eqref{d*T min problem for Wj} and \eqref{d*T  formulation in  step for 3 W} we prove the following lemma:
    \begin{lemma}
        $d\widetilde{W}_j+\widetilde{F_j} \rightarrow d\widetilde{W}+\widetilde{F}$ point-wise almost everywhere in $B_r\left(0\right)$.  
    \end{lemma}
    \begin{proof}
    Take $\eta =\widetilde{W}_j-\widetilde{W}\in W^{1,p}_{d^*,T}\left(B_r\left(0\right)\right)$ in \eqref{d*T min problem for Wj}, \eqref{d*T  formulation in  step for 3 W} and subtracting one equation from the other we get:
    \begin{multline*}
    \int_{B_r\left(0\right)}\left \langle \left \langle G \left(d\widetilde{W}_j+\widetilde{F}_j\right), d\widetilde{W}_j+\widetilde{F}_j\right  \rangle^{\frac{p-2}{2}}G \left(d \widetilde{W}_j+\widetilde{F}_j\right)\right.\\
    \quad\left.-\left \langle G \left(d\widetilde{W}+\widetilde{F}\right), d\widetilde{W}+\widetilde{F}\right  \rangle^{\frac{p-2}{2}}G \left(d \widetilde{W}+\widetilde{F}\right), d\widetilde{W}_j+\widetilde{F}_j-\left(d\widetilde{W}+\widetilde{F}\right) \right \rangle\\
        =\int_{B_r\left(0\right)}\left \langle \left \langle G \left(d\widetilde{W}_j+\widetilde{F}_j\right), d\widetilde{W}_j+\widetilde{F}_j\right  \rangle^{\frac{p-2}{2}}G \left(d \widetilde{W}_j+\widetilde{F}_j\right)\right.\\
    \left.-\left \langle G \left(d\widetilde{W}+\widetilde{F}\right), d\widetilde{W}+\widetilde{F}\right  \rangle^{\frac{p-2}{2}}G \left(d \widetilde{W}+\widetilde{F}\right), \widetilde{F}_j-\widetilde{F}\right \rangle.
   \end{multline*}
   Define,
   \begin{align*}
       &C_1:=\begin{aligned}[t]\int_{B_r\left(0\right)}&\left \langle \left \langle G \left(d\widetilde{W}_j+\widetilde{F}_j\right), d\widetilde{W}_j+\widetilde{F}_j\right  \rangle^{\frac{p-2}{2}}G \left(d \widetilde{W}_j+\widetilde{F}_j\right)\right.\\
    &\quad\left.-\left \langle G \left(d\widetilde{W}+\widetilde{F}\right), d\widetilde{W}+\widetilde{F}\right  \rangle^{\frac{p-2}{2}}G \left(d \widetilde{W}+\widetilde{F}\right), d\widetilde{W}_j+\widetilde{F}_j-\left(d\widetilde{W}+\widetilde{F}\right) \right \rangle,
    \end{aligned}\\
    &C_2:=\begin{aligned}[t]
    \int_{B_r\left(0\right)}&\left \langle \left \langle G \left(d\widetilde{W}_j+\widetilde{F}_j\right), d\widetilde{W}_j+\widetilde{F}_j\right  \rangle^{\frac{p-2}{2}}G \left(d \widetilde{W}_j+\widetilde{F}_j\right)\right.\\&
    \left.-\left \langle G \left(d\widetilde{W}+\widetilde{F}\right), d\widetilde{W}+\widetilde{F}\right  \rangle^{\frac{p-2}{2}}G \left(d \widetilde{W}+\widetilde{F}\right), \widetilde{F}_j-\widetilde{F}\right \rangle.
    \end{aligned}
   \end{align*}
   Thus, $C_1=C_2\leq \left \lvert C_2\right \rvert$.
    Using \eqref{c1p} we get
    \begin{align}\label{C_1}
        \int_{B_r\left(0\right)} \left \lvert V_p\left(d\widetilde{W}_j+\widetilde{F}_j\right)-V_p\left(d\widetilde{W}+\widetilde{F}\right) \right \rvert^2  \lesssim C_1.
    \end{align}
    On the other \eqref{ellipticity} gives the following: 
    \begin{align}\label{C_2}
        \left \lvert C_2\right \rvert \lesssim \int_{B_r\left(0\right)} \left[\left \lvert \left(d\widetilde{W}_j+\widetilde{F}_j\right)\right \rvert+\left \lvert \left(d\widetilde{W}+\widetilde{F}\right)\right \rvert\right]^{p-1}\left \lvert \widetilde{F}-\widetilde{F}_j\right \rvert.
    \end{align}
    Combining \eqref{C_1}, \eqref{C_2} with the fact $C_1\leq \left \lvert C_2\right \rvert$ and by
    H\"{o}lder's Inequality with exponents $\left(\frac{p}{p-1},p\right)$  we arrive at 
    \begin{align}\label{Vfjf}
        &\int_{B_r\left(0\right)} \left \lvert V_p\left(d\widetilde{W}_j+\widetilde{F}_j\right)-V_p\left(d\widetilde{W}+\widetilde{F}\right) \right \rvert^2\notag \\
        &\quad\lesssim \int_{B_r\left(0\right)} \left[\left \lvert \left(d\widetilde{W}_j+\widetilde{F}_j\right)\right \rvert+\left \lvert \left(d\widetilde{W}+\widetilde{F}\right)\right \rvert\right]^{p-1}\left \lvert \widetilde{F}-\widetilde{F}_j\right \rvert\notag \\
        &\quad\leq \left(\int_{B_r\left(0\right)} \left[\left \lvert \left(d\widetilde{W}_j+\widetilde{F}_j\right)\right \rvert+\left \lvert \left(d\widetilde{W}+\widetilde{F}\right)\right \rvert\right]^{p}\right)^{\frac{p-1}{p}}\left(\int_{B_r\left(0\right)}\left \lvert \widetilde{F}-\widetilde{F}_j\right \rvert^p\right)^{\frac{1}{p}} \notag \\
        &\quad\lesssim\left(\int_{B_r\left(0\right)} \left[\left \lvert \left(d\widetilde{W}+\widetilde{F}_j\right)\right \rvert+\left \lvert \left(d\widetilde{W}+\widetilde{F}\right)\right \rvert\right]^{p}\right)^{\frac{p-1}{p}}\left(\int_{B_r\left(0\right)}\left \lvert \widetilde{F}-\widetilde{F}_j\right \rvert^p\right)^{\frac{1}{p}}\text{ [by \eqref{Ij min ineq}]}. 
    \end{align}
    Let $p-2\geq 0$ then
    \begin{align*}
        &\int_{B_r\left(0\right)} \left \lvert \left(d\widetilde{W}_j-d\widetilde{W}\right)+ \left(\widetilde{F}_j-\widetilde{F}\right)\right \rvert^p\\
        &\quad\leq \int_{B_r\left(0\right)} \left[ \left\lvert d\widetilde{W}_j+\widetilde{F}_j\right\rvert+ \left \lvert d\widetilde{W}+\widetilde{F}\right\rvert\right]^{p-2}\left \lvert \left(d\widetilde{W}_j-d\widetilde{W}\right)+ \left(\widetilde{F}_j-\widetilde{F}\right)\right \rvert^2 \\
        &\quad\lesssim \int_{B} \left \lvert V_p\left(d\widetilde{W}_j+\widetilde{F}_j\right)-V_p\left(d\widetilde{W}+\widetilde{F}\right) \right \rvert^2\\
        &\quad\lesssim \left(\int_{B_r\left(0\right)}\left[\left \lvert d\widetilde{W}\right \rvert^p+\left \lvert \widetilde{F}_j\right \rvert^p+ \left \lvert \widetilde{F}\right \rvert^p\right]\right )^{\frac{p-1}{p}}\left(\int_{B_r\left(0\right)}\left \lvert \widetilde{F}_j-\widetilde{F}\right \rvert^p\right)^{\frac{1}{p}}.
    \end{align*}
    Since, $\widetilde{F}_j \rightarrow \widetilde{F}$ in $L^{\widetilde{q}}\left(B_r\left(0\right)\right)$, \eqref{Vfjf}, implies  
        $d\widetilde{W}_j+\widetilde{F}_j \rightarrow d\widetilde{W}+\widetilde{F} \text{ in } L^p\left(B_r\left(0\right)\right)$.\\
    For $1<p<2$ we define,
    \begin{align*}
        &\mathcal{A}:=\left\{ x\in B_r\left(0\right):\left \lvert \left(d\widetilde{W}_j+\widetilde{F}_j\right)- \left(d\widetilde{W}
        +\widetilde{F}\right)\right \rvert> \left \lvert d\widetilde{W}+\widetilde{F}\right \rvert\right\},\\
        &\mathcal{B}:=\left\{x\in B_r\left(0\right): \left \lvert \left(d\widetilde{W}_j+\widetilde{F}_j\right)- \left(d\widetilde{W}
        +\widetilde{F}\right)\right \rvert\leq  \left \lvert d\widetilde{W}+\widetilde{F}\right \rvert\right\}.
    \end{align*}
    Suppose, $x\in \mathcal{A}$, then using triangle inequality we have
    \begin{align*}
         &\left \lvert d\widetilde{W}_j+\widetilde{F}_j\right\rvert+  \left\lvert d\widetilde{W}
        +\widetilde{F}\right \rvert
        \leq  3\left \lvert \left(d\widetilde{W}_j+\widetilde{F}_j\right)- \left(d\widetilde{W}
        +\widetilde{F}\right)\right \rvert,\\
        &\text{ i.e. } 3^{p-2}\left(\left \lvert \left(d\widetilde{W}_j+\widetilde{F}_j\right)- \left(d\widetilde{W}
        +\widetilde{F}\right)\right \rvert\right)^{p-2} \leq  \left(\left \lvert d\widetilde{W}_j+\widetilde{F}_j\right\rvert+  \left\lvert d\widetilde{W}
        +\widetilde{F}\right \rvert\right)^{p-2};
    \end{align*}
    which implies 
    \begin{align*}
        \int_{\mathcal{A}}&\left \lvert \left(d\widetilde{W}_j+\widetilde{F}_j\right)- \left(d\widetilde{W}
        +\widetilde{F}\right)\right \rvert^p\\
        &=\int_{\mathcal{A}}\left \lvert \left(d\widetilde{W}_j+\widetilde{F}_j\right)- \left(d\widetilde{W}
        +\widetilde{F}\right)\right \rvert^{p-2}\left \lvert \left(d\widetilde{W}_j+\widetilde{F}_j\right)- \left(d\widetilde{W}
        +\widetilde{F}\right)\right \rvert\\
        &\lesssim\int_{\mathcal{A}}\left ( \left\lvert d\widetilde{W}_j+\widetilde{F}_j\right\rvert + \left\lvert d\widetilde{W}
        +\widetilde{F}\right\rvert\right)^{p-2}\left \lvert \left(d\widetilde{W}_j+\widetilde{F}_j\right)- \left(d\widetilde{W}
        +\widetilde{F}\right)\right \rvert\\
        &\lesssim \int_{B_r\left(0\right)}\left \lvert V_p \left(d\widetilde{W}_j+\widetilde{F}_j\right)-V_p \left(d\widetilde{W}+\widetilde{F}\right)\right \rvert^2\\
        &\lesssim \left(\int_{B_r\left(0\right)}\left[\left \lvert d\widetilde{W}\right \rvert^p+\left \lvert \widetilde{F}_j\right \rvert^p+ \left \lvert \widetilde{F}\right \rvert^p\right]\right )^{\frac{p-1}{p}}\left(\int_{B_r\left(0\right)}\left \lvert \widetilde{F}_j-\widetilde{F}\right \rvert^p\right)^{\frac{1}{p}}.
    \end{align*}
    On the other hand if $x\in \mathcal{B}$, then
   $\left \lvert d\widetilde{W}_j+\widetilde{F}_j\right\rvert+  \left\lvert d\widetilde{W}
        +\widetilde{F}\right \rvert\leq  3\left \lvert d\widetilde{W}
        +\widetilde{F}\right \rvert$ and hence
    by Young's Inequality with exponent $\left(2,2\right)$ the following estimate follows:
    \begin{align*}
       \int_{\mathcal{B}}&\left \lvert \left(d\widetilde{W}_j+\widetilde{F}_j\right)- \left(d\widetilde{W}
        +\widetilde{F}\right)\right \rvert^p\\
        &=\int_{\mathcal{B}}\left \lvert \left(d\widetilde{W}_j+\widetilde{F}_j\right)- \left(d\widetilde{W}
        +\widetilde{F}\right)\right \rvert^{{p-1}}\left \lvert \left(d\widetilde{W}_j+\widetilde{F}_j\right)- \left(d\widetilde{W}
        +\widetilde{F}\right)\right \rvert\\
        &\leq \int_{\mathcal{B}}\left (\left\lvert d\widetilde{W}_j+\widetilde{F}_j\right\rvert + \left\lvert d\widetilde{W}
        +\widetilde{F}\right\rvert\right)^{{p-1}}\left \lvert \left(d\widetilde{W}_j+\widetilde{F}_j\right)- \left(d\widetilde{W}
        +\widetilde{F}\right)\right \rvert\\
        &\lesssim \int_{\mathcal{B}}\left\lvert d\widetilde{W}
        +\widetilde{F}\right\rvert ^{{p-1}}\left \lvert \left(d\widetilde{W}_j+\widetilde{F}_j\right)- \left(d\widetilde{W}
        +\widetilde{F}\right)\right \rvert\\
        &= \int_{\mathcal{B}} \left(\epsilon^{-\frac{1}{2}}\left \lvert d\widetilde{W}+\widetilde{F}\right \rvert^{\frac{p-2}{2}}\left \lvert  \left(d\widetilde{W}_j+\widetilde{F}_j\right)-\left(d\widetilde{W}+\widetilde{F}\right)\right \rvert\right)\left(\epsilon^{\frac{1}{2}} \left \lvert d\widetilde{W}+\widetilde{F}\right \rvert^{\frac{p}{2}}\right)\\
        &\leq \begin{aligned}[t]\frac{c}{\epsilon}&\int_{B_r\left(0\right)}\left \lvert V_p\left(d\widetilde{W}+\widetilde{F}\right)-V_p\left(d\widetilde{W}_j+\widetilde{F}_j\right)\right \rvert^2+\epsilon\int_{\mathcal{B}}\left\lvert d\widetilde{W}+\widetilde{F}\right\rvert^{p}.
        \end{aligned}
    \end{align*}
    Using \eqref{Vfjf} for the case of  $1<p<2$ we arrive at the following 
    \begin{align*}
       \int_{\mathcal{B}}&\left \lvert \left(d\widetilde{W}_j+\widetilde{F}_j\right)- \left(d\widetilde{W}
        +\widetilde{F}\right)\right \rvert^p\\
        &\leq \epsilon\int_{\mathcal{B}}\left\lvert d\widetilde{W}+\widetilde{F}\right\rvert^{p}+\frac{c}{\epsilon}\left(\int_{B_r\left(0\right)} \left[\left \lvert \widetilde{F}_j\right\rvert^p +\left \lvert d\widetilde{W}\right \rvert^p+\left \lvert \widetilde{F}\right \rvert\right]^{p}\right)^{\frac{p-1}{p}}\left(\int_{B_r\left(0\right)}\left \lvert \widetilde{F}-\widetilde{F}_j\right \rvert^p\right)^{\frac{1}{p}}.
    \end{align*}
    If $\mathcal{B}_1:=\left \{x\in \mathcal{B}: \left \lvert \left(d\widetilde{W}_j+\widetilde{F}_j\right)- \left(d\widetilde{W}
        +\widetilde{F}\right)\right \rvert\neq 0\right\}$ and $\mathcal{B}_2:=\left \{x\in \mathcal{B}: \left\lvert d\widetilde{W}+\widetilde{F}\right\rvert\neq 0\right\}$ both are of positive measure then we define,
    \begin{align*}\epsilon=\frac{1}{2}\times\frac{\int_{\mathcal{B}}\left \lvert \left(d\widetilde{W}_j+\widetilde{F}_j\right)- \left(d\widetilde{W}
        +\widetilde{F}\right)\right \rvert^p}{\int_{\mathcal{B}}\left\lvert d\widetilde{W}+\widetilde{F}\right\rvert^{p}}.
        \end{align*}
        and hence
    \begin{align*}
       \int_{\mathcal{B}}\left \lvert \left(d\widetilde{W}_j+\widetilde{F}_j\right)- \left(d\widetilde{W}
        +\widetilde{F}\right)\right \rvert^p
        &\begin{aligned}[t]\leq &\frac{1}{2}\int_{\mathcal{B}}\left \lvert \left(d\widetilde{W}_j+\widetilde{F}_j\right)- \left(d\widetilde{W}
        +\widetilde{F}\right)\right \rvert^p\\
        &+ \frac{c\int_{\mathcal{B}}\left\lvert d\widetilde{W}+\widetilde{F}\right\rvert^{p}}{\int_{\mathcal{B}}\left \lvert \left(d\widetilde{W}_j+\widetilde{F}_j\right)- \left(d\widetilde{W}
        +\widetilde{F}\right)\right \rvert^p} \left(\int_{B_r\left(0\right)}\left \lvert \widetilde{F}-\widetilde{F}_j\right \rvert^p\right)^{\frac{1}{p}}\\
        &\qquad\qquad\qquad\qquad\times \left(\int_{B_r\left(0\right)} \left[\left \lvert \widetilde{F}_j\right\rvert^p +\left \lvert d\widetilde{W}\right \rvert^p+\left \lvert \widetilde{F}\right \rvert\right]^{p}\right)^{\frac{p-1}{p}}.
        \end{aligned}
    \end{align*}
    Thus
    \begin{align*}
       &\left(\int_{\mathcal{B}}\left \lvert \left(d\widetilde{W}_j+\widetilde{F}_j\right)- \left(d\widetilde{W}
        +\widetilde{F}\right)\right \rvert^p\right)^{2}\\
        &\qquad\lesssim  {\int_{\mathcal{B}}\left\lvert d\widetilde{W}+\widetilde{F}\right\rvert^{p}}\left(\int_{B_r\left(0\right)} \left[\left \lvert \widetilde{F}_j\right\rvert^p +\left \lvert d\widetilde{W}\right \rvert^p+\left \lvert \widetilde{F}\right \rvert\right]^{p}\right)^{\frac{p-1}{p}} \left(\int_{B_r\left(0\right)}\left \lvert \widetilde{F}-\widetilde{F}_j\right \rvert^p\right)^{\frac{1}{p}}
        \end{align*}
    i.e.
    \begin{align*}
       \int_{\mathcal{B}}&\left \lvert \left(d\widetilde{W}_j+\widetilde{F}_j\right)- \left(d\widetilde{W}
        +\widetilde{F}\right)\right \rvert^p\\
        &\lesssim \begin{aligned}[t]\left[\int_{B_r\left(0\right)}\left\lvert d\widetilde{W}+\widetilde{F}\right\rvert^{p}\right]^{\frac{1}{2}}&\left(\int_{B_r\left(0\right)} \left[\left \lvert \widetilde{F}_j\right\rvert^p +\left \lvert d\widetilde{W}\right \rvert^p+\left \lvert \widetilde{F}\right \rvert\right]^{2p}\right)^{\frac{p-1}{2p}} \left(\int_{B_r\left(0\right)}\left \lvert \widetilde{F}-\widetilde{F}_j\right \rvert^p\right)^{\frac{1}{2p}}.
        \end{aligned}
    \end{align*} 
    The last estimate is trivially true if at least one of $\left \lvert \mathcal{B}_1\right \rvert$ and $\left \lvert \mathcal{B}_1\right \rvert$ is zero. Finally,
    \begin{align*}
       &\int_{B_r\left(0\right)=\mathcal{A}\sqcup\mathcal{B}}\left \lvert \left(d\widetilde{W}_j+\widetilde{F}_j\right)- \left(d\widetilde{W}
        +\widetilde{F}\right)\right \rvert^p\\
        &\quad\lesssim \begin{aligned}[t]&\left(\int_{B_r\left(0\right)}\left[\left \lvert d\widetilde{W}\right \rvert^p+\left \lvert \widetilde{F}_j\right \rvert^p+ \left \lvert \widetilde{F}\right \rvert^p\right]\right )^{\frac{p-1}{p}}\left(\int_{B_r\left(0\right)}\left \lvert \widetilde{F}_j-\widetilde{F}\right \rvert^p\right)^{\frac{1}{p}}\\
        &+\left[\int_{B_r\left(0\right)}\left\lvert d\widetilde{W}+\widetilde{F}\right\rvert^{p}\right]^{\frac{1}{2}}\left(\int_{B_r\left(0\right)} \left[\left \lvert \widetilde{F}_j\right\rvert^p +\left \lvert d\widetilde{W}\right \rvert^p+\left \lvert \widetilde{F}\right \rvert\right]^{2p}\right)^{\frac{p-1}{2p}}\left(\int_{B_r\left(0\right)}\left \lvert \widetilde{F}-\widetilde{F}_j\right \rvert^p\right)^{\frac{1}{2p}}.\\
        \end{aligned}
    \end{align*}
    which implies 
        $d\widetilde{W}_j+\widetilde{F}_j \rightarrow d\widetilde{W}+\widetilde{F} \text{ in } L^p\left(B_r\left(0\right)\right)$ and
    after relabeling a subsequence we have
    $d\widetilde{W}_j+\widetilde{F}_j \rightarrow d\widetilde{W}+\widetilde{F}        \text{ a.e in  } B_r\left(0\right)$.
    \end{proof}
    Using \eqref{smooth estimate j} and Fatou's Lemma
    \begin{align*}
        \int_{Q_0}\left \lvert d\widetilde{W}+\widetilde{F}\right \rvert^{\widetilde{q}}& \leq \operatorname{lim inf}\int_{Q_0}\left \lvert d\widetilde{W}_j+\widetilde{F}_j\right \rvert^{\widetilde{q}}\\
        &\leq c\left(n,\tilde{q},p,k,r\right) \liminf_{j\rightarrow \infty} \left[\left[\int_{Q_0}\left \lvert d\widetilde{W}_j+\widetilde{F}_j \right \rvert^p\right]^{\frac{\tilde{q}}{p}}
        +\int_{Q_0}\left \lvert \widetilde{F}_j\right\rvert^{\tilde{q}}\right]\\
        &= c\left(n,\tilde{q},p,k,r\right) \lim_{j\rightarrow \infty} \left[\left[\int_{Q_0}\left \lvert d\widetilde{W}_j+\widetilde{F}_j \right \rvert^p\right]^{\frac{\tilde{q}}{p}}
        +\int_{Q_0}\left \lvert \widetilde{F}_j\right\rvert^{\tilde{q}}\right]\\
        &= c\left(n,\tilde{q},p,k,r\right) \left[\left[\int_{Q_0}\left \lvert d\widetilde{W}+\widetilde{F} \right \rvert^p\right]^{\frac{\tilde{q}}{p}}
        +\int_{Q_0}\left \lvert \widetilde{F}\right\rvert^{\tilde{q}}\right].
    \end{align*}
    Also, $\left \lvert Q_0\right \rvert$  depends only on $n,p,\widetilde{q},k$ and $r$.
    \end{proof}
    \smallskip

    \section{Proof of Theorem \ref{Goal Theorem}}
    
    We have $u\in W^{1,\widetilde{q}}\left(B;\Lambda^k\right)$ where $0\leq k\leq n-1,1<p<\widetilde{q}<\infty$ with $\widetilde{q}$ satisfying \eqref{qtilde}.
    $v\in u+W^{1,p}_{d^*,T}\left(B;\Lambda^{k}\right)$ is the unique minimizer of the functional 
            $\overline{w}\mapsto \int_{B} \left \lvert d\overline{w} \right \rvert^p$ where $\overline{w} \in u+W^{1,p}_{d^*,T}\left(B;\Lambda^{k}\right)$. Thus 
 \begin{align*}
        \int_{B}\left \langle \left \lvert dv\right \rvert^{p-2}dv,d\phi \right \rangle=0 \text{ for $\phi \in W^{1,p}_{d^*,T}\left(B;\Lambda^k\right)$}.
    \end{align*}    
     Define $\widetilde{w}:=v-u\in W^{1,p}_{d^*,T}\left(B;\Lambda^k\right),f:=du$; thus     
    \begin{align*}
        \int_{B}\left \langle \left \lvert d\widetilde{w}+f\right \rvert^{p-2}\left(d\widetilde{w}+f\right),d\phi \right \rangle=0 \text{ for } \phi \in W^{1,p}_{d^*,T}\left(B;\Lambda^k\right).
    \end{align*} 
     Now, $d\widetilde{w}\in L^p\left(B, \Lambda^{k+1}\left(\mathbb{R}^n\right)\right)$ with $d\left(d\widetilde{w}\right)=0$ in $B$ in distributional sense and since $\nu \wedge \widetilde{w}=0$ on $\partial B$ we have  $\nu \wedge d\widetilde{w}=0$ on $\partial B$, by Lemma \ref{poincare} there exists $w\in W^{1,p}_0\left(B, \Lambda^k\left(\mathbb{R}^n\right)\right)$ such that $dw=d\widetilde{w}$ in $ in B$.
    Hence 
    \begin{align} \label{el d*T}
         \int_{B}\left \langle \left \lvert d{w}+f\right \rvert^{p-2}\left(d{w}+f\right),d\phi \right \rangle=0 \text{ for all } \phi\in W^{1,p}_{d^*,T}\left(B;\Lambda^k\right).
    \end{align}
    Again Lemma \ref{poincare} and Lemma \ref{1-div-curl} we get the following equivalent statement:
    \begin{align} \label{el 0}
         \int_{B}\left \langle \left \lvert d{w}+f\right \rvert^{p-2}\left(d{w}+f\right),d\phi \right \rangle=0 \text{ for all } \phi\in W^{1,p}_{0}\left(B;\Lambda^k\right).
    \end{align}
    Since, $C_c^\infty\left(B;\Lambda^k\right)$ is a dense subspace of $W^{1,p}_{0}\left(B;\Lambda^k\right)$ we derive another equivalent weak formulation: 
    \begin{align} \label{el CC}
         \int_{B}\left \langle \left \lvert d{w}+f\right \rvert^{p-2}\left(d{w}+f\right),d\phi \right \rangle=0 \text{ for all } \phi\in C_c^\infty\left(B;\Lambda^k\right).
    \end{align}
We divide our proof into three subsections. By the classical regularity result (see \eqref{sup bound}) we already have the local higher integrability of $dv$ in $B$. To achieve a global higher integrability estimates we first flatten out the a small part of the boundary and via reflection trick we extend the Euler-Lagrange equation in a neighborhood of the boundary $\partial B$ and in the final section we complete the proof of Theorem \ref{Goal Theorem}. We explain the steps in details below: 
\subsection{Flattening out the boundary}\label{Flattening out the boundary}
Fix, $x^0 \in \partial B$ and take $U_{\widetilde{r}}:=B_{\widetilde{r}}\left(x_0\right)\cap B$. Without loss we take $B=B_1\left(e_n\right),x^0=0$ where $e_n=\left(0,\cdots,0,1\right)\in \mathbb{R}^n$ and proceed to  flattening of the boundary.
    For flattening out the boundary we shall use the method mentioned in \cite[APPENDIX C.1.]{evans2010partial}.\\
    
    For $x\in \mathbb{R}^n$ we write $x=\left(x',x_n\right)$ where $x'=\left(x_1,x_2,\cdots,x_{n-1}\right)\in \mathbb{R}^{n-1}$ and whenever $\left \lvert x'\right \rvert^2=\sum_{j=1}^{n-1}x_j^2\leq 1$ we define $\gamma\left(x'\right)=\gamma\left(x_1,...,x_{n-1}\right):=1-\sqrt{1-\left \lvert x'\right \rvert^2}$.\\
    Let $0<\widetilde{r}<\frac{1}{1000
    }$; we define $\Phi:B_{\widetilde{r}}\left(0\right)\cap B_1\left(e^n\right)\rightarrow \mathbb{R}^n_+$ by
    \begin{align}\label{phi}
        \Phi\left(x\right):&=\left(\Phi_1\left(x\right),\cdots,\Phi_{n-1}\left(x\right),\Phi_n\left(x\right)\right)=\left(x', x_n-\gamma\left(x'\right)\right).
    \end{align}
Clearly $\Phi$ is an invertible matrix define its inverse by $\Psi$. Also, $\operatorname{det}\left(D\Phi\right)=1$ in $B_{\widetilde{r}}\left(0\right)\cap B_1\left(e^n\right)$ and 
\begin{lemma}\label{Phi Lemma after def.}
    \begin{align}\label{del phi i del x j converging to zero}
\left \lVert\delta_{ij}-\frac{\partial \Phi_j}{\partial x_i}\right \rVert_{L^\infty\left(B_{\widetilde{r}}\left(0\right)\cap B\right)}\rightarrow 0 \text{ as } \widetilde{r}\rightarrow 0^+. 
    \end{align}   
    \end{lemma}
    \begin{proof}
    By direct computations we get,
    \begin{align*}
        &\frac{\partial \Phi_j}{\partial x_i}=\delta_{ij} \text{ for } j\neq n, \frac{\partial \Phi_n}{\partial x_i}=\frac{-x_i}{\sqrt{1-\left \lvert x'\right \rvert^2}}dx_i \text{ for } i\neq n \text{ and } \frac{\partial \Phi_n}{\partial x_n}=1.
    \end{align*}
    For $l\neq n$ and $x\in B_{\widetilde{r}}\left(0\right)\cap B_1\left(e^n\right)$ we have
    $\left \lvert \frac{x_l}{\sqrt{1-\left(\sum_{j=1}^{n-1} x_j^2\right)}}\right \rvert \leq \frac{\widetilde{r}}{\sqrt{1-\widetilde{r}^2}}$ and \eqref{del phi i del x j converging to zero} follows.
    \end{proof}
    
    Notice that, $\Phi\left(B_{\widetilde{r}}\left(0\right)\cap B\right)\subset \mathbb{R}^n_+$. 
    Fix small enough $r$ (depending on $\widetilde{r}$) with $B_r^+\left(0\right)\subset \Phi\left(B_{\widetilde{r}}\left(0\right)\cap B_1\left(e^n\right)\right)$. 
    Also notice that,
    \begin{align}\label{tilde r and r}
        r\rightarrow 0+ \text{ whenever } \widetilde{r} \rightarrow 0+
    \end{align}We define
    \begin{align}
    \left \lbrace
    \begin{aligned}&W:=\Psi^*\left(w\right),F:=\Psi^*\left(f\right)
    \text{ in } \Phi\left(B_{\widetilde{r}}\left(0\right)\cap B_1\left(e^n\right)\right),\\
    &\text{ i.e. } w=\Phi^*\left(W\right), f=\Phi^*\left(F\right) \text{ in } B_{\widetilde{r}}\left(0\right)\cap B_1\left(e^n\right).
    \end{aligned}
    \right.
    \end{align}
    One can also notice that that $W=0$ on $\overline{B^+_{{r}}\left(0\right)}\cap \left\{x_n=0\right\}$ in trace sense.\\
    Choose, $\eta \in C_c^\infty\left(B^+_r\left(0\right), \Lambda^k\right)$, then $\Phi^*\left(\eta\right)\in C^\infty_c\left(B_{\widetilde{r}}\left(x^0\right)\cap B_1\left(e_n\right); \Lambda^k\right)$ and hence by \eqref{el CC} we have
    \begin{align*}
        \int_{B_{\widetilde{r}}\left(x^0\right)\cap B}\left \langle \left[\left \lvert \Phi^*\left(dW+F\right)\right \rvert^{p-2}\Phi^*\left( dW+F\right)\right], \Phi^*\left(d\eta\right)\right \rangle=0.
    \end{align*}
    Change of variable gives 
    \begin{align}\label{flattening section last eqn}
        \int_{ B^+_r\left(0\right)}\left \langle \left \lvert \left[ \Phi^*\left(dW+F\right)\right]\circ \Psi\right \rvert^{p-2}\left[\Phi^*\left( dW+F\right)\right]\circ \Psi, \left[\Phi^*\left(d\eta\right)\right]\circ \Psi\right \rangle=0,
    \end{align}
    for each $\eta \in C^\infty_c\left(B^+_r\left(0\right)\right)$.\\
 (For the notations like $\Phi^*\left(\mathcal{G}\right)$ and $\mathcal{G}\circ \Phi$ see Definition \ref{Pullback definition}).
    \subsection{Reflection}

    For $\xi =\sum_J \xi_Jdx^J\in \Lambda^{k+1}$ and $x\in B^+_{{r}}\left(0\right)$ we define
    \begin{align}\label{TPhi def}
        T_{\Phi}\left(x\right) \xi
        := \sum_{1\leq i_1<...<i_{k+1}\leq n}\xi_{i_1...i_{k+1}} \left(d\Phi^{i_1}\circ \Psi\right)\left(x\right)\wedge... \wedge \left(d\Phi^{i_{k+1}}\circ \Psi\right)\left(x\right)
    \end{align}
    Hence for $x\in B^+_{{r}}\left(0\right), T_\Phi\left(x\right)$ is understood as an element of $ \mathcal{L}\left(\Lambda^{k+1};\Lambda^{k+1}\right)$.
    Now for $k+1$-differential form $\mathcal{G}$ with domain in $B^+_r\left(0\right)$ we have,
    \begin{align}\label{Tphi g}
        &\left[\left(\Phi^*\left(\mathcal{G}\right)\right)\circ \Psi\right]\left(x\right)\notag\\
        &=\left[\left(\sum_{1\leq i_1<...<i_{k+1}\leq n}\left(\mathcal{G}_{i_1...i_{k+1}}\circ \Phi\right)  d\Phi^{i_1}\wedge... \wedge d\Phi^{i_{k+1}}\right)\circ \Psi\right]\left(x\right)\notag\\
        &=\sum_{1\leq i_1<...<i_{k+1}\leq n}\mathcal{G}_{i_1...i_k}\left(x\right) \left(d\Phi^{i_1}\circ\Psi\right)\left(x\right)\wedge\cdots \wedge \left(d\Phi^{i_{k+1}}\circ\Psi\right)\left(x\right) \notag\\
        &=T_\Phi\left(x\right)\mathcal{G} \left(x\right).
    \end{align}
    Hence from \eqref{flattening section last eqn} and \eqref{Tphi g} we get
    \begin{align}\label{Gphi + phi star}
        0=&\int_{ B^+_r\left(0\right)}\left \langle \left \lvert \left[ \Phi^*\left(dW+F\right)\right]\circ \Psi\right \rvert^{p-2}\left[\Phi^*\left( dW+F\right)\right]\circ \Psi, \left[\Phi^*\left(d\eta\right)\right]\circ \Psi\right \rangle\notag \\
        &=\int_{ B^+_r\left(0\right)}\left \langle \left \langle \left[T_\Phi\right]^T T_\Phi\left(dW+F\right), dW+F\right \rangle^{\frac{p-2}{2}}\left[T_\Phi\right]^T T_\Phi\left(dW+F\right), d\eta \right \rangle\notag \\
        &=\int_{ B^+_r\left(0\right)}\left \langle \left \langle G^+_\Phi\left(d{W}+{F}\right), d{W}+{F}\right \rangle^{\frac{p-2}{2}}G^+_\Phi\left(d{W}+{F}\right), d\eta \right \rangle,
    \end{align}
    where $G^+_\Phi:=\left[T_\Phi\right]^T T_\Phi$ for $x\in B^+_r\left(0\right)$.
    Define $S\left(x',x_n\right):=\left(x',-x_n\right)$ and 
    \begin{align*}
        \widetilde{W}:=\begin{cases}
            W &\text{ in } B^+_r\left(0\right),\\
            -S^*\left(W\right) &\text{ in } B^-_r\left(0\right).
        \end{cases}
    \end{align*}
    Thus for $x$ in $B^{-}_r\left(0\right)$ we have
    \begin{align*}
        \widetilde{W}\left(x\right)=&-\sum_{1\leq i_1< \cdots< i_{k+1}<n}W_{i_1\cdots i_{k+1}}\left(x',-x_{n}\right)dx^{i_1}\wedge \cdots dx^{i_{k+1}}\\
        &\quad+\sum_{1\leq i_1< \cdots< i_{k}<n}W_{i_1\cdots i_{k}n}\left(x',-x_{n}\right)dx^{i_1}\wedge \cdots dx^{i_{k}}\wedge dx^{n}.
    \end{align*}
     Clearly, $\widetilde{W}\in W^{1,p}\left(B_r\left(0\right)\right)$ since $W=0$ on $\left \{x_n=0\right\}\cap \overline{B_r^+\left(0\right)}$ in trace sense. 
   Similarly we define
    \begin{align*}
        \widetilde{F}:=\begin{cases}
            F &\text{ in } B^+_r\left(0\right),\\
            -S^*\left(F\right) &\text{ in } B^-_r\left(0\right).
        \end{cases}
    \end{align*}
    For $x\in B_r^-\left(0\right)$ and $\xi\in \Lambda^{k+1}$ define
    \begin{align}\label{RPhi definition}
        R_{\Phi\left(x\right)}&\xi \notag \\
        &=\begin{aligned}[t]
            &\sum_{1\leq i_1<...<i_{k+1}<n}\xi_{i_1...i_{k+1}} dx^{i_1}\wedge \cdots \wedge dx^{i_{k+1}}\\
            &\qquad \qquad-\sum_{1\leq i_1<...<i_{k}<n}\xi_{i_1...i_{k}n} dx^{i_1}\wedge \cdots \wedge dx^{i_{k}}\wedge \left[d\Phi_n \left(\Psi\left(S\left(x\right)\right)\right)\right].
        \end{aligned}
    \end{align}
    Comparing with the case of $T_{\Phi}$ one can see that 
      for $x\in B^-_{{r}}\left(0\right)$, $R_\Phi\left(x\right)\in \mathcal{L}\left(\Lambda^{k+1};\Lambda^{k+1}\right)$.
    Hence for a $k+1$-differential form $\mathcal{G'}$ with domain $B^-_r\left(0\right)$ we have,
    \begin{align}\label{RPhi def}
        &\left[\left[\Phi^*S^*\left(\mathcal{G'}\right)\right]\circ \Psi\right]\circ S\notag \\
        &=\sum_{i_1<...<i_{k+1}}\mathcal{G'}_{i_1...i_{k+1}} \left[d\left(S\circ\Phi\right)^{i_1}\circ \left(\Psi \circ S\right)\right]\wedge... \wedge \left[d\left(S\circ\Phi\right)^{i_{k+1}}\circ \left(\Psi \circ S\right)\right]\notag \\
        &=\begin{aligned}[t]
            &\sum_{1\leq i_1<...<i_{k+1}<n}\mathcal{G}'_{i_1...i_{k+1}} dx^{i_1}\wedge \cdots \wedge dx^{i_{k+1}}\notag \\
            &\qquad \qquad-\sum_{1\leq i_1<...<i_{k}<n}\mathcal{G}'_{i_1...i_{k}n} dx^{i_1}\wedge \cdots \wedge dx^{i_{k}}\wedge \left[d\Phi_n \left(\Psi\left(S\left(x\right)\right)\right)\right]
        \end{aligned}\\&=R_{\Phi\left(x\right)}\mathcal{G}'\left(x\right).
    \end{align}
    Again for $\eta \in C^\infty_c\left(B^-_r\left(0\right)\right)$ we have the following
    \begin{align}\label{B- 0}
    \int_{ B^-_r\left(0\right)}&\left \langle \left \langle R_\Phi^T R_\Phi\left(d\widetilde{W}+\widetilde{F}\right), d\widetilde{W}+\widetilde{F}\right \rangle^{\frac{p-2}{2}}R_\Phi^T R_\Phi\left(d\widetilde{W}+\widetilde{F}\right), d\eta\right \rangle\notag \\
        &\begin{aligned}[t]=\int_{ B^-_r\left(0\right)}\left \langle \left \lvert R_\Phi\left(d\widetilde{W}+\widetilde{F}\right)\right \rvert^{p-2}R_\Phi\left(d\widetilde{W}+\widetilde{F}\right), R_\Phi\left(d\eta\right)\right \rangle
        \end{aligned}\notag\\
        &\begin{aligned}[t]=\int_{ B^-_r\left(0\right)}&\left \lvert \left[\Phi^*S^*\left(d\widetilde{W}+\widetilde{F}\right)\right]\circ \Psi \circ S\right \rvert^{p-2}\\
        &\left \langle \left[\Phi^*S^*\left(d\widetilde{W}+\widetilde{F}\right)\right]\circ \Psi \circ S,\notag  d\left( \Phi^*\left(S^*\eta\right)\right)\circ \Psi\circ S\right \rangle \text{[see \eqref{RPhi def}] }
        \end{aligned}\notag\\
        &\begin{aligned}[t]=-\int_{ S\left(B^+_r\left(0\right)\right)}\left \langle \left \lvert \Phi^*\left(d{W}+{F}\right)\circ \Psi \circ S\right \rvert^{p-2}\right.&\Phi^*\left(d{W}+{F}\right)\circ \Psi \circ S,\notag\\
        &\left. d\left( \Phi^*\left(S^*\eta\right)\right)\circ \Psi\circ S\right \rangle
        \end{aligned}\notag\\
        &\begin{aligned}[t]=-&\int_{ B^+_r\left(0\right)}\left \lvert \Phi^*\left(d{W}+{F}\right)\circ \Psi \right \rvert^{p-2}\\
        &\qquad \times \left \langle \Phi^*\left(d{W}+{F}\right)\circ \Psi, \Phi^*\left(d\left(S^*\eta\right)\right)\circ \Psi\right \rangle=0\text{ [by \eqref{flattening section last eqn}]}.
        \end{aligned}
    \end{align}
    We define $G^-_\Phi\left(x\right):=\left[R_\Phi\right]^T R_\Phi$ for $x\in B^-_r\left(0\right)$.
     So for $\eta \in C^\infty_c\left(B^-_r\left(0\right);\Lambda^k\right)$
    \begin{align}\label{Gphi - phi star}
    &\int_{ B^-_r\left(0\right)}\left \langle \left \langle G^-_\Phi\left(x\right)\left(d\widetilde{W}+\widetilde{F}\right), d\widetilde{W}+\widetilde{F}\right \rangle^{\frac{p-2}{2}}G^-_\Phi\left(x\right)\left(d\widetilde{W}+\widetilde{F}\right), d\eta\right \rangle\notag \\
        &\quad \begin{aligned}[t]=-\int_{ B^+_r\left(0\right)}\left \langle \left \lvert \Phi^*\left(d{W}+{F}\right)\circ \Psi \right \rvert^{p-2}\right. &\Phi^*\left(d{W}+{F}\right)\circ \Psi,\left. \Phi^*\left(d\left(S^*\eta\right)\right)\circ \Psi\right \rangle=0.
        \end{aligned}
    \end{align}
    Define
    \begin{align}\label{Gphi}
        G_{\Phi}:=\begin{cases}
            G_\Phi^+ & \text{ in } B^+_r\left(0\right)\\
            G_\Phi^- & \text{ in } B^-_r\left(0\right)
        \end{cases}
    \end{align}
    \begin{lemma}\label{Extension of G_Phi eqn.}
        For each $\eta \in C_c^\infty\left(B_r\left(0\right)\right)$ we have the following
        \begin{align}\label{Gphi extension in eqn.}
        \int_{B_r\left(0\right)}\left\langle \left \langle G_\Phi\left(d\widetilde{W}+\widetilde{F}\right), d\widetilde{W}+\widetilde{F}\right \rangle^{\frac{p-2}{2}}G_\Phi\left(d\widetilde{W}+\widetilde{F}\right),d\eta \right \rangle=0 .
    \end{align}
    \end{lemma}
    \begin{proof}
    Choose a smooth function $\zeta : \mathbb{R}\rightarrow [0,1]$ such that
        $\zeta \equiv 1 $ in $\mathbb{R}\setminus [-2,2],\operatorname{supp}\left(\zeta\right) \subseteq \mathbb{R}\setminus [-1,1]$ and $\zeta\left(t\right)=\zeta\left(-t\right)$ for all $t \in \mathbb{R}$.\\ 
        For $\eta \in C_c^\infty\left(B_r\left(0\right)\right)$ we define $\zeta_m\left(t\right)=\zeta\left(mt\right)$ and $ \eta_m\left(x\right)=\eta\left(x\right)\zeta_m\left(x_n \right)$ where $x=\left(x_1,\cdots, x_n\right)$, thus
    \begin{align*}
        0&=\int_{B_r\left(0\right)} \left \langle \left \langle G_\Phi\left(d\widetilde{W}+\widetilde{F}\right), d\widetilde{W}+\widetilde{F}\right \rangle^{\frac{p-2}{2}}G_\Phi\left(d\widetilde{W}+\widetilde{F}\right),d\eta_m\right \rangle\\
        &=\int_{B_r\left(0\right)} \left \langle \left \langle G_\Phi\left(d\widetilde{W}+\widetilde{F}\right), d\widetilde{W}+\widetilde{F}\right \rangle^{\frac{p-2}{2}}G_\Phi\left(d\widetilde{W}+\widetilde{F}\right),d\eta\right \rangle\zeta_m\left(x_n\right)\\
        &\quad+\int_{B_r\left(0\right)} \left \langle \left \langle G_\Phi\left(d\widetilde{W}+\widetilde{F}\right), d\widetilde{W}+\widetilde{F}\right \rangle^{\frac{p-2}{2}}G_\Phi\left(d\widetilde{W}+\widetilde{F}\right),dx_n\wedge \eta\right \rangle\zeta'_m\left(x_n\right).
    \end{align*}
    Define, 
    \begin{align*}
        &K^{\left(m\right)}_1:=\int_{B_r\left(0\right)} \left \langle \left \langle G_\Phi\left(d\widetilde{W}+\widetilde{F}\right), d\widetilde{W}+\widetilde{F}\right \rangle^{\frac{p-2}{2}}G_\Phi\left(d\widetilde{W}+\widetilde{F}\right),d\eta\right \rangle\zeta_m\left(x_n\right),\\
        &K^{\left(m\right)}_2:= \int_{B_r\left(0\right)} \left \langle \left \langle G_\Phi\left(d\widetilde{W}+\widetilde{F}\right), d\widetilde{W}+\widetilde{F}\right \rangle^{\frac{p-2}{2}}G_\Phi\left(d\widetilde{W}+\widetilde{F}\right),dx_n\wedge \eta\right \rangle\zeta'_m\left(x_n\right).
    \end{align*}
    Now $\zeta_m \equiv 1 $ in $\mathbb{R}\setminus [-\frac{2}{m},\frac{2}{k}],\operatorname{supp}\left(\zeta\right) \subseteq \mathbb{R}\setminus [-\frac{1}{m},\frac{1}{m}]$
    hence using DCT we get that
    \begin{align*} K^{\left(m\right)}_1&\rightarrow \int_{B_r\left(0\right)} \left \langle \left \langle G_\Phi\left(d\widetilde{W}+\widetilde{F}\right), d\widetilde{W}+\widetilde{F}\right \rangle^{\frac{p-2}{2}}G_\Phi\left(d\widetilde{W}+\widetilde{F}\right),d\eta\right \rangle  \text{ as } m \rightarrow \infty.
    \end{align*}
    On the other hand by straight forward calculation
    \begin{align}\label{zeta prime}
    \left \lbrace
    \begin{aligned}
        K^{\left(m\right)}_2=&\begin{aligned}[t]&\int_{B^+_r\left(0\right)}\left \langle \left \lvert \left[\Phi^*\left(d{W}+{F}\right)\right]\circ \Psi\right \rvert^{p-2} \left[\Phi^*\left(d{W}+{F}\right)\right]\circ \Psi,\right.& \\
        &\left.\qquad\qquad \left[\Phi^*\left(dx_n\wedge \eta\right)\right]\circ \Psi+
        \left[\Phi^*S^*\left(dx_n\wedge \eta\right)\right]\circ \Psi\right \rangle \zeta'_m\left(x_n\right).
        \end{aligned}
        \end{aligned}
        \right.
    \end{align}
   Notice that $\left(S\circ \Phi\right)_n=-\Phi^n$; thus
    \begin{align*}
        &\left[\left[\Phi^*\left(dx_n\wedge \eta\right)\right]\circ \Psi\right]\left(x\right)=\sum_{1\leq i_1<...<i_k\leq n-1} \eta_{i_1...i_k} \left(x\right)\left[d\Phi^n \circ \Psi\right]\wedge dx^{i_1}\wedge...\wedge dx^{i_k};\\
        &\begin{aligned}[t]&\left[\left[\Phi^*S^*\left(dx_n\wedge \eta\right)\right]\circ \Psi\right]\left(x\right)\\
        &\qquad\qquad=-\sum_{1\leq i_1<...<i_k\leq n-1} \eta_{i_1...i_k} \left(x',-x_n\right)\left[d\Phi^n \circ \Psi\right]\wedge dx^{i_1}\wedge...\wedge dx^{i_k}.
        \end{aligned}
    \end{align*}
    Now, $\zeta'_m\left(x_n\right)=m \zeta'\left(mx_m\right)$, on the other hand
    $\left \lvert \eta \left(x',x_n\right) -\eta\left(x',-x_n\right) \right \rvert \leq 2\left \lVert\nabla \eta \right \rVert_{L^\infty\left(B_r\left(0\right)\right)}\left \lvert x_n\right \rvert$.
    Hence, by \eqref{zeta prime}  letting $m \rightarrow \infty$ we have
    \begin{multline*}
        \left \lvert K^{\left(m\right)}_2\right \rvert 
        \lesssim m\left \lVert\nabla \eta \right\rVert_{L^\infty\left(B_r\left(0\right)\right)}\left\lVert\zeta'\right\rVert_{L^\infty\left(\mathbb{R}\right)}\int_{B^+_r\left(0\right)\cap \left\{x_n \leq \frac{1}{m}\right\}} \left \lvert \left[\Phi^*\left(d{W}+{F}\right)\right]\circ \Psi\right \rvert^{p-1} \left \lvert x_n\right \rvert\\
        \leq \left \lVert\nabla \eta \right\rVert_{L^\infty\left(B_r\left(0\right)\right)}\left\lVert\zeta'\right\rVert_{L^\infty\left(\mathbb{R}\right)}\int_{B^+_r\left(0\right)\cap \left\{x_n \leq \frac{1}{m}\right\}} \left \lvert \left[\Phi^*\left(d{W}+{F}\right)\right]\circ \Psi\right \rvert^{p-1}\rightarrow 0.
    \end{multline*}
    and finally we get
    \begin{align*}
        \int_{B_r\left(0\right)}\left\langle \left \langle G_\Phi\left(d\widetilde{W}+\widetilde{F}\right), d\widetilde{W}+\widetilde{F}\right \rangle^{\frac{p-2}{2}}G_\Phi\left(d\widetilde{W}+\widetilde{F}\right),d\eta \right \rangle=0. 
    \end{align*}
    \end{proof}
    To use Theorem \ref{G Theorem} we need \eqref{r} and \eqref{1<p<2 and p<2lambda first mention} of $*$-condition i.e. we need $\left \lVert G_\Phi-I \right \rVert_{L^\infty\left(B_r\left(0\right)\right) }$ small enough for suitably fixed small $r>0$; which is guaranteed by the following lemma:
        \begin{lemma}\label{* condition Lemma}
            There exists $\widetilde{r}>0$ small enough such that $*$-condition holds; (see Subsection \ref{Flattening out the boundary} for $\widetilde{r}$ and see \eqref{delta},\eqref{r},\eqref{1<p<2 and p<2lambda first mention} for $*$-condition).
        \end{lemma}
        \begin{proof}
         Using \eqref{del phi i del x j converging to zero} we have 
    \begin{align*}
        \left \lVert d\Phi^i-dx^i\right \rVert_{L^{\infty}\left(B_{\widetilde{r}}\left(0\right)\cap B\right)}&=\left \lVert \sum_{j=1}^n\frac{\partial \Phi^i}{\partial x_j}dx^j-dx^i\right \rVert_{L^{\infty}\left(B_{\widetilde{r}}\left(0\right)\cap B\right)}\\
        &=\left \lVert \sum_{j=1}^n\left(\frac{\partial \Phi^i}{\partial x_j}-\delta_{ij}\right)dx^j\right \rVert_{L^{\infty}\left(B_{\widetilde{r}}\left(0\right)\cap B\right)}\\
        &\leq \sum_{j=1}^n\left \lVert \frac{\partial \Phi^i}{\partial x_j}-\delta_{ij}\right \rVert_{L^{\infty}\left(B_{\widetilde{r}}\left(0\right)\cap B\right)}\rightarrow 0 \qquad\text{ as } \widetilde{r}\rightarrow 0^+.
    \end{align*}
    Therefore
    \begin{align}\label{dPhi convergence}
        \left \lVert d\Phi^i\circ \Psi-dx^i\right \rVert_{L^{\infty}\left(B_r^+\left(0\right)\right)}=\left \lVert d\Phi^i\circ \Psi-dx^i\circ \Psi\right \rVert_{L^{\infty}\left(B_r^+\left(0\right)\right)}\rightarrow 0 \text{ as } \widetilde{r}\rightarrow 0^+,
    \end{align}
     and by \eqref{Tphi g},
    $\left \lVert T_{\Phi}- I\right \rVert_{L^{\infty}\left(B_r^+\left(0\right)\right)}$ as $\widetilde{r} \rightarrow 0^+$. Recall that $G^+_\Phi:=\left[T_\Phi\right]^T T_\Phi$ for $x\in B^+_r\left(0\right)$
    and hence $\left \lVert G^+_{\Phi}- I\right \rVert_{L^{\infty}\left(B_r^+\left(0\right)\right)}$ as $\widetilde{r} \rightarrow 0^+$.\\
    On the other hand by \eqref{dPhi convergence}
    \begin{align*}
        &\left \lVert d\left(S\circ\Phi\right)^{n}\circ \left(\Psi \circ S\right)+ dx^n\right \rVert_{L^\infty\left(B^-_r\left(0\right)\right)}\\
        &\qquad\qquad=\left \lVert -d\Phi_n\circ \left(\Psi \circ S\right)+ dx^n\right \rVert_{L^\infty\left(B^-_r\left(0\right)\right)}\\
        &\qquad \qquad\qquad\qquad=\left \lVert -d\Phi_n\circ \Psi + dx^n\right \rVert_{L^\infty\left(B^+_r\left(0\right)\right)}\rightarrow 0 \text{ as } \widetilde{r} \rightarrow 0+.
    \end{align*}
    Recall the definition of $R_\Phi$ in \eqref{RPhi definition} and we have $\lim_{\widetilde{r}\rightarrow 0+}\left \lVert R_{\Phi}- \mathcal{D}\right \rVert_{L^{\infty}\left(B_r^-\left(0\right)\right)}=0$ where $\mathcal{D}$ is a diagonal matrix with diagonal entries $+1$ or $-1$. Moreover $G^-_\Phi=R_\Phi^T R_\Phi$ which implies  $\lim_{\widetilde{r}\rightarrow0+}
     \left \lVert G^-_{\Phi}- I\right \rVert_{L^{\infty}\left(B_r^-\left(0\right)\right)}=0$. Recall the definition of $G_\Phi$ in \eqref{Gphi}
     and we get
     $\left \lVert G_{\Phi}- I\right \rVert_{L^{\infty}\left(B_r\left(0\right)\right)}\text{ as } \widetilde{r} \rightarrow 0^+$. As an immediate consequence, we fix $\widetilde{r}$ small enough so that \eqref{delta}, \eqref{r} and \eqref{1<p<2 and p<2lambda first mention} of $*$-condition hold.
     \end{proof}
Now we are in a position to conclude the proof of Theorem \ref{Goal Theorem}.     
\begin{proof}
We fixed $\widetilde{r}$ small enough such that \eqref{delta}, \eqref{r} and \eqref{1<p<2 and p<2lambda first mention} of $*$-condition holds as guaranteed by Lemma \ref{* condition Lemma} and using Theorem \ref{G Theorem} we have the following:
     \begin{align}\label{last dependence on r}
        &\int_{Q^+_0}\left \lvert d{W}+{F} \right \rvert^{\widetilde{q}}  \leq \begin{aligned}[t]
        c\left(n,\tilde{q},p,k,r\right) \left[\left[\int_{Q^+_0}\left \lvert d{W}+{F} \right \rvert^p\right]^{\frac{\tilde{q}}{p}}
        +\int_{Q^+_0}\left \lvert F\right\rvert^{\tilde{q}}\right].
        \end{aligned}
     \end{align}
Using Lemma \ref{* condition Lemma} we  choose $\widetilde{r}$ small enough such that \eqref{r} and \eqref{1<p<2 and p<2lambda first mention} of $*$-condition holds. So $\widetilde{r}$ depends only on $n,p,\widetilde{q},\alpha_1, c_U$ and $k$, on the other hand $c_U,\alpha_1$ depend only on $n,p,k$. Thus $\widetilde{r}$ depends only on $n,p,\widetilde{q},k$. Also, we fixed $r>0$ small enough so that $B_r^+\left(0\right)\subset \Phi\left(B_{\widetilde{r}}\left(0\right)\cap B_1\left(e^n\right)\right)$; thus $r$ also depends on $\widetilde{r}$, hence on $n,p,\widetilde{q},k$. Therefore the estimate \eqref{last dependence on r} becomes
    \begin{align*}
        \int_{Q^+_0}\left \lvert \Psi^*\left(dv\right)\right \rvert^{\widetilde{q}}  \leq \begin{aligned}[t]
        c\left(n,\tilde{q},p,\alpha_1,c_U,k\right) \left[\left[\int_{Q^+_0}\left \lvert \Psi^*\left(dv\right) \right \rvert^p\right]^{\frac{\tilde{q}}{p}}
        +\int_{Q^+_0}\left \lvert \Psi^*\left(f\right)\right\rvert^{\tilde{q}}\right].
        \end{aligned}
    \end{align*}
Standard change of variable formula gives the following:
        \begin{align*}
        &\int_{\Psi\left(Q^+_0\right)}\left \lvert \left[\Psi^*\left(dv\right)\right]\circ\Phi \right \rvert^{\widetilde{q}}\\  &\quad\leq \begin{aligned}[t]
        &c\left(n,\tilde{q},p,\alpha_1,c_U,k\right) \left[\left[\int_{\Psi\left(Q^+_0\right)}\left \lvert \left[\Psi^*\left(dv\right)\right]\circ \Phi \right \rvert^p\right]^{\frac{\tilde{q}}{p}}
        +\int_{\Psi\left(Q^+_0\right)}\left \lvert \left[\Psi^*\left(f\right)\right]\circ \Phi\right\rvert^{\tilde{q}}\right].
        \end{aligned}
    \end{align*}
      Immediately we get
        \begin{align*}
        \int_{\Psi\left(Q^+_0\right)}\left \lvert dv \right \rvert^{\widetilde{q}}&\leq \begin{aligned}[t]
        &c\left(n,\tilde{q},p,k\right) \left[\left[\int_{\Psi\left(Q^+_0\right)}\left \lvert dv \right \rvert^p\right]^{\frac{\tilde{q}}{p}}
        +\int_{\Psi\left(Q^+_0\right)}\left \lvert f\right\rvert^{\tilde{q}}\right]
        \end{aligned}\\
        &\leq \begin{aligned}[t]
        &c\left(n,\tilde{q},p,k\right) \left[\left[\int_{B}\left \lvert dv \right \rvert^p\right]^{\frac{\tilde{q}}{p}}
        +\int_{B}\left \lvert f\right\rvert^{\tilde{q}}\right].
        \end{aligned}
    \end{align*}
    But $f=du$ in $B$ and by definition $\int_{B} \left \lvert dv\right \rvert^p\leq \int_{B} \left \lvert du\right \rvert^p$, as a consequence
    \begin{align*}
        \int_{\Psi\left(Q^+_0\right)}\left \lvert dv \right \rvert^{\widetilde{q}}\leq 
        c\left(n,\tilde{q},p,\alpha_1,c_U,k\right)\int_{B}\left \lvert dv \right \rvert^{\tilde{q}}
    \end{align*}
        Therefore for $0\in \partial B$ there exists a neighborhood of $0$, $Z_0:=\Psi\left(Q^+_0\right)$ such that  
     \begin{align*}
        \int_{Z_0}\left \lvert dv \right \rvert^{\widetilde{q}}
        &\leq c\left(n,\tilde{q},p,\alpha_1,c_U,k\right)\int_{B}\left \lvert du \right \rvert^{\tilde{q}}.
    \end{align*}
    Recall $\left \lvert Q_0\right \rvert$  depends only on $n,p,\widetilde{q},k$ and $r$ but both $r$ depends on $n,p,\tilde{q},k$. So, $\left \lvert Z_0\right \rvert$ depends only on $n,p,\widetilde{q},k$.
    Similarly for each $x\in \partial B$ there exists a neighborhood of $x$ say $Z_x$ such that  
     \begin{align*}
        \int_{Z_x}\left \lvert dv \right \rvert^{\widetilde{q}}
        &\leq c\left(n,\tilde{q},p,\alpha_1,c_U,k\right)\int_{B}\left \lvert du \right \rvert^{\tilde{q}},
    \end{align*} 
    and by rotational symmetry for $x,y\in \partial B$ we have $\left \lvert Z_x\right \rvert=\left \lvert Z_y\right \rvert$ depending on $n,p,\tilde{q},k,\alpha_1,c_U$.
    But $\partial B$ is compact; hence there are only finitely many $Z_{x_1},\cdots,Z_{x_N}$ forming a cover for $\partial B$ where $N\equiv N\left(n,p,\tilde{q},k\right) \in \mathbb{N}$.\\
    
    Cover $B\setminus \sqcup_{i=1}^N Z_{x_i}$ by balls $B_l\left(y_j\right)$ such that $B_{2l}\left(y_j\right)\subset \subset B$ where $j=1,\cdots,M$; one can notice that both $l$ and $M$ depend on $N$, hence on $n$,$p$,$\tilde{q}$,$k$. Thus for $j=1,2,\cdots, M$
     using \eqref{sup bound} we get the following
     \begin{align*}
         \int_{B_l\left(y_j\right)}\left \lvert dv\right \rvert^{\tilde{q}}\leq  \omega_n l^n \left(\sup_{B_l\left(y_j\right)}  \left \lvert dv\right \rvert^p\right)^{\frac{\tilde{q}}{p}}&\leq \omega_n l^n \left(\bar{c}\fint_{B_{2l}\left(y_j\right)}  \left \lvert dv\right \rvert^p\right)^{\frac{\tilde{q}}{p}}\\&\leq c\left(n,p,\widetilde{q},k\right)l^{n-\frac{n\widetilde{q}}{p}}  \left(\int_{B_{2l}\left(y_j\right)}  \left \lvert dv\right \rvert^p\right)^{\frac{\widetilde{q}}{p}}\\
         &\leq   c\left(n,p,\widetilde{q},k\right)\left(\int_{B}  \left \lvert dv\right \rvert^p\right)^{\frac{\widetilde{q}}{p}}.
     \end{align*}
     For the last estimate we have used the fact that, $l$ depends only on $n,p,\widetilde{q},k,\alpha_1,c_U$ 
     As mentioned earlier $M,N$ depend on $n,\tilde{q},p,\alpha_1,c_U,k$ we finally get the estimate:
    \begin{align*}
        \int_{{B}}\left \lvert dv \right \rvert^{\widetilde{q}}\notag
        &\leq \sum_{i=1}^N\int_{Z_{x_i}}\left \lvert dv \right \rvert^{\widetilde{q}}+\sum_{j=1}^M\int_{{B_{l}\left(y_j\right)}}\left \lvert dv \right \rvert^{\widetilde{\widetilde{q}}}\\
        &\leq \begin{aligned}[t]
        &c\left(n,\tilde{q},p,k\right) \left(M+N\right) \int_{B}\left \lvert du \right \rvert^{\tilde{q}}
        \end{aligned}\notag\leq \begin{aligned}[t]c\left(n,\tilde{q},p,k\right)\int_{B}\left \lvert du \right \rvert^{\tilde{q}}.
        \end{aligned}
    \end{align*}
    \end{proof} 
\textbf{AI Declaration:} The results and proofs in this paper were not generated by AI.\\

\textbf{Acknowledgment:}
The author warmly thanks Dr. Swarnendu Sil for many valuable discussions and suggestions. The research of the author is supported by the NBHM Doctoral Fellowship (Ref. No. 0203/14/2023/R\&D-II/13951). 
\bibliographystyle{plain}
\bibliography{references}

@article {MR3294408,
    AUTHOR = {Colombo, Maria and Mingione, Giuseppe},
     TITLE = {Regularity for double phase variational problems},
   JOURNAL = {Arch. Ration. Mech. Anal.},
  FJOURNAL = {Archive for Rational Mechanics and Analysis},
    VOLUME = {215},
      YEAR = {2015},
    NUMBER = {2},
     PAGES = {443--496},
      ISSN = {0003-9527,1432-0673},
   MRCLASS = {49N60 (35B27 35B65)},
  MRNUMBER = {3294408},
MRREVIEWER = {Eugen\ Viszus},
       DOI = {10.1007/s00205-014-0785-2},
       URL = {https://doi.org/10.1007/s00205-014-0785-2},
}

@article {MR3720379,
    AUTHOR = {Sil, Swarnendu},
     TITLE = {Regularity for elliptic systems of differential forms and
              applications},
   JOURNAL = {Calc. Var. Partial Differential Equations},
  FJOURNAL = {Calculus of Variations and Partial Differential Equations},
    VOLUME = {56},
      YEAR = {2017},
    NUMBER = {6},
     PAGES = {Paper No. 172, 35},
      ISSN = {0944-2669,1432-0835},
   MRCLASS = {35J57 (35B65)},
  MRNUMBER = {3720379},
MRREVIEWER = {Giovanni\ S.\ Alberti},
       DOI = {10.1007/s00526-017-1269-7},
       URL = {https://doi.org/10.1007/s00526-017-1269-7},
}

@article {MR1179331,
    AUTHOR = {Hamburger, C.},
     TITLE = {Regularity of differential forms minimizing degenerate
              elliptic functionals},
   JOURNAL = {J. Reine Angew. Math.},
  FJOURNAL = {Journal f\"ur die Reine und Angewandte Mathematik. [Crelle's
              Journal]},
    VOLUME = {431},
      YEAR = {1992},
     PAGES = {7--64},
      ISSN = {0075-4102,1435-5345},
   MRCLASS = {49N60 (58E15)},
  MRNUMBER = {1179331},
MRREVIEWER = {Jean-Pierre\ Raymond},
       DOI = {10.1515/crll.1992.431.7},
       URL = {https://doi.org/10.1515/crll.1992.431.7},
}

@book{csato2011pullback,
  title={The pullback equation for differential forms},
  author={Csat{\'o}, Gyula and Dacorogna, Bernard and Kneuss, Olivier},
  volume={83},
  year={2011},
  publisher={Springer Science \& Business Media}
}

@book{brezis2010functional,
  title={Functional Analysis, Sobolev Spaces and Partial Differential Equations},
  author={Brezis, H.},
  isbn={9780387709130},
  lccn={2010938382},
  series={Universitext},
  url={https://books.google.co.in/books?id=GAA2XqOIIGoC},
  year={2010},
  publisher={Springer New York}
}

@phdthesis{sil2016calculus,
  title={Calculus of variations for differential forms},
  author={Sil, Swarnendu},
  year={2016},
  school={EPFL}
}

@book{giaquinta2013introduction,
  title={An introduction to the regularity theory for elliptic systems, harmonic maps and minimal graphs},
  author={Giaquinta, Mariano and Martinazzi, Luca},
  year={2013},
  publisher={Springer Science \& Business Media}
}

@article {MR474389,
    AUTHOR = {Uhlenbeck, K.},
     TITLE = {Regularity for a class of non-linear elliptic systems},
   JOURNAL = {Acta Math.},
  FJOURNAL = {Acta Mathematica},
    VOLUME = {138},
      YEAR = {1977},
    NUMBER = {3-4},
     PAGES = {219--240},
      ISSN = {0001-5962,1871-2509},
   MRCLASS = {58E99 (35J60)},
  MRNUMBER = {474389},
MRREVIEWER = {A.\ J.\ Tromba},
       DOI = {10.1007/BF02392316},
       URL = {https://doi.org/10.1007/BF02392316},
}

@article{bwmeta1.element.bwnjournal-article-smv75i1p22bwm,
     author = {Tadeusz Iwaniec},
     title = {Projections onto gradient fields and $L^{p}$-estimates for degenerated elliptic operators},
     journal = {Studia Mathematica},
     volume = {76},
     year = {1983},
     pages = {293-312},
     zbl = {0552.35034},
     language = {en},
     url = {http://dml.mathdoc.fr/item/bwmeta1.element.bwnjournal-article-smv75i1p22bwm}
}

@book{evans2010partial,
  title={Partial Differential Equations},
  author={Evans, L.C.},
  isbn={9780821849743},
  lccn={2009044716},
  series={Graduate studies in mathematics},
  url={https://books.google.co.in/books?id=Xnu0o_EJrCQC},
  year={2010},
  publisher={American Mathematical Society}
}

@article {MR997847,
    AUTHOR = {Acerbi, E. and Fusco, N.},
     TITLE = {Regularity for minimizers of nonquadratic functionals: the
              case {$1<p<2$}},
   JOURNAL = {J. Math. Anal. Appl.},
  FJOURNAL = {Journal of Mathematical Analysis and Applications},
    VOLUME = {140},
      YEAR = {1989},
    NUMBER = {1},
     PAGES = {115--135},
      ISSN = {0022-247X},
   MRCLASS = {49B21 (35J60 58E15)},
  MRNUMBER = {997847},
MRREVIEWER = {T.\ Giec},
       DOI = {10.1016/0022-247X(89)90098-X},
       URL = {https://doi.org/10.1016/0022-247X(89)90098-X},
}

@book{8bbc344b-39d3-31c4-89fe-4b604b1a5849,
 ISBN = {9780691080796},
 URL = {http://www.jstor.org/stable/j.ctt1bpmb07},
 author = {ELIAS M. STEIN},
 publisher = {Princeton University Press},
 title = {Singular Integrals and Differentiability Properties of Functions (PMS-30)},
 urldate = {2026-07-03},
 year = {1970}
}

@article {MR3989958,
    AUTHOR = {Sil, Swarnendu},
     TITLE = {Nonlinear {S}tein theorem for differential forms},
   JOURNAL = {Calc. Var. Partial Differential Equations},
  FJOURNAL = {Calculus of Variations and Partial Differential Equations},
    VOLUME = {58},
      YEAR = {2019},
    NUMBER = {4},
     PAGES = {Paper No. 154, 32},
      ISSN = {0944-2669,1432-0835},
   MRCLASS = {58J05 (35B65 35J47)},
  MRNUMBER = {3989958},
MRREVIEWER = {Matteo\ Cozzi},
       DOI = {10.1007/s00526-019-1598-9},
       URL = {https://doi.org/10.1007/s00526-019-1598-9},
}

@book{dacorogna2007direct,
  title={Direct Methods in the Calculus of Variations},
  author={Dacorogna, B.},
  isbn={9780387552491},
  lccn={92009560},
  series={Applied Mathematical Sciences},
  url={https://books.google.co.ls/books?id=79no0AEACAAJ},
  year={2007},
  publisher={Springer New York}
}
\end{document}